\documentclass{article}
\usepackage[utf8]{inputenc}
\usepackage[a4paper]{geometry}
\usepackage{amsthm}
\usepackage{amsfonts}
\usepackage{epsfig}
\usepackage{color}
\newcommand{\bibd}{bib/}
\newcommand{\eeq}{\end{equation}}
\newcommand{\bee}{\begin{eqnarray}}
\newcommand{\eee}{\end{eqnarray}}
\newcommand{\been}{\begin{eqnarray*}}
\newcommand{\eeen}{\end{eqnarray*}}

\newcommand{\qq}[1]{\footnote{ \textcolor{red}{#1} }}  
\newcommand{\qqqq}[1]{\footnote{ \textcolor{green}{#1} }}  
\newcommand{\qqq}[1]{\footnote{ \textcolor{blue}{#1} }}  
\providecommand{\noglossaryignore}[1]{}
\newcommand{\globalglossaryentry}[3]{\makebox[1.5in][l]{\tt $\backslash${#1}} 
\makebox[1.1in][l]{{$#2$}} \makebox[2.5in][l]{{#3}}\newline} 
\newcommand{\newcommandabbreviation}[3]{\newcommand{#1}{#2}%
\noglossaryignore{\globalglossaryentry{#3}{#2}{}}}
\newcommand{\renewcommandabbreviation}[3]{\renewcommand{#1}{#2}%
\noglossaryignore{\globalglossaryentry{#3}{#2}{}}}
\newcommand{\newcommandmacro}[4]{\newcommand{#1}{#2}%
\noglossaryignore{\globalglossaryentry{#3}{#2}{#4}}}
\newcommand{\gge}[3]{\noglossaryignore{\globalglossaryentry{#1}{#2}{#3}}}
\newcommand{\myaddress}%
{\parbox{3in}{\footnotesize \begin{center} 
Mathematics Department, City University, \\  
Northampton Square, London EC1V 0HB, UK.\end{center}}}
    
\newcounter{minidef}[section]
\renewcommand{\theminidef}{\thesection.\arabic{minidef}}
\newcommand{\mdef}{\refstepcounter{minidef} 
\medskip \noindent ({\bf \theminidef}) }

\newtheoremstyle{pu}
{7pt}%
{7pt}%
{\it}
{}
{}
{.}
{ }
{\thmnumber{({\bf #2}) }\thmname{\textsc{#1}}\thmnote{#3}}
\newtheoremstyle{puu}
{3pt}%
{3pt}%
{\rm}
{}
{}
{.}
{ }
{\thmnumber{({\bf #2}) }\thmname{\textsc{#1}}\thmnote{#3}}

\theoremstyle{pu}
\newtheorem{mmpr}[minidef]{Proposition}

\newtheorem{mmco}[minidef]{Corollary}

\theoremstyle{puu}

\newcounter{minicapt}

\theoremstyle{pu}
\newtheorem{theo}[minidef]{Theorem}      
\newtheorem{cor}{Corollary}[minidef]   
\newtheorem{pr}[minidef]{Proposition} 

\newtheorem{lem}[minidef]{Lemma} 

\theoremstyle{puu}
\newtheorem{de}[minidef]{Definition}

\noglossaryignore{GREEK ETC.\newline}
\newcommandabbreviation{\e}{\epsilon}{e}        
\newcommandabbreviation{\lam}{\lambda}{lam}  
\newcommandabbreviation{\la}{\langle}{la}        
\newcommandabbreviation{\ran}{\rangle}{ran}
\newcommandabbreviation{\ha}{\#}{ha}             
\newcommandabbreviation{\rmap}{\rightarrow}{rmap}
\newcommandabbreviation{\aaa}{\alpha}{aaa}        
\newcommandabbreviation{\ab}{\alpha,\beta}{ab}
\newcommandabbreviation{\aab}{a(\ab )}{aab}       
\noglossaryignore{\newline RINGS\newline}
\newcommandabbreviation{\HH}{H \!\!\! I}{HH}               
\newcommandabbreviation{\C}{\mathbb C}{C}
\newcommandabbreviation{\N}{\mathbb N}{N}   
\newcommandabbreviation{\Z}{\mathbb Z}{Z}      
\renewcommandabbreviation{\Re}{\mathbb R}{Re}
\newcommandabbreviation{\R}{{\mathbb R}}{R}
\newcommandabbreviation{\Q}{\mathbb Q }{Q}
\renewcommandabbreviation{\H}{\mathbb H }{H}
\noglossaryignore{\newline SYMMETRIC GROUP\newline}
\def\Sym(#1){\Sigma(#1)}                   
\gge{Sym(-)}{\Sym(-)}{symmetric group on - objects}
\def\Sy(#1){\Sigma_{#1}}                   
\gge{Sy(-)}{\Sy(-)}{symmetric group irreducible -}
\def\sym(#1){\mbox{\LARGE s}(#1)}        
\gge{sym(-)}{\sym(-)}{symmetric group on - objects (variant)}
\def\sy(#1){\mbox{\LARGE s}({#1})}        
\gge{sy(-)}{\sy(-)}{symmetric group irreducible - (variant)}
\newcommandmacro{\cs}{\C \, \sy(n)}{cs}{symmetric group algebra over $\C$}
\noglossaryignore{\newline PARTITIONS/SETS\newline}
\newcommand{\Nset}[1]{\underline{#1}}
\gge{Nset\{-\}}{\Nset{-}}{set of natural numbers to -} 
\def\nset(#1){ \{ #1 \}_{ \underline{n} }} 
\gge{nset(-)}{\nset(-)}{a set $-\times\Nset$}
\def\ul(#1){_{\underline{#1}}}             
\gge{ul(-)}{{}\ul(-)}{subscript underline -}
\def\Ee(#1){{\bf E}_{#1}}                  
\gge{Ee(-)}{\Ee(-)}{set of equivalence relations on set -}
\def\Eee(#1){{\bf E}_{\{ #1 \}_{\underline{n}}}}   
\gge{Eee(-)}{\Eee(-)}{ditto for nset}
\def\Een(#1,#2){{\bf E}_{\{ #1 \}_{\underline{#2}}}}   
\def\Ssn(#1,#2){{\bf S}_{\{ #1 \}_{\underline{#2}}}}   
\def\Ss(#1){{\bf S}_{#1}}                  
\def\Sss(#1){{\bf S}_{\{ #1 \}_{\underline{n}}}}   
\def\bbc(#1){((\beta_1)(\beta_2)...(\beta_{#1}))}      
\newcommandmacro{\Ln}{{\Gamma}^{n}}{Ln}{large index set}
\newcommandmacro{\LnQ}{{\Gamma}^{n}_Q}{LnQ}{index set}
\newcommandmacro{\Zz}{\zeta}{Zz}{`shape' function}
\noglossaryignore{\newline PARTITION ALGEBRA\newline}
\def\ka(#1){\kappa_{#1}}                   
\def\Sm(#1){\Sigma_{#1}}                   
\newcommandmacro{\com}{\bullet}{com}{bullet composition}
\newcommandmacro{\enm}{\; e^n(\! m\! ) \;}{enm}{product of idempotents}
\def\Ai(#1){ A^{ #1 \cdot } }              
\def\Aij(#1,#2){ A^{ #1  #2 } }            
\newcommandmacro{\One}{\mbox{\bf $1 \!\!\! 1$}}{One}{algebra unit 1}
\newcommandmacro{\Bp}{B_p}{Bp}{partition basis}
\def\Bb(#1){B_p[#1]}                       
\def\Pp(#1){P_n[#1]}                       
\def\Ps(#1){P_n[#1] \! /}                  
\newcommandmacro{\Ph}{\hat{P}}{Ph}{P hat  algebra}
\def\Is(#1){\sim^{#1}}                     
\noglossaryignore{\newline MODULES\newline}
\def\Wm(#1){{\cal S}_{#1}}                 
\gge{Wm(-)}{\Wm(-)}{Weyl module with index -}
\def\wm(#1,#2){{}_{#1}{\cal S}_{#2}}       
\gge{wm(-1,-)}{\wm(-1,-)}{Weyl module with index -}
\def\Ind(#1,#2,#3){\mbox{Ind}_{#1}^{#2}#3} 
\gge{Ind(-1,-2,-)}{\Ind(-1,-2,-)}{induction}
\def\Res(#1,#2,#3){\mbox{Res}_{#1}^{#2}#3} 
\gge{Res(-1,-2,-)}{\Res(-1,-2,-)}{restriction}
\newcommandabbreviation{\weyl}{standard}{weyl}
\newcommandabbreviation{\mod}{\mbox{mod}}{mod}
\newcommandabbreviation{\head}{\mbox{head }}{head}
\newcommandabbreviation{\Weyl}{Weyl}{Weyl}
\def\SS(#1){{\cal S}_{#1}}                 
\gge{SS(-)}{\SS(-)}{Specht/Weyl module index -}
\def\LL(#1){{\cal L}_{#1}}                 
\gge{LL(-)}{\LL(-)}{simple module index -}
\noglossaryignore{\newline FUNCTORS/MAPS\newline}
\newcommandmacro{\Gg}{{\cal G}}{Gg}{G Functor}
\newcommandmacro{\Fg}{{\cal F}}{Fg}{F Functor}
\newcommandmacro{\ra}{\rightarrow}{ra}{}
\def\ses(#1,#2,#3){0\ra #1 \ra #2 \ra #3 \ra 0}   
\gge{ses(1,2,-)}{\ses(1,2,-)}{\hspace{.5in} short exact sequence}
\def\starr(#1){ \stackrel{ #1 }{\longrightarrow} }
\gge{starr(-)}{\starr(-)}{}
\newcommandmacro{\doublerightarrow}{\; -\!\!\! -\!\!\!\!\!\! \gg \;}
{doublerightarrow}{}
\noglossaryignore{\newline PARTITION ALGEBRA MAPS\newline}
\newcommandmacro{\smap}{s}{smap}{`inclusion' map}
\newcommandmacro{\tmap}{t}{tmap}{$ P_n -> S_n$}
\newcommandmacro{\pmap}{\psi}{pmap}{$ S_n -> P_n $}
\noglossaryignore{\newline MISC.\newline}
\def\Amap(#1){{\cal A}_{#1}}               
\gge{Amap(-)}{\Amap(-)}{}
\def\Rr(#1){R_{#1}}                        
\gge{Rr(-)}{\Rr(-)}{restriction of E}
\def\Cr(#1){C_{#1}}                        
\gge{Cr(-)}{\Cr(-)}{restriction of E to N}
\newcommandmacro{\Tm}{{\cal T}}{Tm}{Transfer Matrix}
\def\On(#1){{\cal I}_{#1}}
\gge{On(-)}{\On(-)}{}
\newcommandmacro{\UU}{\underline{\sqcup}}{UU}{}  
\newcommandmacro{\UUU}{\sqcup}{UUU}{}  
\newcommandmacro{\Vq}{V_Q^{\otimes n}}{Vq}{Potts config. space}
\def\bs(#1,#2){\mbox{{\Large $\ast$}}^{#1}_{#2}}  
\gge{bs(-,-)}{\bs(-,-)}{general plumbing multiplier}
\newcommand{\ignore}[1]{}
\gge{ignore\{-\}}{\ignore{-}}{ignore argument!}
\def\choo(#1,#2){ \left( \begin{array}{c} #1 \\ #2 \end{array} \right) } 
\gge{choo(-1,-)}{\choo(-1,-)}{choose}
\newcommand{\Qed}{$\Box$}
\gge{Qed}{\mbox{\Qed}}{QED}
\def\staq(#1){\stackrel{#1}{=}}            
\gge{staq(-)}{\staq(-)}{}
\def\stam(#1){\stackrel{#1}{\rightarrow}}  
\gge{stam(-)}{\stam(-)}{}
\def\mat{ \left( \begin{array} }    
\def\tam{ \end{array}  \right) }
\gge{mat/tam}{...}{matrix delimiters}
\newcommand{\beq}{\begin{equation} }
\def\eql(#1){ \begin{equation} \label{#1} 
}
\newcommand{\eq}{\end{equation} }
\def\eqal(#1){\begin{eqnarray} \label{#1} }
\def\eqa{\end{eqnarray} }
\def\lab(#1){\label{#1}
}
\def\prl(#1){ \begin{pr} \label{#1} 
}
\def\del(#1){ \begin{de} \label{#1} 
}
\gge{smeq\{-\}}{...}{small equation}
\gge{fneq\{-\}}{...}{very small equation}
\noglossaryignore{\newline HECKE/BLOB\newline}
\newcommandmacro{\Hnq}{H_n(q)}{Hnq}{ * freestanding symbol}
\newcommandmacro{\Hn}{H_n}{Hn}{      *-mod etc.}
\newcommandmacro{\A}{{\cal A}}{A}{}
\newcommandmacro{\Cwts}{C}{Cwts}{}
\newcommandmacro{\CA}{{\cal A}}{CA}{}

\newcommandmacro{\calA}{{\cal A}}{calA}{}
\newcommandmacro{\modi}{\mbox{Mod} }{modi}{was mod not modi!}
\newcommandmacro{\Wgen}{{\Bbb S}}{Wgen}{}
\def\ol(#1){\overline{#1}}
\newcommandmacro{\st}{\mbox{St}}{st}{}
\def\CMult(#1,#2){(#1:#2)}
\def\CM(#1,#2){( #1 : #2 )}
\def\FMult#1,#2{(#1:#2)}
\def\CF#1,#2{(#1:#2)}

\newcommandmacro{\Top}{\mbox{Top}}{Top}{}
\newcommandmacro{\Soc}{\mbox{Soc}}{Soc}{}
\newcommandmacro{\Head}{\mbox{Head}}{Head}{}
\newcommandmacro{\Filt}{{\cal F}}{Filt}{}
\newcommandmacro{\Mod}{\mbox{mod}}{Mod}{}
\newcommandmacro{\Resi}{\mbox{Res }}{Resi}{was without i!}
\newcommandmacro{\Indi}{\mbox{Ind }}{Indi}{was without i!}
\def\RR(#1,#2){R^{#1}_{#2}}   
\def\TT(#1,#2){T^{#1}_{#2}}   

\newcommandmacro{\Ann}{\mbox{Ann}}{Ann}{}
\newcommandmacro{\Cen}{\mbox{Cen}}{Cen}{}
\newcommandmacro{\End}{\mbox{End}}{End}{}
\newcommandabbreviation{\semisimple}{semisimple}{semisimple}
\newcommandabbreviation{\Bratteli}{Bratteli}{Bratteli}
\newcommandabbreviation{\JBC}{Jones Basic Construction}{JBC}
\newcommandabbreviation{\pa}{partition algebra}{pa}
\newcommandabbreviation{\TM}{transfer matrix}{TM}
\newcommandabbreviation{\PM}{Potts model}{PM}
\newcommandabbreviation{\QSC}{quantum spin chain}{QSC}
\newcommandabbreviation{\Hamiltonian}{Hamiltonian}{Hamiltonian}
\newcommandabbreviation{\YS}{Young symmetrizer}{YS}

\newcommand{\BB}{{\mathcal B }}   
\newcommand{\JJJ}{{\mathcal J}}  
\newcommand{\BBl}[1]{{\mathcal J}_{#1}}   
\newcommand{\BBBl}[1]{{{\mathcal J}_{#1}^1}}  
\newcommand{\JJ}{J}  
\newcommand{\TTT}{{\mathcal T }}   
\newcommand{\bb}{{\mathfrak b}}   
\newcommand{\bB}{{{\mathsf J}^{\bullet}}}  
\newcommand{\deltap}{\delta'}    

\newcommand{\bBd}{\bB}  
\newcommand{\Jd}{{J}}  

\newcommand{\paref}[1]{(\ref{#1})} 

\newcommand{\ppm}[1]{\textcolor{red}{#1}}
\newcommand{\red}[1]{\textcolor{red}{#1}}
\newcommand{\redx}[1]{}

\newcommand{\kmy}{\cite{\kamy}} 
\newcommand{\kamy}{KadarMartinYu}
\newcommand{\ms}{MartinSaleur94a}

\begin{document}
\title{
Geometric partition categories: \\
  On short Brauer algebras and their blob subalgebras}
\author{Z. K\'AD\'AR and P. P. MARTIN
\\ School of Mathematics, University of Leeds}
\date{}
\maketitle

\noindent {\small
Abstract:
The main result here  
gives an algebra(/linear category) 
isomorphism between a
geometrically defined subcategory
$\JJJ^1_0$
of a
short Brauer category $\JJJ_0$ and a
certain one-parameter specialisation of the blob category $\bb$.
That is,  
we prove
the Conjecture in Remark 6.7 of \cite{\kamy}.
We also define a sequence of generalisations $\JJJ^i_{i-1}$ of the
category $\JJJ^1_0$.
The connection of $ \JJJ_0$ with the blob category inspires a search for
connections also with its beautiful representation theory.
Here we obtain formulae determining the non-semisimplicity condition
(generalising the classical `root-of-unity' condition).
}

\medskip

\noindent {\small
{\em Keywords}: diagram algebra, topological spin chain.}

\medskip

\section{Introduction}

A motivating aim here is to study the structure  
of the $k$-linear categories $\BBl{l}$ from \cite{\kamy},
and in particular the representation theory of the 
corresponding $k$-algebras
(with $k$ a field)
$\JJJ_{l,n}$
in the non-semisimple cases. 
These structures are of intrinsic interest
(cf. \cite{jk,james,Martin0915,CoxDeVisscher});
and see also \cite{\kamy} for a discussion of some
of the extrinsic motivations for this study 
--- in short one seeks generalisations 
of the intriguing examples of Kazhdan--Lusztig theory
\cite{KazhdanLusztig79,Soergel97a,AndersenJantzenSoergel94} 
observed  \cite{Martin0915}  in the representation theory
of the Brauer category $\BB = \BBl{\infty}$ \cite{Brauer37}.
Another motivating aim is to study module categories over monoidal
categories
(see e.g. \cite{Ostrik01} for a review)
beyond the usual `semisimple' setting.

The study strategy
in Part~1 (\S\ref{ss:pre1}-\ref{ss:main01})
can be seen as trying 
to relate the problem to the representation theory of
the blob category $\bb$ and the blob algebra $\bb_n$ \cite{\ms}, 
which is contrastingly very well understood (see e.g. \cite{CoxGrahamMartin03}),
itself with deep and 
tantalising connections to Kazhdan--Lusztig theory
\cite{MartinWoodcock03}.
(More recently see e.g. \cite{BowmanCoxSpeyer}.)
This also allows us to make contact with the original physical
motivations for these algebras, as the algebras of physical systems
with boundaries and interfaces \cite{\ms}.
Indeed the blob algebras have been of renewed interest recently in
several areas, not only of physics but also for example the study of
KLR algebras \cite{KLRI,KLRII},
Soergel bimodules \cite{Soergel07} and monoidal categories \cite{JoyalStreet}.

As we shall see, in the simplest non-trivial case the algebras are
(at least) related by inclusions of the form 
$\bb_m \hookrightarrow \JJJ_{0,n}$. 
Inclusion is not in general a directly helpful
relationship in representation theory. 
(For example the Temperley--Lieb algebra $T_n$ \cite{tl} is a subalgebra of
$\bb_n$, but the representation theories of these algebras are radically
different: cf. \cite{Martin91} and \cite{CoxGrahamMartin03}.) 
However the inclusion here is of
`high index', so there is hope that it will indeed shed light on the
open problem.

In Part~2 (\S\ref{ss:other1}) we include some indicative results on
$\JJJ_{0,n}$ representation theory.
These are obtained by working directly with $\JJJ_{0,n}$,
but serve 
as a first step in this direction
(full analysis of these results is demoted to a separate paper).

\medskip

In Section~\ref{ss:pre1} we introduce concepts and notations.
In  
\S\ref{ss:s4} we 
define for each category $\BBl{l}$ a new subcategory.  
In \S\ref{ss:main01} we  examine the relationship to the blob
category. 
In particular in Section~\ref{ss:main01} we state and prove the main theorem.
In section~\ref{ss:other1} we consider consequences for the algebras
$\BB^l_n  =\JJJ_{l,n}$ themselves. 
In Section~\ref{ss:discusstar} 
we discuss  
related open problems.

\vspace{.51cm} 

\section{Preliminary definitions} \label{ss:pre1}


\newcommand{\Vnm}{V^n_m}  

Define $\underline{n} = \{ 1,2,...,n \}$,
$\underline{n}' = \{ 1',2',...,n' \}$
and so on. 
Define $\Vnm = \underline{n}\cup\underline{m}'$.
Write $J(n,m)$ for the set of set partitions of
$\underline{n}\cup\underline{m}'$
into subsets of order 2.
%
Fix a commutative ring $k$ and $\delta\in k$. 
In \S\ref{ss:bb1}-\ref{ss:bb12}   we
recall 
the  Brauer partition category 
$$
{\BB} = ({\mathbb N}_0 , k J(n,m), *)
$$
with loop-parameter $\delta$.
%
%
The category $\BB$ 
has an infinite family of subcategories $\BBl{l}$ introduced in \cite{\kamy}. 
The key  
ingredient in \cite{\kamy}
is the definition of the {\em left-height} of a partition.
We recall it here in \S\ref{ss:short1}.
%
In \S\ref{ss:blob1}, 
we recall the blob category \cite{\ms}.


\subsection{
  Brauer picture calculus} \label{ss:bb1}

%
\begin{figure}
\begin{center}\includegraphics[width=4cm]{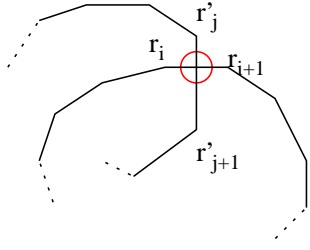}\end{center}
\vspace{-.21in}
\caption{\label{inter}
  A regular intersection between lines
(circled).
}
\end{figure}
%

We 
recall the definition of multiplication $*$ in ${\BB}$
and relate to
`pictures' of partitions. 




\mdef \label{de:line}
Fix a dimension $d>0$.
Given an ordered list $r = (r^1, r^2 , ..., r^l)$ of points in $\R^d$,
let
\[
(r) = \bigcup_{i=1}^{l-1} [r^i , r^{i+1} ]
\]
where $[r^i , r^{i+1} ] \subset \R^d$ is the straight line between these points.

A {\em line}  or {\em polygonal path}
in the plane $\R^2$  is  a
subset
of form $(r)$ 
such that each point $r^i$ lies in at most two of the straight lines.  
%

\mdef Let $R \subset \R^2$ be a rectangle, with frame $\partial R$ and
interior $(R)= R \setminus \partial R$. 
A set of lines is {\em regular in $R$} if:
\\
(R1) each line $(r) \subset R$ touches  $\partial R$ only if $r^1 \neq r^l$
and then only at its end-points $r^1, r^l$;
\\
(R2)
the point-list
$r^1, r^2 , ..., r^l$
of each line has no intersection with any other line.
See Fig.~\ref{inter}.

\mdef Remark. \label{re:fib}
Given  a regular set of lines $D$,
consider the subset $p(D)$ of $R$ defined by  $D$.
Note that we can recover the points of $D$ from $d=p(D)$ except for those
points $r^i$ that are colinear with their neighbours $r^{i-1},r^{i+1}$.
Note that colinearity is not a generic condition.
The fibre $p^{-1}(d)$ of regular sets of lines over $d$
includes a representative $D'$ in which no line has an $i$ with
$r^{i-1},r^i,r^{i+1}$ colinear.
The fibre 
consists of sets obtained by inserting such colinear points in lines.

Note that $d$ completely determines $D$ up to such inserted colinear points.

\begin{figure}
\begin{center}
\includegraphics[width=3.4cm]{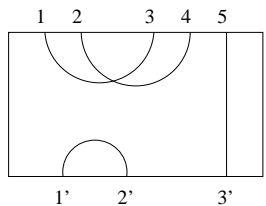}
\end{center}
\vspace{-.21in}
\caption{\label{fig:nc1}
A picture of $\{ \{ 3,1 \}, \{ 5,3' \}, \{ 4,2 \},\{ 2',1' \}\}\in J(5,3)$. 
N.B. Drawings of piecewise-smooth approximations to piecewise-linear
embeddings are safe to use here, provided that crossings remain
manifestly transversal.
}
\end{figure}

\mdef 
\label{de:picture} 
A {\em Brauer picture} $d$ of a partition $p$ in $\JJ(n,m)$ is 
a rectangle $R\subset \R^2$ with 
$n$ points
(called `vertices')
labeled $1,2,...,n$
on the northern edge,
and $m$ on the southern edge
(as in Figure~\ref{fig:nc1}); 
and a subset of $R$ consisting of  
a regular set of lines $(r)$ in $R$
as follows. 
\ignore{{
where a line is an ordered list of finitely many points 
$r^1, r^2 , ..., r^l$ in $R$ and the straight lines between each pair
$r^i , r^{i+1}$, such that $r^i$ lies in at most two straight lines.
}}
Each line is either a loop ($r^1 = r^l$) or else connects vertices
pairwise in accordance with the pairs in $p$.
\ignore{{
the set of lines  is {\em regular}
(cf. \cite{kmy})%
:
\\ (R1) a line $(r)$ touches the frame $\partial R$ at most at its
end-points $r^1, r^l$;
\\ (R2) 
the $r$-list of each line has no intersection with any other line.
See Fig.~\ref{inter}.
}}


\begin{figure}
\begin{center}\includegraphics[width=9.7cm]{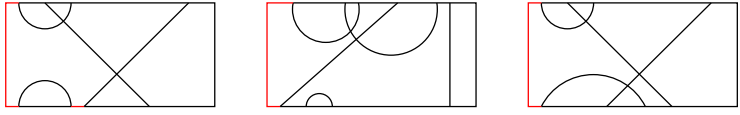}\end{center}
\vspace{-.12in}
\caption{\label{fig:ls1} 
 `Lowest-height' pictures of elements of $J(4,4)$.}
\end{figure}

See Fig.\ref{fig:nc1} 
and Fig.\ref{fig:ls1} for  examples. 
%
Remark:
As far as physically drawn figures are concerned,
piecewise linear and piecewise smooth lines are effectively
indistinguishable,
by \cite[§6]{Moise77} for example.

\mdef \label{de:pi}
Note that 
the construction ensures a well-defined map $\pi$ from pictures to
partitions: 
for each $i \in \Vnm$ the 
 line from $i$ in $d$ may be followed unambiguously  to the other end;
and one takes a line with end-points $i,j$ to $\{ i,j \}$.


\mdef
For $p \in J(n,m)$ write $[p]$ for the set of  pictures $d$ such that
$\pi(d)=p$.
Note that every $[p]$ is non-empty.
Write
$\#(d)$ for the number of loops in  picture $d$. 

%

\mdef
Recall
$\delta \in k$. 
We extend the $\pi$-map to  
$\Pi(d)=\delta^{\#(d)} \pi({d}) \in k J(m,n)$. 

For $d$ a picture let $\hat{d}$ denote  $d$   with all loops removed. 
Thus $\Pi(d) = \delta^{\#(d)} \Pi(\hat{d})$. 


\medskip

\mdef
Consider pictures
$d_1$ for $p_1\in J(m,n)$ and $d_2$ for $p_2\in J(n,q)$.  


We have not specified exactly where the $n$ vertices lie on the
southern frame of $R$ in $d_1$
(or indeed where the southern frame lies in $\R^2$),
but it will be clear that there are representative pictures
of $p_1, p_2$ for which 
the $n$ vertices from $d_1$ match up with the $n$ from $d_2$.
These pictures can then be stacked  (with $d_1$ over $d_2$)
so that the $n$ vertex sets in each
`factor' coincide.
Note that the concatenation is again a Brauer picture. 
It is denoted by $d_1\,|\,d_2$. 

Note then (1) that 
$d_1 | d_2$ can be seen as a picture of an element of $k J(m,q)$;
\\
(2) that if $d_2$ and $d_3$ also stack then 
\beq \label{eq:assocx}
(d_1 | d_2)|d_3 = d_1|(d_2 | d_3)
\eq


\subsection{The Brauer partition category $\BB$} \label{ss:bb12}

Given a set $p$ of symbols, let $p'$ be the set obtained by adding
primes to symbols in $p$
(note that this can work recursively).
For example
$ \{ \{ 1,2 \},\{3,2' \}, \{1',3' \}\}'
= \{\{ 1',2' \},\{ 3',2''\},\{ 1'',3''\}\}$.

\mdef \label{de:star}
Let $p_1 \in J(m,n)$ and $p_2 \in J(n,q)$.
Now 
form $p_1 \cup p_2'$.
Note that $p_1$ and $p_2'$ may not be quite disjoint in general.
When a primed pair in $p_1$ meets the corresponding unprimed pair
from $p_2$ the union  `flattens' this to  a single pair. 
Note each single-primed element appears twice
(or once if flattened)
and others once.
Consider a maximal chain
$\{ a_0, a_1 \},\{ a_1, a_2 \},...,\{a_{k-1}, a_k \}$
in $p_1 \cup p_2'$.
The chain is either all primed,
and $a_0 = a_k$ or $k=1$;
or else 
$\{ a_0, a_k \} $ lies in 
$\underline{m}\cup\underline{q}''$
(to which we may apply $i'' \leadsto i'$ to obtain an element of
$\underline{m}\cup\underline{q}'$).
In this way $p_1 \cup p_2'$ determines an element  
$p_1 . p_2$ of $J(m,q)$; and
also
a number $\#$ of all-primed `closed' chains.
Define
$$
p_1 * p_2 = \; \delta^\# p_1 . p_2  \;  \in kJ(m,q) .
$$
Note  that if $d_i \in [p_i]$ and $d_1 | d_2$ then
\\
(I)
each maximal chain above corresponds to  vertices in a line 
component of $d_1 | d_2$.  
\\
(II)
the pair $\{ a_0, a_k \} $
corresponds to the vertices at the ends of the 
line component
of $d_1 | d_2$.  
Thus the image in $\pi(d_1 | d_2)$
if $a_0 \neq a_k$ is $\{a_0, a_k \}$
(more precisely $\{ a_0, a_k \}$ with any $i'' \leadsto i'$).



\mdef \label{de:stary}
From the 
notes (I,II) in (\ref{de:star}) we have,
independently of the choice of pictures,
$$
p_1 * p_2 = \Pi(d_1\,|\,d_2)   .  
$$ 
For example, from Figure~\ref{fig:wid1} (ignoring the red line for now)
\[
\{\{ 4',2\},\{3,5'\},\{1,3'\},\{ 1',2' \}\} *
 \{\{2,1\},\{4,5\},\{ 1',3 \}\}
\; = \; \delta
\{\{ 3,2 \} , \{ 1,1'\}\} \; \in kJ(3,1)\ .
\]
Note we use the convention that a picture $d_1$ for $p_1\in J(m,n)$
has $m$ (resp. $n$) points on its northern  (southern) edge. 
Thus $d_1\,|\,d_2$ is the concatenation of $d_1$ on top of $d_2$.  


{\theo{\label{th:B} {\rm \kmy}
Composition $*$ is 
associative. 
This defines the category $\BB$.
}}
\proof
By the 
points (I,II) 
 in (\ref{de:star})
above we may obtain $p_1 * p_2$ from
$d_1 | d_2$ (independently of the choices of these pictures).
Existence of constructs of form $d_1 | (d_2 | d_3) $ will  be evident.
Associativity follows since $(d_1 | d_2)|d_3 = d_1|(d_2 | d_3)$. 
(Alternatively we may stay with the initial
machinery of (\ref{de:star}) and simply
introduce more primes. The pictures can be seen as bookkeeping the primes.)
\qed

\mdef \label{de:moncat}
Define
$\otimes : J(n,m) \times J(n',m')
\rightarrow J(n+n', m+m')$
as the composition corresponding to side-by-side concatenation of
pictures. Note the following. 

\begin{lem} \label{lem:moncancel}
  (I) The composition $\otimes$ makes $\BB$ into a monoidal category.
  \\
  (II) 
  Let $w \in J(n,m)$ for some $n,m$,
  with $J(n,m) \hookrightarrow kJ(n,m)$ in the natural way.
  Then $p\otimes w = q\otimes w$
  implies $p=q$. 
  \qed
\end{lem}


\subsection{The subcategory $\BBl{l}$ of category $\BB$} \label{ss:short1}

The regularity property 
(R2) of a picture $d$ means that the number, 
and position, of crossings of lines in $d$ is well-defined. 
A {\em path} through a picture $d$ is a 
further line in $R$ that
satisfies (R2).

Fix a picture $d$, let $x$ be a point in $d$,
and consider all paths from $x$ to the left edge. 
Then 
`left-height' 
$h_d(x)$  
is the minimum number of crossings with lines of the picture
among such paths. 

The {\em left-height} $h(d)$ 
of a picture $d$ is the maximum $h_d(x)$ of all crossing
points $x$ of lines in $d$; or, if there are no crossings,
then $h(d)=-1$. 

For $p \in J(n,m)$ let $h(p)$ denote
the minimum $h(d)$ among pictures in $[p] = \{d: \pi(d)=p\}$. 
\ignore{{
Here we will omit ``left-'' and just use  ``height''. 
We denote the height function on partitions by 
$h:J(m,n)\to \{-1,0,1,2\dots\}$. 
The value $-1$ is
a convention for pair partitions, which have a  
picture without crossing points, these are arrows of the
Temperley-Lieb subcategory ${\cal T}$.
}}
Define
\[
\JJ_{\leq l}(n,m) = \{ p \in \JJ(n,m) \; | \; h(p) \leq l \}
\] 

Define $[p]'$ as the subset of $[p]$ of pictures of $p$ of the
minimum height. That is, 
$$
[p]' = \{ d \in [p] \; | \; h(d)=h(p) \} .
$$
See Fig.\ref{fig:ls1} for examples of minimum height pictures.


As shown in \kmy,
$h(p_1*p_2)\leq \max(h(p_1),(h(p_2))$. 
A consequence is the following.
{\theo{{\rm \kmy} 
Category ${\BB}$ has a subcategory 
$\BBl{l} = ({\mathbb N}_0 , k \JJ_{\leq l}(n,m), *)$.
\qed
}}


\mdef 
\label{de:alcy}
Given a subset $S$ of a rectangle $R$, an {\em alcove} of $S$ is a connected
component of $R\setminus S$.


\mdef Remark.  \label{re:bhi}
Let $d \in [p]$.
Note that alcoves of $d$ have well-defined
left-height.
Note that the
left-heights of the intervals of the frame of $R$ 
are determined by $p$, and are otherwise independent of $d$.

Proof:
Note that there exist  paths $w$ from the left edge
of $R$ to points on $\partial R$
such that $w$  lies 
in  a neighbourhood of $\partial R$.
Note
(e.g. from the Jordan Curve Theorem)
that there are such paths that  have the lowest number of crossings.
By construction two pictures $d,d' \in [p]$ are close to
identical in a neighbourhood of $\partial R$.
In particular
there are paths to points on $\partial R$ that lie 
in such a neighbourhood; have the lowest number of crossings;
and that have the same
number of crossings in $d$ and $d'$.
\qed



\subsection{The blob category $\bb$} \label{ss:blob1}
\newcommand{\bx}{\circ}  

\begin{figure}
\begin{center}\includegraphics[width=7.9cm]{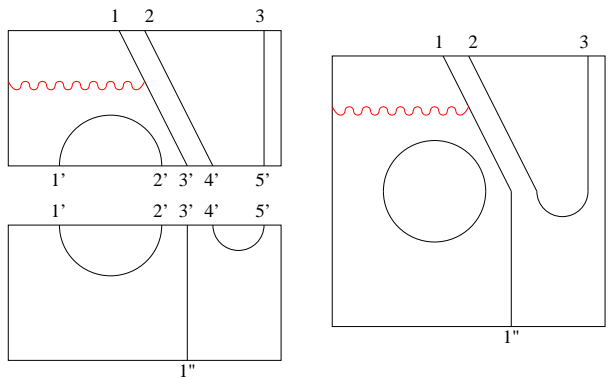}\end{center}
\caption{ A left-exposed line is a part of
  a left-exposed line after concatenation.
  \label{fig:wid1}}
\end{figure}


\begin{figure}
  \begin{center}
    \includegraphics[width=5cm]{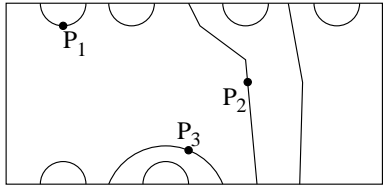}
  \end{center}
\caption{ \label{fig:many2} \label{fig:which1}   
   An illustrative blob picture.
}
\end{figure}

Now we  recall the blob category \cite{\ms}.
Note that for  $p\in J_{-1}(m,n)$ a pictures $d$  in $[p]'$
has no crossings. 
Thus for each pair $v$ in $p$ the corresponding
line $l_v$ in $d \in [p]'$ has the property that every point
$x$ on $l_v$ has the same value of  
$h_d(x)$.
Furthermore this $h_d(l_v) = h_d(x)$ depends only on
$v$ in $p$ and not on $d$.
Thus  
$v$ has a well-defined $h(v) =h_d(l_v)$.
A {\em left-exposed} pair in $p\in J_{-1}(m,n)$ is a pair $v$
with $h_{}(v)=0$.
(That is, there exists a path from the  
left edge of $R$ to $l_v$, which path does not intersect any line of
$d$,
for $d \in [p]'$. 
The line $l_v$ is literally 
left-exposed.)   



\lem{\label{prp}
Let 
  $p_1\in J_{-1}(m,n)$ and 
$p_2\in J_{-1}(n,q)$.
Let $d_i \in [p_i]'$. 
Let $v$ be a left-exposed pair in $p_1$ or $p_2$
and let $l_v$ be the corresponding line in $d_i$. 
Then the   line  $l$ 
in $d_1\,|\,d_2$
containing $l_v$ 
 has  $h_d(l)=0$.  
} 
\proof
{A path connecting the line $l_v$ 
  to the left
edge without intersection in $d_1$ or $d_2$ also connects the line
in $d_1|d_2$ which contains $l_v$ to the left edge without
intersection
(cf. Fig.\ref{fig:wid1}).
%
\qed}     


\newcommand{\SL}{S^L}    

\mdef \label{de:2.18}
For $p \in  J_{-1}(m,n)$ let $\SL_p$ denote the subset of  left-exposed pairs.
Define
\beq  \label{eq:defblob}
\bB(m,n) \; =\; \{(p,s)\,|\,p\in J_{-1}(m,n), \; s\subset \SL_p \}
\eeq
%



\mdef
Let $(p_1, s_1) \in \bB(m,n)$ and $(p_2 , s_2) \in \bB(n,q)$.
Recall $p_1 . p_2$ and  $p_1 * p_2$ from (\ref{de:star}).
Define 
$\overline{s_1 s_2}$ to be the set of those
pairs $\{a_0,a_k\}$ in $p_1 . p_2$,
from chains
$\{a_0,a_1\},\{a_1,a_2\},\dots,\{a_{k-1},a_k\}$
where at least one pair comes
(in the obvious sense)
from $s_1$ or $s_2$.   
Define $\#$ as the number of closed chains with no pair from
$s_1, s_2$,   
and $\#'$ as the number of remaining closed chains.

We fix $\delta,\delta'\in k$ and define
the composition $\bx: \bB(m,n)\times \bB(n,q)\to k\bB(m,q)$
by
$$
(p_1,s_1) \bx (p_2,s_2)
 = \; \delta^\# \delta'^{\#'}  (p_1 . p_2,\overline{s_1 s_2}) . 
$$






{\theo{ Fix a commutative ring $k$ and $\delta,\delta' \in k$. 
    Then $\bb = (\N_0, k \bB(n,m), \bx )$ is a category.
}}

\proof We require to prove associativity of the product, and this
follows analogously to (\ref{de:stary}) and (\ref{th:B}).
Again the bookkeeping of extra primes may be seen from suitable pictures.


\mdef A picture with blobs is a pair $(d,b)$ where $d$ is a picture
in the sense of (\ref{de:picture});
and $b$ is a set of points (called {\em  blobs})
in the interiors of  lines of $d$.
We require that each blob lies in exactly one line
(if $d \in [p]'$ with $p \in J_{-1}$ then lines here are
non-crossing and this is automatic).

\mdef
Given a picture with blobs $(d,b)$
such that $h_d(x \in b) =0$ 
we define 
$$
\pi':(d,b) \mapsto (\pi(d),s(b)) \; \in \bB(m,n)
$$ 
where
$\pi$ is 
as in (\ref{de:pi}) 
and $s(b)$ is  the set of pairs
associated to the lines  decorated by $b$.

\mdef
An element $(p,s)\in \bB(m,n)$ can be represented by a pair $(d,b)$,  
where $d$ is a no-loop picture of $p$;  
and $b$
consists of at least one point in the interior of each line of $d$
corresponding to a pair in $s$. 
See Fig.\ref{fig:many2}.
%


\newcommand{\sss}[1]{\{ #1 \}}  

For an example note that Figure~\ref{fig:which1}
 is a picture $(d,\{P_1,P_2,P_3\})$ for an element
$$
(p,s) = 
(\{
   \sss{2',1'},
   \sss{9,10},
   \sss{3,4},\sss{7',5},
   \sss{1,2},\sss{8,8'},\sss{3',6'},\sss{5',4'},\sss{6,7} \},
   \;
     \{ \sss{5,7'}, \sss{1,2}, \sss{6',3'} \})
$$  
$\in \bB(10,8)$.
The blobs $P_i,i=1,2,3$ are mapped to the
pairs
$\sss{2,1}$, $\sss{5,7'}$, $\sss{6',3'}$
respectively. 


\mdef \label{de:Pipic}
Define
$$
\Pi(d,b) \;\; := \;\; \delta^{\#(d,b)} \delta'^{\#'(d,b)}  \;
\pi'({d},{b})
\;\; \in k \bB(m,n),
$$ 
where 
$\#(d,b)$
(respectively  
$\#'(d,b)$)
is the number of loops without (respectively with) blobs.

 
Let $d_i\in [p_i]'$ and $(d_i, b_i)$ be
no-loop
pictures with appropriate
blobs
(hereafter we just write $d_i$, including blobs).
By Lem.\ref{prp}
the concatenated picture
 $d_1\,|\,d_2$
is a picture of some 
$p \in J(m,q)$ plus possible loops and blobs. 
%
%
By an argument similar to (\ref{de:stary}) we have, independently of choices,
$$
(p_1,s_1) \bx (p_2,s_2) = \; \Pi(d_1 | d_2)  .
$$

\redx{
The construction along with Lemma \ref{prp} establishes  the following. 
{\lem{
For any pair of elements $(p_1,s_1)\in \bB(m,n)$ and $(p_2,s_2)\in
\bB(n,q)$ and pictures $(d_1,b_1),(d_2,b_2)$ with  
$\pi'(d_i,b_i)=(p_i,s_i)$  we have 
\[
\Pi(d_1\,|\,d_2,\overline{b_1\cup b_2}) = \Pi((d_1,b_1)) \bx \Pi((d_2,b_2))
\]
\qed
}}
}


Existence of constructs of form $d_1 | (d_2 | d_3) $ will again be evident.
Since again $d_1 | (d_2 | d_3) = (d_1 | d_2) | d_3$ we are done.
\hfill \qed


{\cor{ The End sets $k\bB(m,m)$ have the structure of an associative algebra. }} 

We may denote these algebras by
$J^{\bullet\,\delta,\delta'}_m$,
or simply $\bb_m$,
indicating the fixed 
parameters from $k$ in the definition of the 
multiplication only when needed for clarity.  


\begin{figure}
  \begin{center}
    \includegraphics[width=1.682cm]{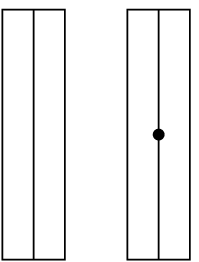} \hspace{1.3in}
    \includegraphics[width=6cm]{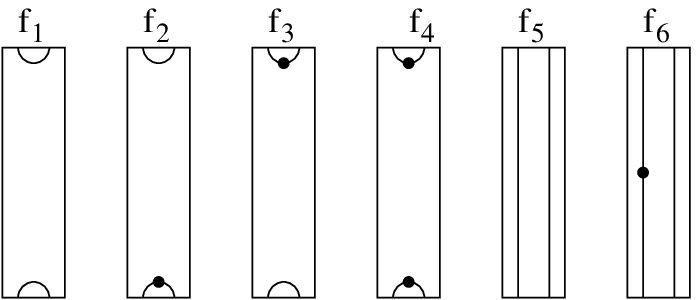}
  \end{center}
\caption{ 
(a) The  basis  $\bB(1,1)$. \hspace{.86in}
(b) The basis  $\bB(2,2)$. 
\label{fig:ffee} 
\ignore{{
\newline
Applying the rules of multiplication
outlined in (\ref{}), we have e.g., $e_1e_2=e_2\,e_1=e_2^2=e_2$,
$f_1^2=\delta f_1$, $f_2^2=\delta'f_2$, $f_2 f_1=f_2 f_3=\delta' f_1$,  
$f_6 f_1=f_6 f_3=f_3$ and so on. 
}}
}
\end{figure}


\mdef  
Examples. Consider Fig.\ref{fig:ffee}(a,b). 
Applying the rules of multiplication
we have e.g.,
$1 e = e 1 = e e = e$,
$\;  f_1^2=\delta f_1$, $f_2^2=\delta'f_2$, $f_2 f_1=f_2 f_3=\delta' f_1$,  
$f_6 f_1=f_6 f_3=f_3$ and so on. 






\subsection{Generators and relations} \label{ss:gnr}


{\theo{\label{th:blobgen1} {\rm \cite{\ms,Martin0706}}
Consider the algebra defined by generators 
$U^e  = \{ e,U_1,U_2,...,U_{n-1} \}$
and relations 
\\
$\tau' = \{ \mbox{
$U_i^2 = \delta U_i$,
$U_i U_{i\pm 1} U_i = U_i$, 
$U_i U_j = U_j U_i$, $j\neq i\pm 1$, 
$e e = e$,
$U_1 e U_1 = \deltap U_1$,
$U_i e = e U_i$, $i>1$
} \}$.
(I)
The map 
$$
U_i \mapsto   
\{ \underbrace{ \{\{1,1' \} \{2,2' \}, ... , 
\{i , i+1 \},\{ i',i+1' \}, ...,
  \{n,n' \}\}}_{p}, \underbrace{  \emptyset }_s \}
$$
$$
e \mapsto \{ \underbrace{ 
      \{ \{ 1,1' \},\{2,2'\}, ...,\{n,n' \}\} }_p,
 \underbrace{  \{ \{1,1'\}\}  }_s  \}
\hspace{1in}
$$
extends to an algeba isomorphism 
$k\langle U^e \rangle/\tau' \cong \bb_n$.
\\
(II) Every element of the partition basis $\bBd(n,n)$ can be expressed
as a word in these generators.
In particular $(p,s) \in \bBd(n,n)$ can be expressed as a word in
which $e$ appears as a factor $|s|$ times.
\qed
}}




\[ \]




\subsection{Generator algebras $B_{l,m,n}$ and spin chain Physics}\label{ss:q}
The original motivation for the blob algebra \cite{\ms} was to study XXZ spin
chains (and other related spin chains
\cite{AlcarazQuispell,LiebMattis66}) 
with various boundary conditions via representation theory. 
In the simplest formulation
(see \cite{PasquierSaleur90,Levy91,\ms} for details) 
one notes that there are boundary
conditions for which the $n$-site XXZ chain Hamiltonian may be expressed in the
form
\[
H = \sum_{i=1}^{n-1} U_i
\]
where $U_i$ acts on $(\C^2)^n$ by 
\[
U_i = 1_2 \otimes 1_2 \otimes ... \otimes 
       \mat{cccc} 0 \\ &q&1 \\ &1&q^{-1} \\ &&&0 \tam 
        \otimes 1_2 \otimes ... \otimes 1_2
\]
It is known that these matrices give a faithful representation of the
Temperley--Lieb algebra. In this sense the algebra controls the
spectrum of $H$. A centraliser algebra is $U_q sl_2$, and so this
can equivalently be seen as controlling the spectrum.

There are several reasons for wanting to generalise away from this
particular choice of boundary conditions. For example:
(1) periodic boundary conditions may be desirable for reasons of
computability or to minimise boundary effects at finite size.
(2) one may be interested in critical 
    bulk physics in the presence of a doped boundary.
(3) one may be interested in physics on the surface of a bulk system
(typically perhaps a 2D surface in 3D, but {\em modelled} 
more simply by a thickened finite interval at the end of an infinite
line). 

The generalisation required for the periodic case requires quite
delicate tuning --- see \cite{\ms}. But 
the simplest form of such generalisation is simply 
to modify the first or last
operator in the chain:
\[
H' = U'_0 + \sum_{i=1}^{n-1} U_i
\]
where $U'_0$ acts only on the first tensor factor
(and then of course to take the thermodynamic limit) \cite{Levy91}. 
The corresponding extension of the Temperley--Lieb algebra is covered
by the blob algebra. 
Another  
challenge is to dope with a more complex operator at the boundary. 
Algebraically this generalisation can get difficult quite quickly. 
A version which at least lies within the Brauer algebra is to
have a Temperley--Lieb chain with some permutation operators at the
end. 
That is, one first considers the Brauer algebra $B_n$ as the algebra
generated by its sub-symmetric group and Temperley--Lieb 
Coxeter generator elements
\begin{equation} \label{eq:Bn}
B_n = \langle \sigma_i , U_i \; : \; i=1,2,...,n-1 \rangle
\end{equation}
Note that this is not a minimal generating set. 
For example, all but one $U_i$ can be discarded.

\mdef \label{de:coxeter}
Trivially one can then define for each $n$ a `Coxeter subalgebra'
\[
B_{l,m,n} = \langle \sigma_i ,  \; : \; i=1,2,...,l-1, \;\; 
U_i ,  \; : \; i=m,m+1,...,n-1 
\rangle
\]
It is clear that various values of $l,m$ reduce to known cases. 
For example if $m>l$ then we just have a product of 
$S_l$ and $T_{n-m}$. So the interesting cases are $m\leq l$. ...

The Hamiltonians for such systems have been considered
\cite{Bondersan}, but in the present work 
we focus on the abstract algebraic aspects. 

It is clear that $B_{l,1,n} \hookrightarrow \BBl{l,n}$;
and that $B_{2,2,n} \hookrightarrow \BBBl{0,n}$. 
It is conjectured that these inclusions are isomorphisms. ...
And in this spirit we can ask about a geometrical 
and categorical characterisation of
$B_{l,l,n}$. 

\medskip




\newcommand{\alc}[2]{A_{#1}(#2)}
\newcommand{\LS}{left-simple}  

\begin{figure}
\[ 
\includegraphics[width=2.08in]{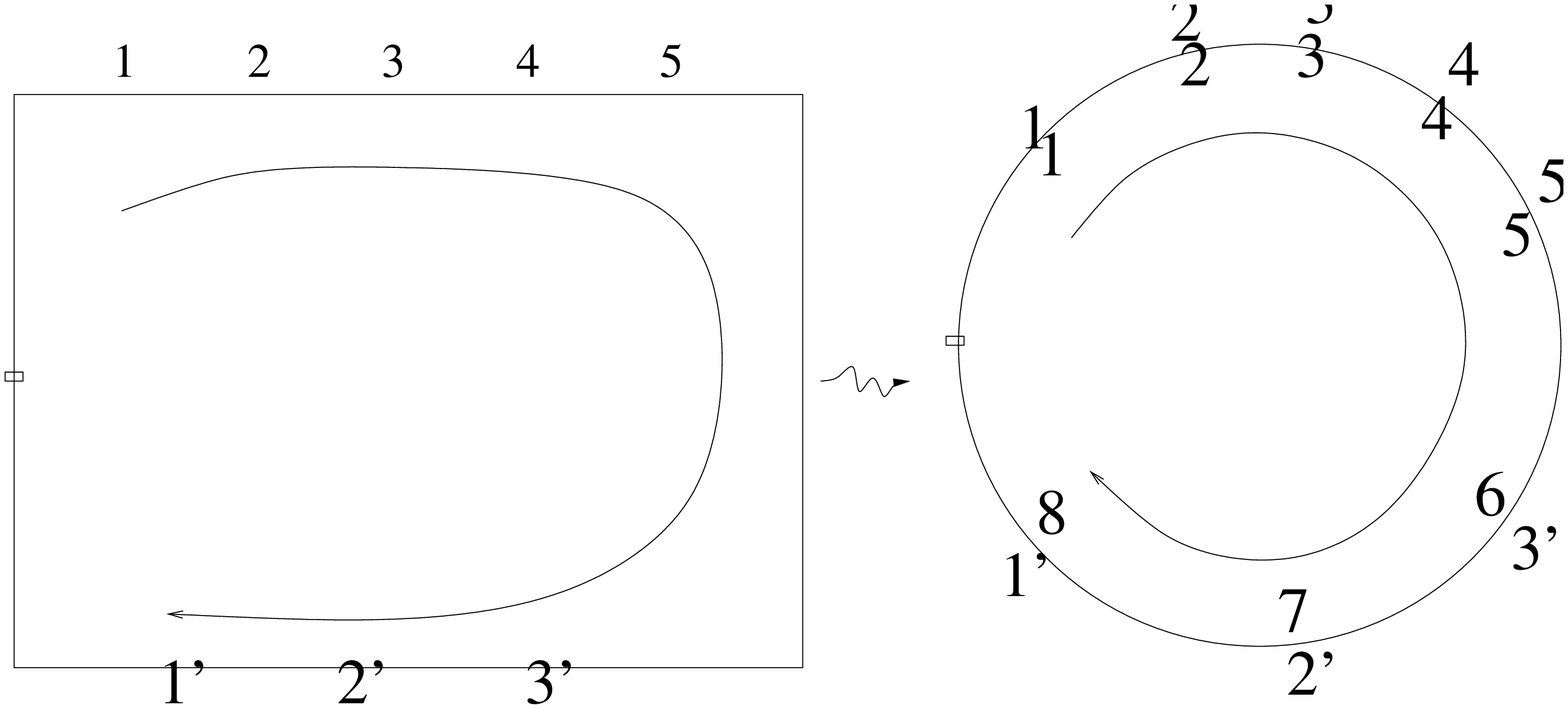}
\]
\caption{Disk order on vertices. \label{fig:disko1}}
\end{figure}

\mdef \label{def:diskorder}
The {\em disk order} on elements of $\{ 1,2 ,...,m \} \cup
\{1',2',...,n' \}$ is given by renumbering $i' \mapsto m+n+1-i$.
In our convention for pictures of $p$ 
this is clockwise order on the topological marked disk $R$ --- see
fig.\ref{fig:disko1}.
For $\{ i,j \}$ a pair in $p \in \JJ(m,n)$, 
with $i <j$ in the disk order, we understand by $[i,j]$
the interval from $i$ to $j$ with respect to the disk order.



\section{
  The subcategory $\BBl{i-1}^i$ of $\BBl{i-1}$} \label{ss:s4}
\newcommand{\Jimi}{J_{i-1,m}^i}
\newcommand{\Jonip}{J_{0,n+1}^1}
\newcommand{\Biim}{B_{i+1,i+1,m}}
\newcommand{\Biin}{B_{i+1,i+1,n}}
\newcommand{\defeq}{:= }
\newcommand{\li}{$Li$-simple}
\newcommand{\lo}{$L1$-simple}
\newcommand{\lic}{$Li$-chain}
\newcommand{\loc}{$L1$-chain}
\newcommand{\Jim}{J^{i}_{\leq i-1}(m,m)}
\newcommand{\Jic}{J^{i,c}_{\leq i-1}(m,m)}
\newcommand{\Jica}{J^{i,c+1}_{\leq i-1}(m,m)}


\mdef \label{de:crossing}
Let $S$ be a totally ordered set and $p$ a partition of $S$ into
pairs (we have in mind the disk order as in (\ref{def:diskorder})).
Via the total order,
the restriction of $p$ to any two  pairs induces a partition of  $\{1,2,3,4 \}$.
The set of such partitions is  
$ \{
\{\{1,3\},\{2,4\}\}, \; \{\{1,2\},\{3,4\}\}, \;   \{\{1,4\},\{2,3\}\}
\}$.

A pair in $p$ is {\em crossing in $p$} if
there is another pair in $p$ such that the partition induced
by the restriction of $p$ to these two pairs 
 is     $\{\{1,3\},\{2,4\}\}$.

 Remark: The point of this terminology is that
 if the total order is the disk order then
 in every picture $d$ of $p$ 
 the line for the `crossing' pair must cross another line.
(This follows from the Jordan curve theorem \cite{Moise77}.)

\ignore{{
Recall that given a collection $p$ of disjoint pairs from a totally
ordered set, a subset $p'$ of these pairs is
      {\em non-crossing}  if the partition of $\{1,2,3,4 \}$ induced
      by the 
    restriction of $p$ to any two  pairs with at least one from $p'$ is either 
    $\{\{1,2\},\{3,4\}\}$ or   $\{\{1,4\},\{2,3\}\}$. 

Remark: The point of this terminology is that $p$ has a picture $d$ in
which the line for each pair in $p'$ does not cross any other line.
(This follows from the Jordan curve theorem.)
}}


\subsection{\label{chi} The crossing number $\chi_p$ of a partition}

\begin{figure}\[\includegraphics{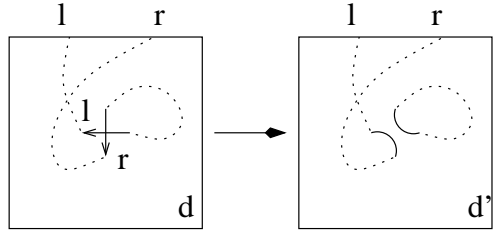}\]
\caption{\label{esi}
Schematic for removing a line self-crossing in a picture.
} 
\end{figure}

%
\begin{figure}\[\includegraphics{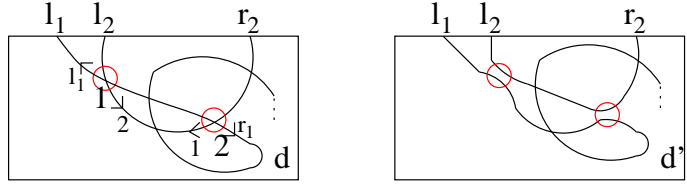}\]
\caption{\label{e2i} Cancellation:
Small neighbourhoods of two crossings of two lines in 
a picture $d\in [p]$ modified such that the resulting picture
satisfies $d'\in [p]$ with 
two fewer crossings.  
}  
\end{figure}

\mdef We denote the number of crossings of a picture $d\in [p]$ by $\#^c(d)$

\mdef \label{de:pcross}
Let $p$ be a Brauer partition and let
$\{\{l_1 , r_1 \},\{ l_2 , r_2 \} \}$ be two pairs in $p$.
Note by the Jordan Curve Theorem that if precisely one of $l_2,r_2$
lies in the (disk)-interval from $l_1$ to $r_1$ then a picture of $p$
must have at least one crossing of the lines corresponding to these
two pairs --- we say the two pairs `cross'.
We thus have a lower bound on the number of crossings in a
picture $d$ of $p$:
\[
\# \geq \chi_p
\]
where
\[
\chi_p = \sum_{ \{\{l_1 , r_1 \},\{ l_2 , r_2 \} \} \subset p}
                      \delta(\{\{l_1 , r_1 \},\{ l_2 , r_2 \} \} )
\]
where the sum is over pairs of pairs and $\delta(-)$ is 1 if they
cross and zero otherwise.


{\lem{ \label{minex} Consider $p \in J(n,m)$.
    There exist pictures $d$ of $p$ achieving the minimum
    $\#^c(d) = \chi_p$.}}

\proof  First consider a picture $d \in [p]$ with a line self-crossing. 
Then note from (\ref{de:line}) that removing the open 
segment of the line corresponding to the loop starting and ending at
the self-crossing produces  another  picture $d^{\circ}$, 
such that $\pi(d^{\circ})=\pi(d)$.
Next, see
Fig.\ref{e2i}.
It shows that whenever there are two crossings
of the same two lines in $d\in [p]$ these
crossing neighbourhoods may be
`cancelled' 
to make $d' \in [p]$ with two fewer crossings.
\qed

\lem{\label{pmin} Consider $p \in J(n,m)$.
  There exist low-height pictures of $p$ with $\chi_p$ crossings.}
\proof{

We claim that we do not increase the height by either $d \leadsto d^{\circ}$ or by
$d \leadsto d'$ in (\ref{minex}) above.
The first is due to Lemma \red{2.31} of \cite{\kamy}.
For $d \leadsto d'$ we proceed as follows. Consider a 
crossing $x$ of $d$ that remains a crossing in $d'$.
Consider a specific path in $d$ from $x$ 
 to the left edge. 
Note that this path is also a path in $d'$ provided that the 
differing neighbourhoods are taken small enough so that the path avoids them. 
Therefore the  
number of crossings of the path with segments of the picture
is identical in $d'$.
This puts an upper bound on the height of this crossing point in
$d'$. 
Since this applies for each crossing point, the height of $d'$ is
bounded above by the height of $d$. 
In particular if $d$ is low-height then so must $d'$ be.
%
%
\qed
}

\mdef{
  Let  $[p]''$ denote the set of low height pictures of $p$ with minimum
  number of crossings.} 


\mdef \label{de:gardenpath}
The following procedure will be useful later. 
Given a picture $d$ for a non-empty non-crossing partition $p$ where $d$ has a loop,
as when a loop appears in composition,
we can modify the picture by a `garden path' as here:
\[
\includegraphics[width=1in]{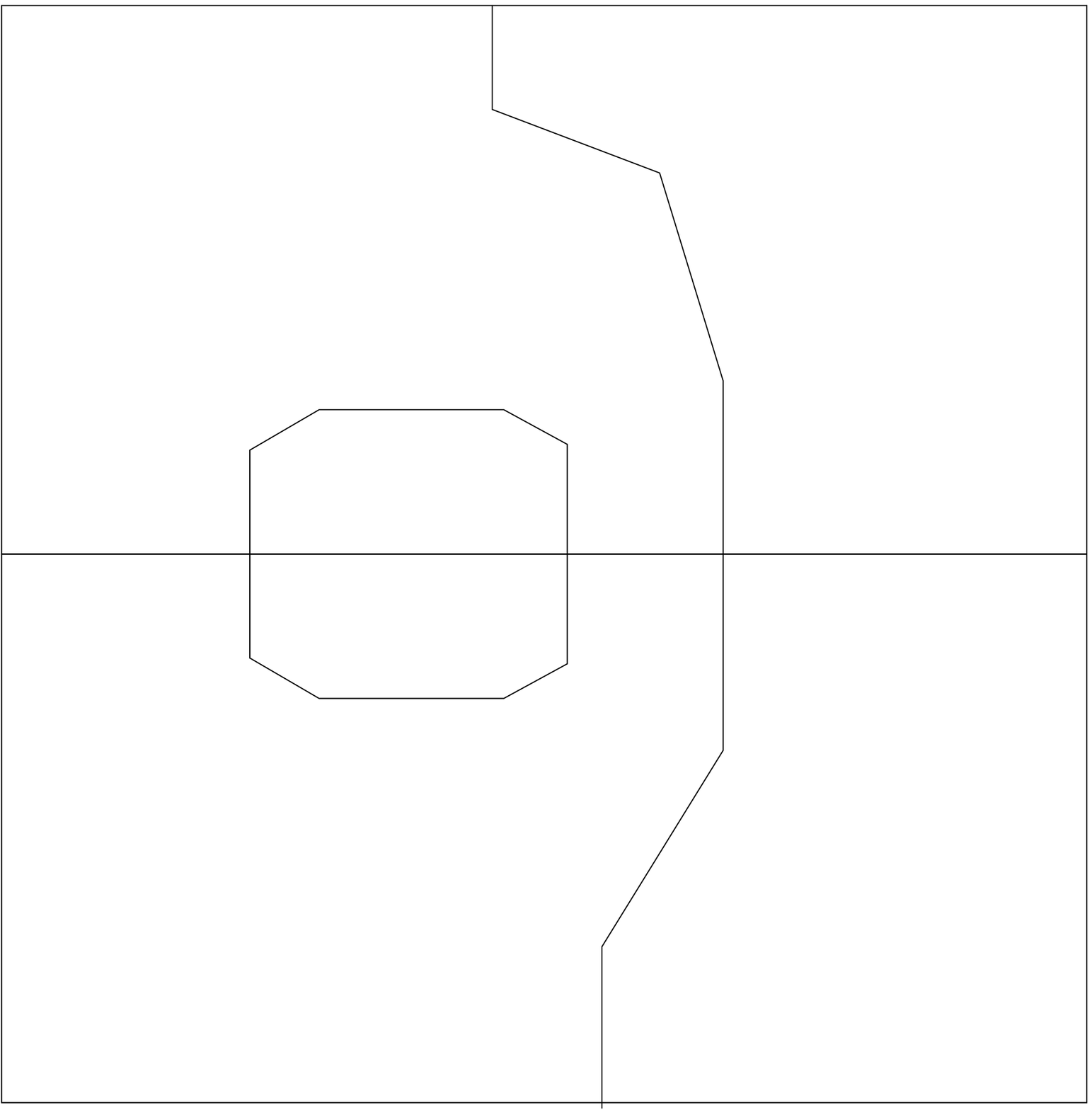}
\leadsto
\includegraphics[width=1in]{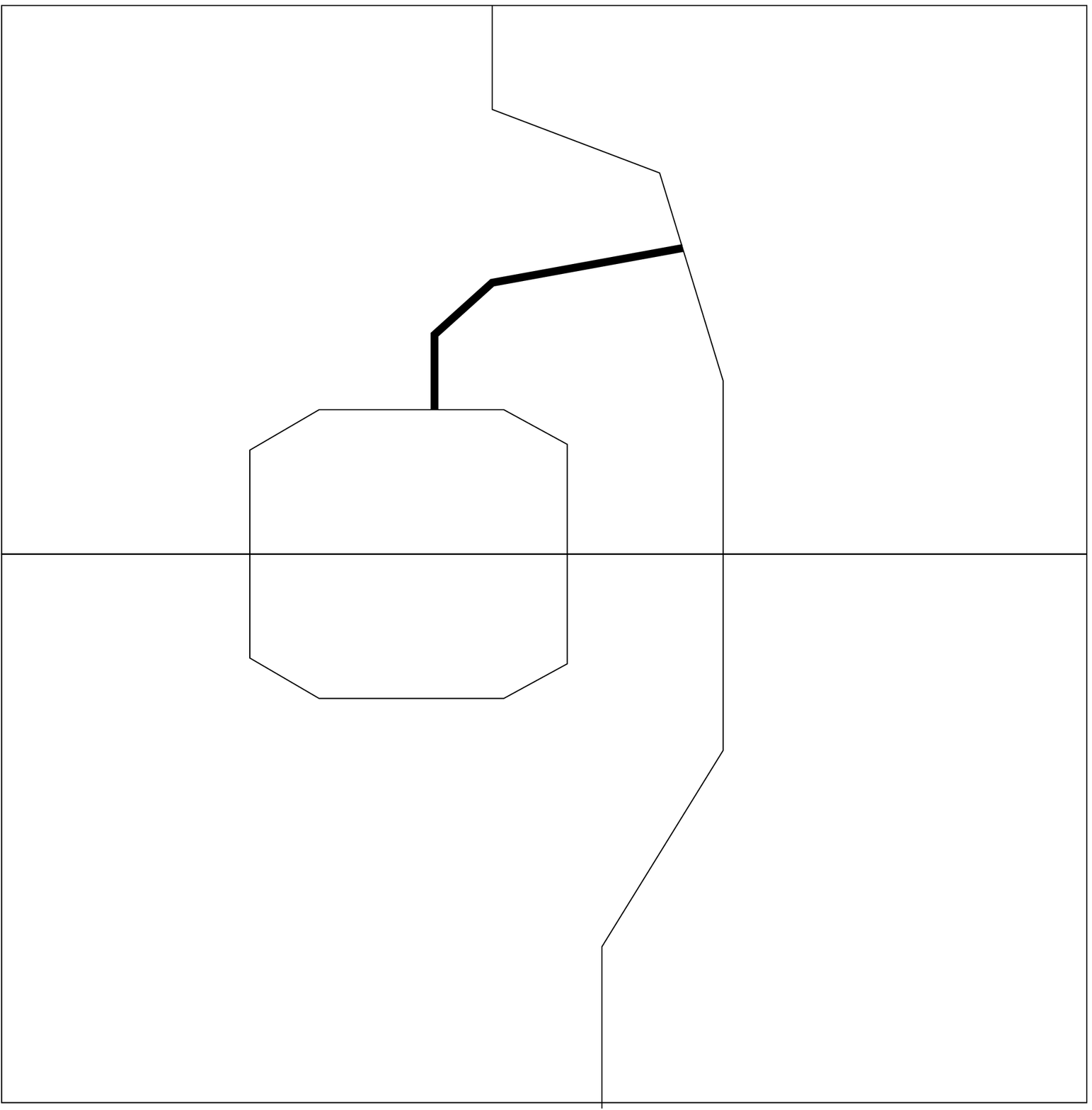}
\leadsto
\includegraphics[width=1in]{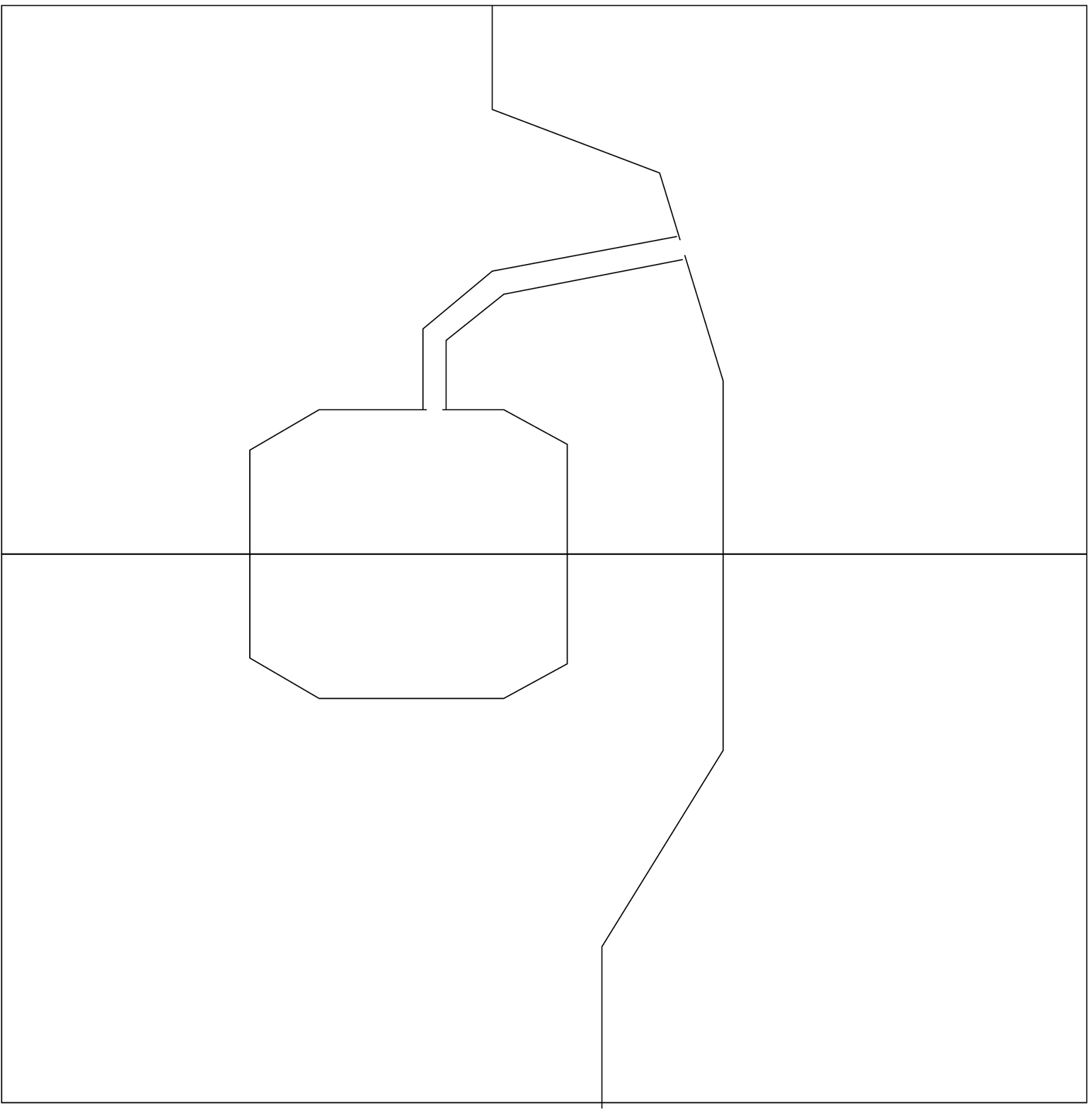}
\]
so the new picture is non-crossing, for the same $p$, and has no loop.


\subsection{\loc\ partitions}
\redx{Before all this, we should define chains and exclusiveness,
  maybe by taking over pics and paragraphs from p46.  
A pic would also be nice to illustrate the first lemma, but I think
it's best to do after we agreed/did above definitions.} 


\newcommand{\fd}{paul46/}
\newcommand{\figd}{./}

\newcommand{\kk}{k} 
\newcommand{\cnp}{c_p}  

\ignore{{

\begin{figure}
\[
\includegraphics[width=7.3cm]{\fd testc2.eps}
\]
\caption{ \label{fig:chainc1} }
\end{figure}

}}

\mdef \label{de:chain}
Consider the underlying set of a set partition, 
equipped with a total order --- for example 
recall the disk order on vertices as in (\ref{def:diskorder}).

Let $p$ be  a Brauer partition.
An ordered subset 
$\{ \{ l_1,r_1 \} , \{ l_2, r_2 \} , ... , \{ l_j, r_j\} \}$ 
of $p$ is called a
{\em chain} (of length $j$)
if   
 $l_{i-1} < l_i < r_{i-1} < r_i$
for all $i$. 

\mdef \label{pa:once}
Note that a chain $c$ has a `canonical'
(isotopy class of) picture in which
the lines for each
successive pair of links cross exactly once.
Indeed this is the unique $\chi$-minimal class.


\mdef \label{de:01210}
A chain $c$ in a partition $p$
divides $\partial R$ (and the corresponding ordered vertex set)
into intervals of height 0,1,2 with
respect to this chain only --- denoted $ht_c$.
In particular the sequence of these heights
from the top-left of $\partial R$ is
$$
seq_c \; = \; 01(21)^{j-1} 0
$$
\\
In the example in fig.\ref{fig:ecd}(a) the $ht$s are given for the
thin line
chain only.

These heights agree with the heights of point in $\partial R$ in a canonical
picture $d$ in the usual sense,
but once $c$ in $p$ is given then 
we may also consider them as invariants of $p$
--- see (\ref{re:bhi}). 
\redx{This para needs tidying.
Our def of pictures should ensure that heights on the boundary are
canonical
and indeed combinatorial. Consider a path to the relevant alcove that
stays very close to the frame...
See also my online notes \cite{Martin07ol}.}

\mdef
A partition $p$ is called {\em \loc} if $ht(p) \leq 0$ and 
there is a chain with $l_1 = 1$ and $r_j = 1'$. 

\mdef
A partition $p$ is called {\em \lo} if 
there is a chain with $l_1 = 1$ and $r_j = 1'$. 

Write $J^1(n,m) \subset J(n,m)$ for the set of \lo\ partitions.



{\lem{ \label{le:ppnc}
Every \loc\ partition $p$ consists of a unique 
chain from $1$ to $1'$ together with a
collection of pairs that are non-crossing in $p$.
 \ignore{{
(in particular every other chain has length $j=1$).

    (a2) If $p$ is \loc\ then it has a
(`canonical' isotopy class of)
   low-height (and $\chi$-minimal) picture with exactly one
crossing between the lines for each successive pair of `links' in the
chain; and no other crossings.
}}
}}
\proof
Firstly $p$ has at least one chain from 1 to $1'$ by definition.
Pick one such, and 
start with a canonical picture $d_c$
of this chain $c$ and the alcoves it defines.
Schematically an interval of the chain looks like:
\[
\includegraphics[width=1.682in]{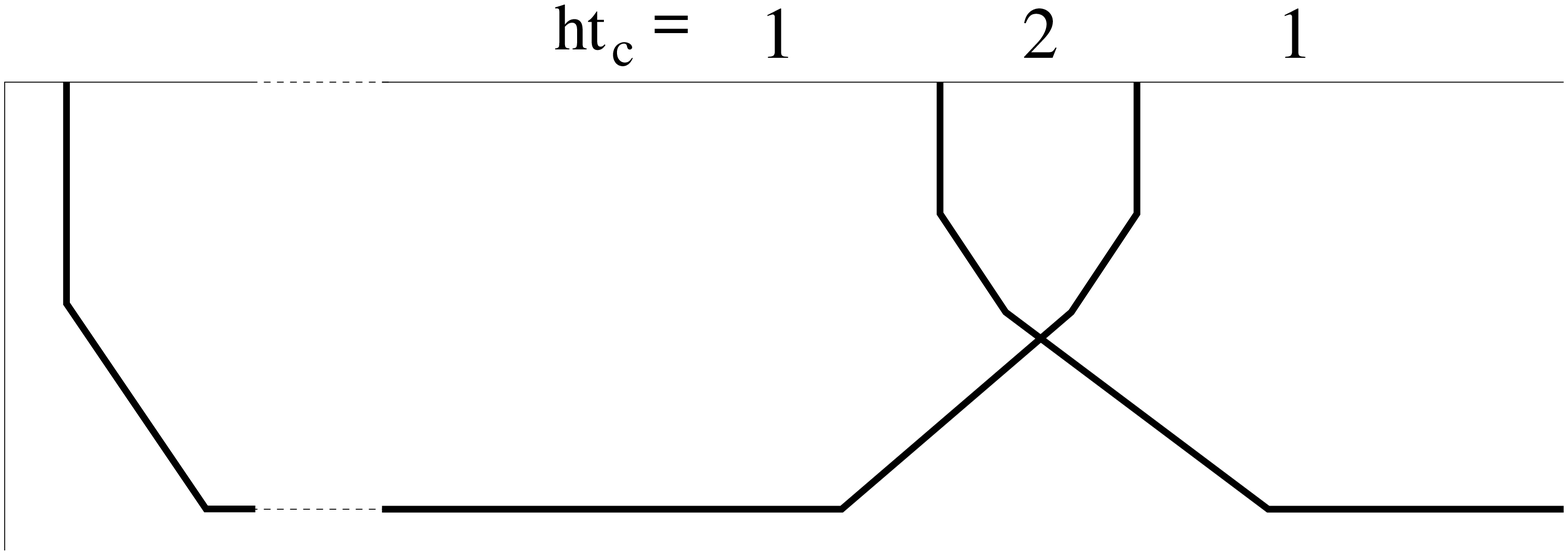}
\]
Consider a pair in $p$ that is not part of the chain.
By (\ref{de:01210}) a path 
for this pair,
in any picture extending $d_c$,
starts in an alcove of $ht_c$ value 1 or 2.
Such a path itself defines at least one new alcove, so the
true heights in some region over the corresponding interval are 2 or 3.
One sees:
\[
\includegraphics[width=1.4682in]{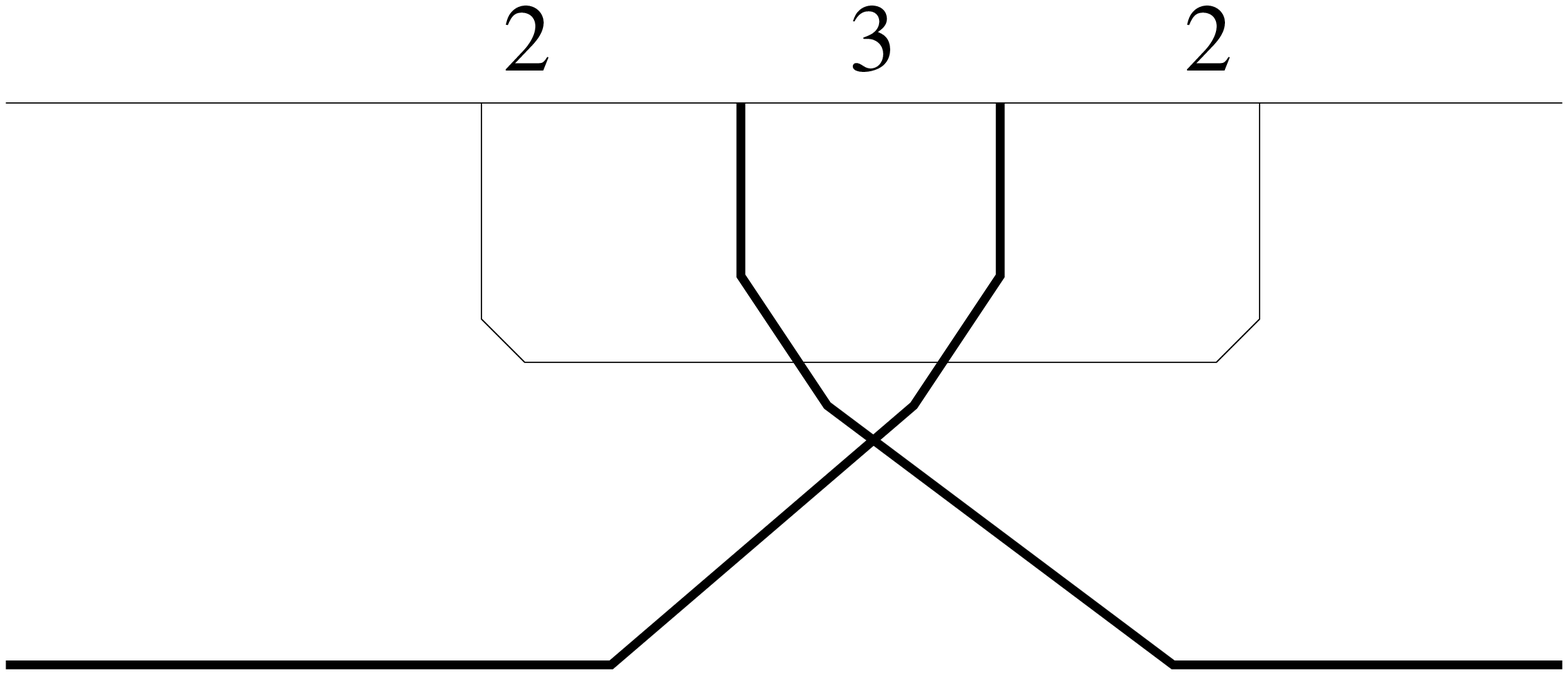}
\hspace{.521in}
\includegraphics[width=1.4682in]{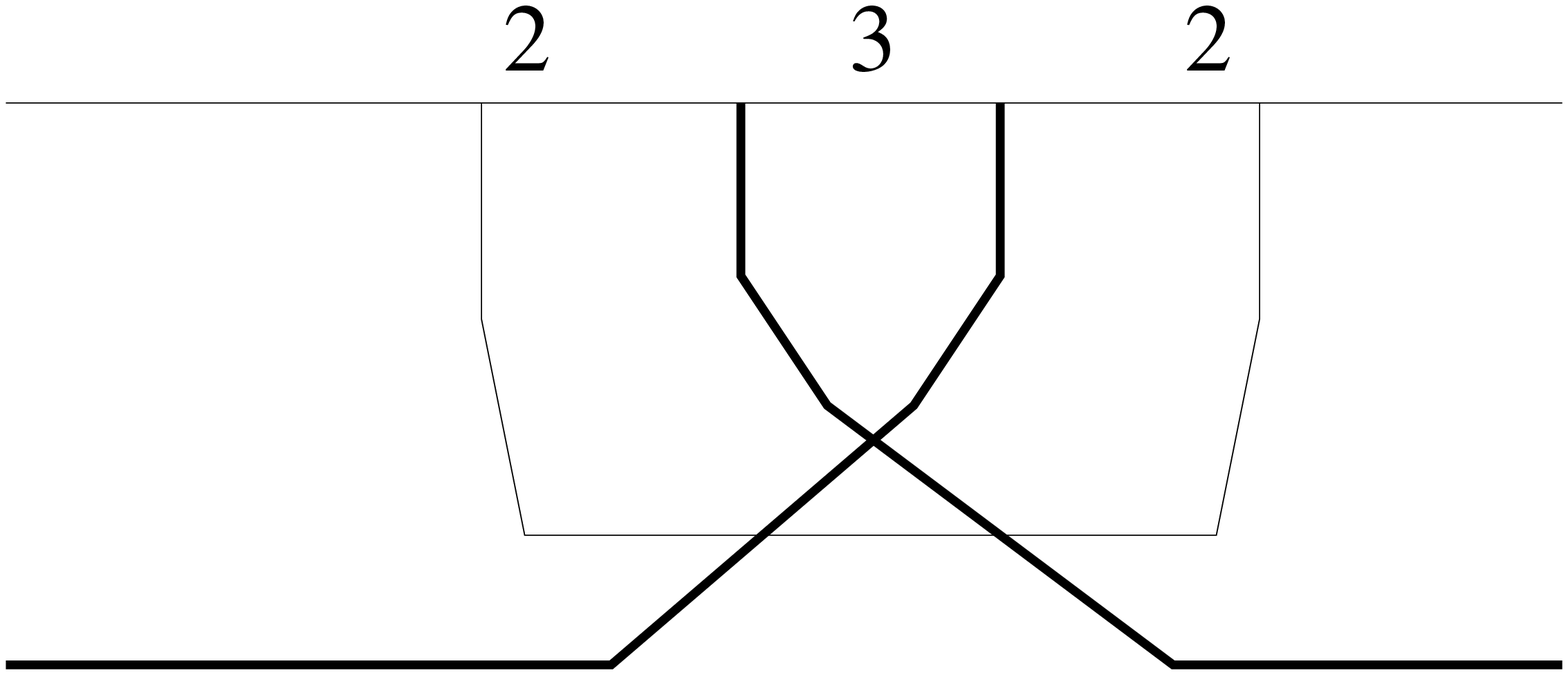}
\hspace{.521in}
\includegraphics[width=1.4682in]{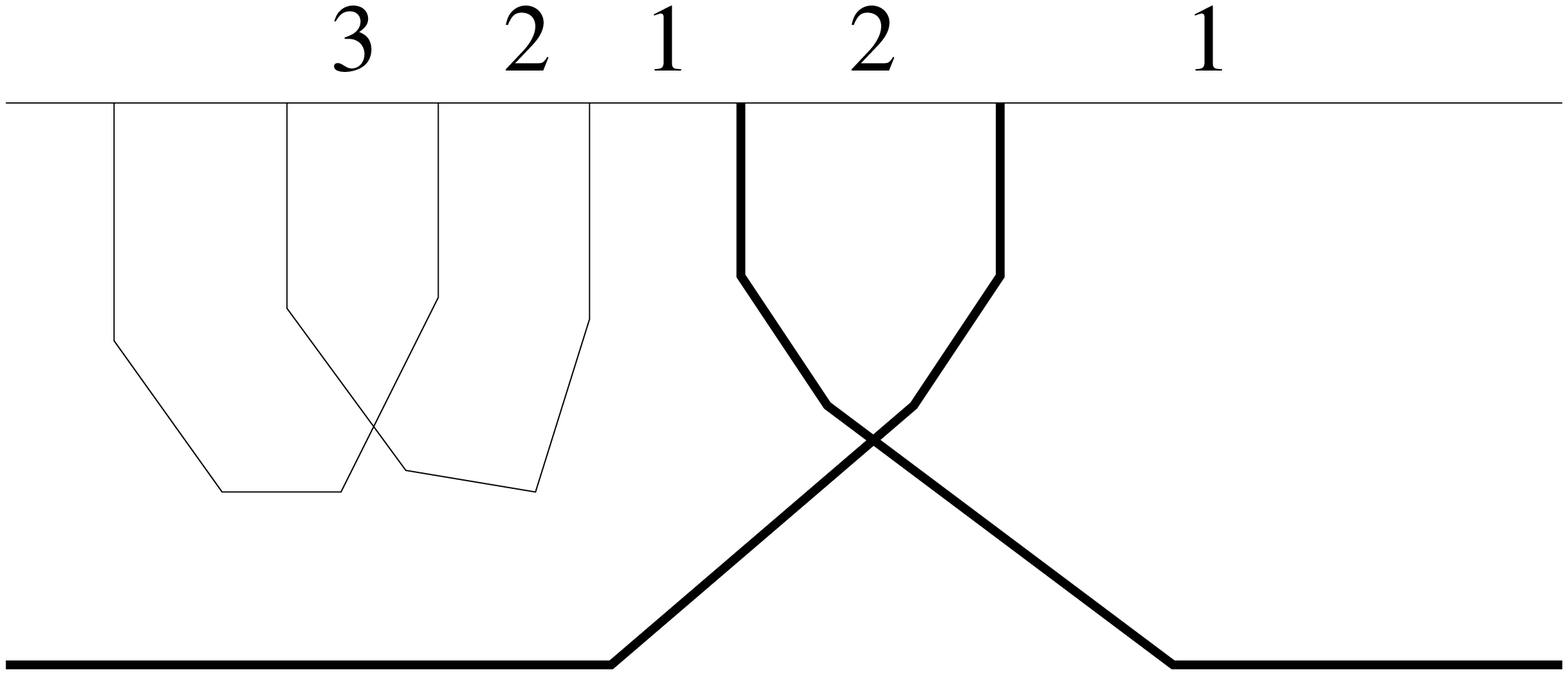}
\]
that such a pair, if it is crossing in $p$,
must give rise to a crossing touching height 3.
Note more generally that if a crossing touches an alcove of height $a$, say, then it has
height $\geq a-2$.
Thus a crossing touching height 3
contradicts the $ht(p) \leq 0$ condition.
Thus a pair in $p\setminus c$ is non-crossing. 
Finally no non-crossing pair can be part of a `long' chain, so $c$ is
unique.
\qed

\ignore{{
Finally consider a pair whose boundary interval lies in a $ht_c =1$
region. If this is crossing with another such then the true height in
the intersection of their boundary intervals is 3, leading again to a
contradiction of the $ht(p) \leq 0$ condition. 

\red{hmm. this is not working so well.
}

...

If it starts in height two then any crossing would violate the height
condition, so it is non-crossing.
If it starts in height one then 

\red{old version:}
(a1): this follows from (a2).
(a2): Start with a canonical picture of the chain, as in
(\ref{pa:once}),
and consider the alcove decomposition of the rectangle $R$ given by this.
By Lem.\ref{chi7} there exists a low-height $\chi$-picture $d$ of $p$
(and by the proof of 
Lem.\ref{chi6}
we can achieve this starting with a canonical picture of the chain).
Consider such a picture, and 
let $x,y$ be a pair in $p$ not in the chain.
If they lie in different alcoves then a path between them must cross
the chain path. It cannot cross at height 1, but if it crossed at
height 0 then this path covers a chain crossing, contradicting the
height 0 condition. Thus they lie in the same alcove.
That is, both lie between the same pair in the chain sequence ().
But then there is a path between them within a single alcove
(else $d$ is not $\chi$-minimal). This cannot cross the path for any other
such pair by the height 0 condition.
\qed
}}


\mdef Note that the
minimum number of crossings in a picture $d\in [p]$ here is $j-1$.


\ignore{{

\mdef Remark: 
The definition of \loc\ is equivalent to 
left-simple (as in (\ref{de:leftsimple})) and height 0. 

\proof It is clear that L-simple implies left-simple. The
converse is not true. However
\red{finish me}.

}}

\ignore{{

\mdef Claim: $B_{2,2,n} = \kk J^1_{\leq 0}(n,n)$

proof: We aim to prove that every partition in $ J^1_{\leq 0}(n,n)$
lies in  $B_{2,2,n} $. 
(This is enough since it means that $ \kk J^1_{\leq 0}(n,n)$
is closed under composition; and since the generators of  $B_{2,2,n} $
are clearly in, we have equality.) 

Work by induction on crossing number. The base case is 0.
The L-simple partitions of this form have a pair $\{ 1,1' \}$ and the
remaining pairs form a non-crossing partition.
It is well-known that these are generated by the corresponding TL
generators. 

So, assume that all partitions in $ J^1_{\leq 0}(n,n)$ with 
$\cnp <l$ lie in   $B_{2,2,n} $. 
Now consider a partition with $\cnp =l$. 
We have the following cases:
\[
pics
\] 

}}

\subsection{\lic\ partitions}
\begin{figure}
\[
\raisebox{-.051in}{
  \includegraphics[width=7cm]{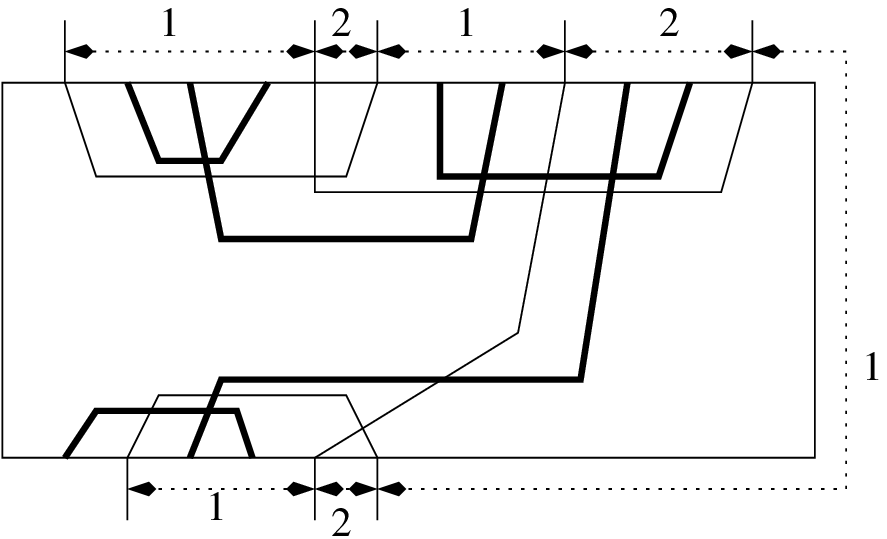}}
\hspace{.753in}
\includegraphics[width=4.35cm]{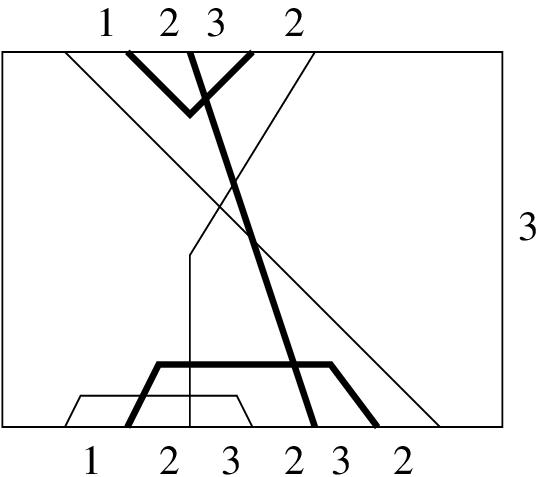}
\]
\caption{\label{fig:ecd} \label{fig:ecd2} 
Pictures of pair-partitions decomposable as two chains.
(a) Picture with $ht$ sequence shown for a single chain (thin lines).
(b) Picture with combined $ht$ sequences for both chains.
}
\end{figure}

\mdef
A pair of chains in $p \in J(n,m)$ is {\em exclusive}
if their individual $ht$ 2 intervals do not intersect.
\\
The example in fig.\ref{fig:ecd}(a) is not exclusive.
The example in fig.\ref{fig:ecd2}(b) is exclusive.

\mdef $ht$ with respect to a pair of chains together is
$ht_{c_1} +ht_{c_2}$. Thus exclusive pair gives hts up to 3,
and hence crossing hts up to 1.

\label{L2uniq}
Consider an exclusive pair whose initial points are adjacent and
whose final points are adjacent. Their combined boundary height
sequence is of the form
$$
seq_{c_1 \cup c_2} \; = \; 01(23)^l 210
$$
for some $l$.
A $ht$ 3 region is necessarily a link region
(a $ht$ 2 region with respect to some single chain);
and all links arise this way.
Thus if partition $p$ is an exclusive pair of this form
then there is only one way
in which it is an exclusive pair. 


\mdef
A {\em canonical drawing} of an exclusive pair $c_1,c_2$
has arcs of $c_2$ that contain ht 2 regions of $c_1$ drawn
{\em over} the corresponding crossings.


\mdef Define $\#_i : J(n,m) \rightarrow \N_0$
so that $\#_i (p)$ is the number of pairs that start in $\{ 1,2,...,i \}$
and end in  $\{ 1',2',...,i' \}$.

\mdef
A partition $p \in J(n,m)$ is {\em \lic} if it has $ht < i$
and 
$i$ chains that are pairwise exclusive,
starting in $\{ 1,2,...,i \}$ and ending in
$\{ 1', 2',...,i' \}$.

N.B. We consider pairs in $p$ that meet both
$\{ 1,2,...,i \}$ and $\{ 1',2',...,i' \}$
as chains of length 1, and hence trivially exclusive with any other
chain.

\mdef
Write $J^i_{< i}(n,m)$
or  $J^i_{\leq i-1}(n,m)$
for the subset of $J(n,m)$ of \lic\ partitions.


\beq \label{fig:ind1b}
\raisebox{-.31in}{
\includegraphics[width=1.72in]{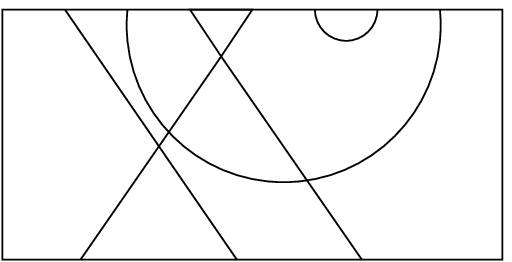}}
\hspace{.71in}
p = \{
\underbrace{\{ 1, 2' \}}_{c_1} ,
\underbrace{\{2,7\} , \{3,3' \}}_{c_2}, 
\underbrace{\{ 4, 1' \}}_{c_3} ,
\{5,6\} \}
\hspace{.2in}
\eq
\mdef Example. 
The figure in (\ref{fig:ind1b}) above gives a $p$ 
which is $L3$-chain, 
in $J^3_{<3}(7,3)$.

\mdef \label{de:uniq}
Given an \lic\ decomposition of a \lic\ partition $p$:
\[
p = \cup_{j=1}^{i} c_j \cup p'
\]
where the $c_j$s are  pairwise exclusive chains from $[1,i]$ to $[i',1']$,
then the restriction to
any two chains obeys the $L2$-chain property.
Thus the boundary height sequence of $p$ is of the form
$$
seq_{\cup_i c_i} \; = \; 012...i(i+1 \; i)^l ...210  ,
$$
independent of any decomposition.
Hence, as in  \ref{L2uniq}, the link positions are determined, and the 
decomposition of $p$ is unique. 


{\lem \label{lem:ppnc}
Let $p\in J^i_{i-1}(m,n)$ with $ m,n\geq i$, so that 
$p=\cup_{j=1}^i  c_j\cup p'$
as above. 
Then $p'$ is a set of pairs not crossing with each other or any chain.
}

\proof Noting (\ref{de:uniq}), this is analogous to Lem.\ref{le:ppnc}.
\qed

\subsection{A $\Biim$-module of \lic\ partitions}
\begin{figure}
\begin{center}
  \includegraphics[width=5cm]{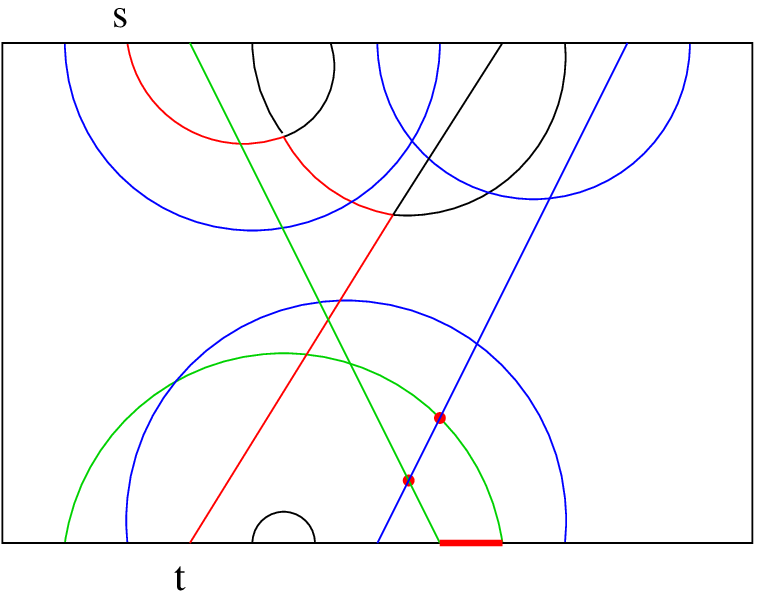}
  \caption{\label{jimi}
    Picture of a partition that is {\em not} L3-chain (blue and green
    chains not exclusive).
  }
\end{center}
\end{figure}

\medskip

Note that (for $m>i$) the algebra $ k J_{\leq i-1}(m,m)$ has a natural subalgebra
isomorphic to the symmetric group algebra $k S_{i+1}$.

{\lem\label{mj} 
The space  $kJ^i_{\leq i-1}(m,n)$   
is closed under the left action of $\Biim$.
}

{\proof
  Consider   $p\in J^i_{i-1}(m,n)$, and let $c_1, c_2,...,c_i $ be
the unique
  chain
decomposition of $p$ as in (\ref{de:uniq}). 
  
We can write out chains as sequences of pairs in the
chain order: $\{ j_1,j_2 \}, \{ j_3, j_4 \}, ... $, or even as
lists $j_1, j_2, j_3, j_4, ...$.

Consider  the action of generators
$\sigma_1, ... , \sigma_i$ and $U_j$ ($j>i$) on $p$ as follows:
\\
Case 1: $\sigma_j$ with $j<i$.
This changes only $c_j $ and $c_{j+1}$, swapping their first terms.
Since these terms are adjacent it follows that the \lic\ property is
preserved. 

\noindent
Case 2: $\sigma_i$. We can subdivide into three possibilities here.
\\
(i) vertex $i+1$ not in any $c_j$. In this case $i+1$ is in a
non-crossing pair by the \lic\ property (Lemma~\ref{lem:ppnc}).
Schematically, drawing only chain $c_i$ and the $i+1$ pair:
\[ 
\includegraphics[width=2in]{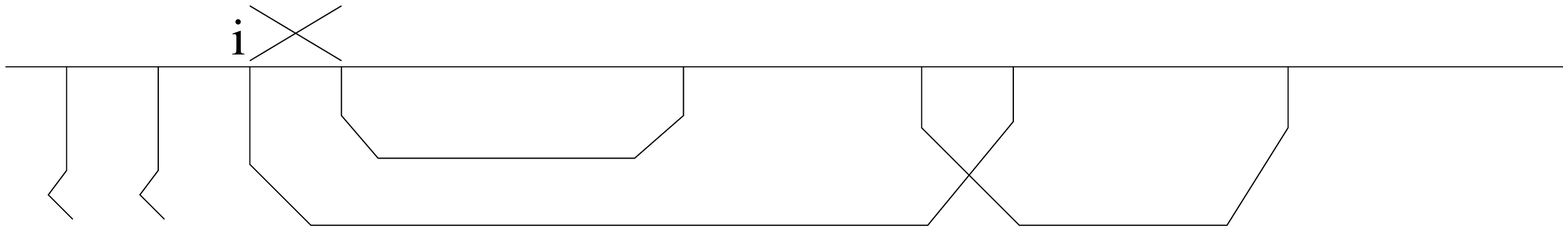}
= \;
\includegraphics[width=2in]{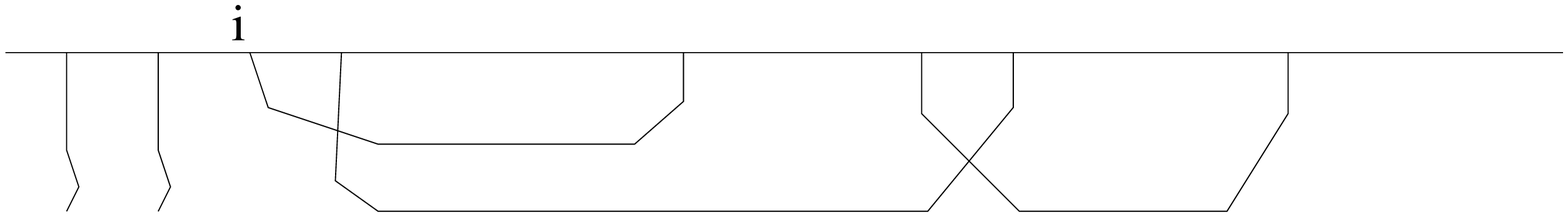}
\]
That is, the partition $\sigma_i p$ has a chain $c_i$ with an extra link.
Note that the new exclusive region for this chain comes from the
non-crossing part of the original partition. Thus it cannot overlap an
exclusive region for any of the undrawn and unchanged chains,
and the exclusive property is preserved.


\smallskip
\noindent
(iii) vertex $i+1$ in $c_{i}$.
Here there is a chain $\{ i,j \}, \{ i+1,l \},\{ k,m \},...$.
This becomes $\{ i,l \},\{ k,m \},...$ and leaves $\{ i+1, j \}$. 
The first of these is a chain from $i$. The second a non-crossing pair.
Schematically:
\[
\includegraphics[width=2in]{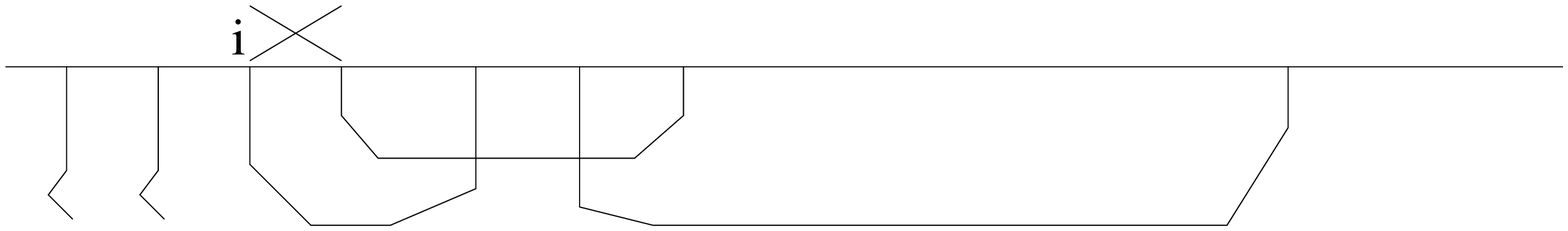}
= \;\;\;
\includegraphics[width=2in]{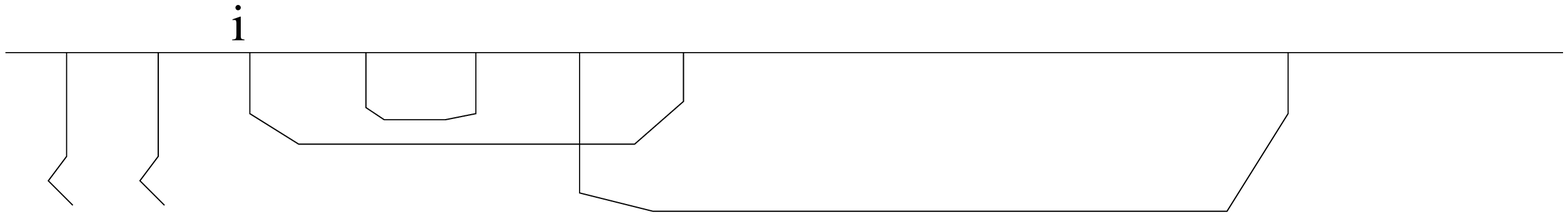}
\]
Here one link region is removed and so there is nothing to check for
exclusivity of $\sigma_i p$.


\smallskip
\noindent
(ii) vertex $i+1$ in $c_{j<i}$.
Note that $i+1$ cannot be in a pair directly with $j$ since then there is no
way to on-link this chain, so it must be the first vertex in the
second pair of $c_j$.
So here there is a chain $\{ j,k \}, \{ i+1,l \},...$ and a chain
$\{ i,r \},...$.
Schematically:
\[ 
\includegraphics[width=2in]{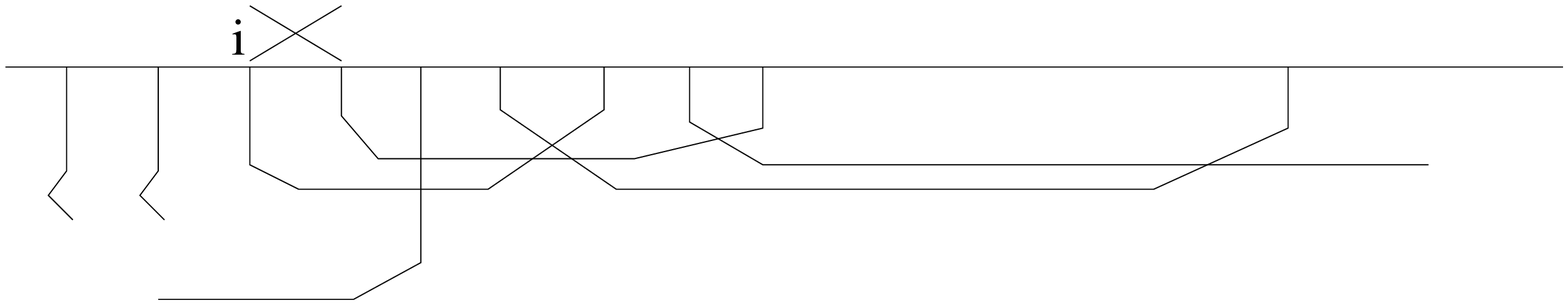}
= \;\;\;
\includegraphics[width=2in]{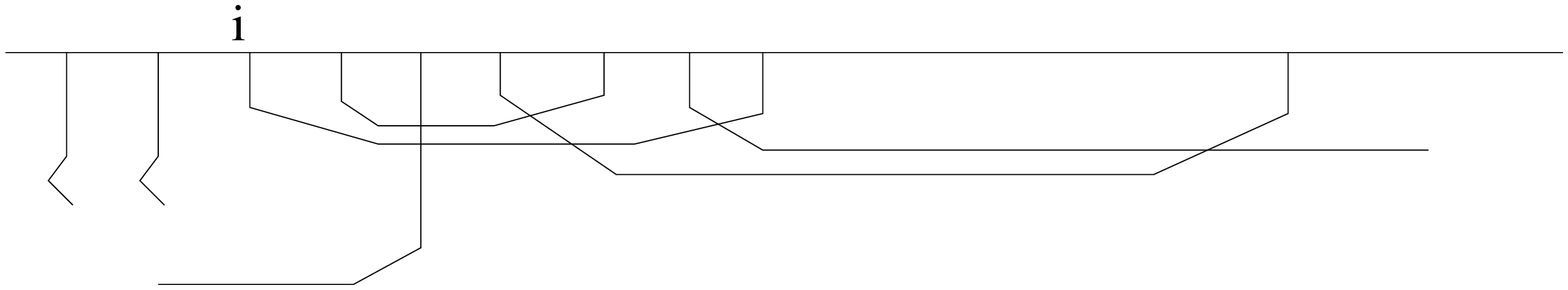}
\]
The new chains are $\{ i,l \},...$ and $\{ j,k \}, \{ i+1, r \}, ...$.
Note  that these are chains from the same starting points.
It is clear that they are pairwise exclusive;
and that pairwise exclusivity with other chains is not affected.


\smallskip
\noindent
Case 3: $U_j$ with $j>i$. The subcases here are
(i) $U_j$ touches no chain $c_k$. Schematically:
\[
\includegraphics[width=2in]{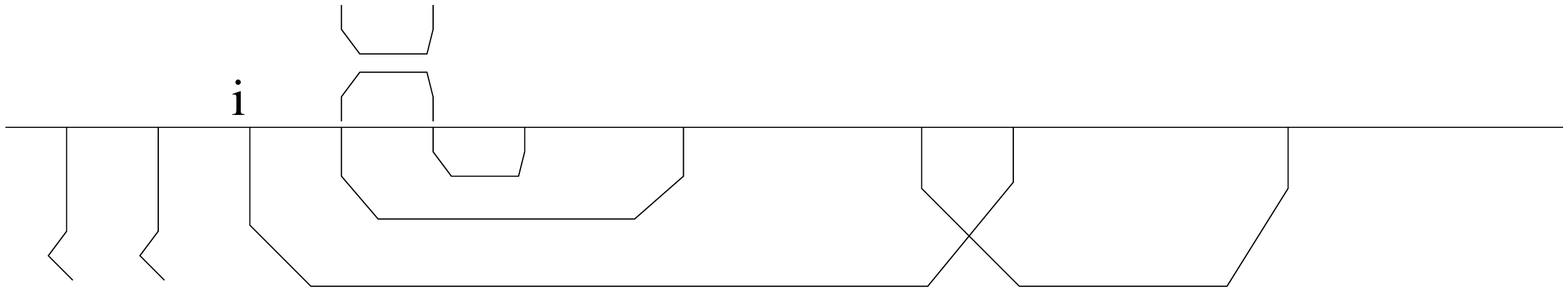}
= \;\;\;
\includegraphics[width=2in]{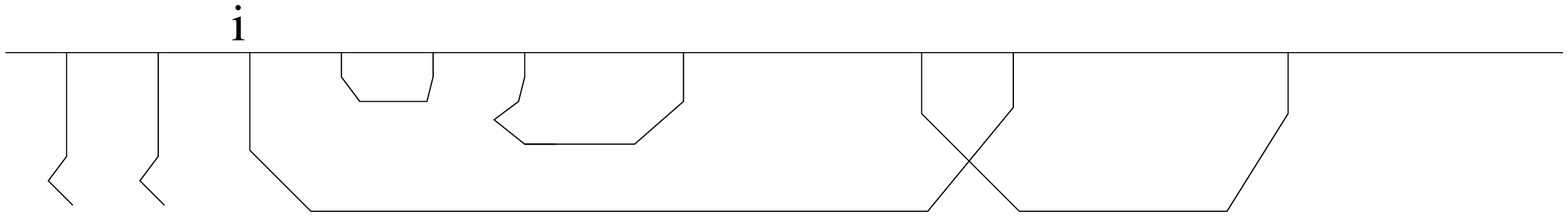}
\]
Note here that there is nothing to check for exclusivity of the new
partition.
\\
(ii) $U_j$ touches 
one vertex of a pair in some chain $c_k$, and no other:
\[
\includegraphics[width=2in]{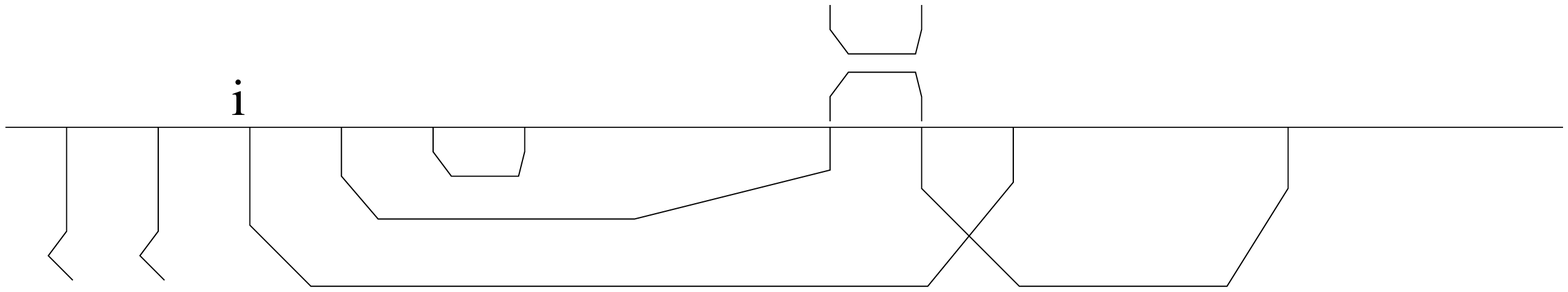}
= \;\;\;
\includegraphics[width=2in]{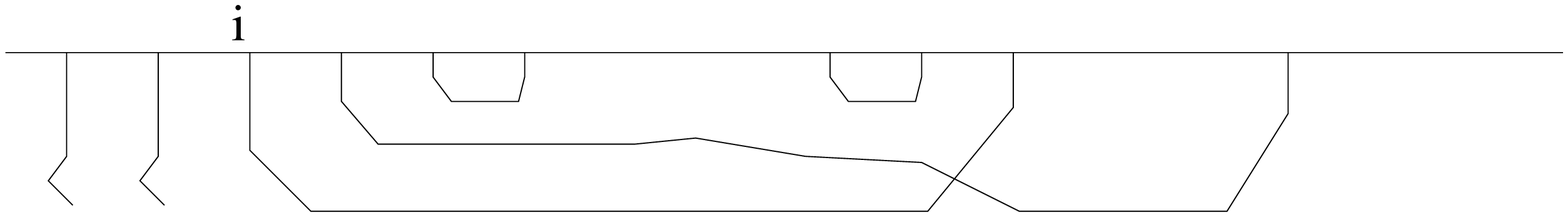}
\]
Here note that the size of a link region changes, but only into or out
of a non-crossing region. Thus exclusivity is not affected. 
\\
(iii)  $U_j$ touches some chain $c_k$ at two points:
\[
\includegraphics[width=2in]{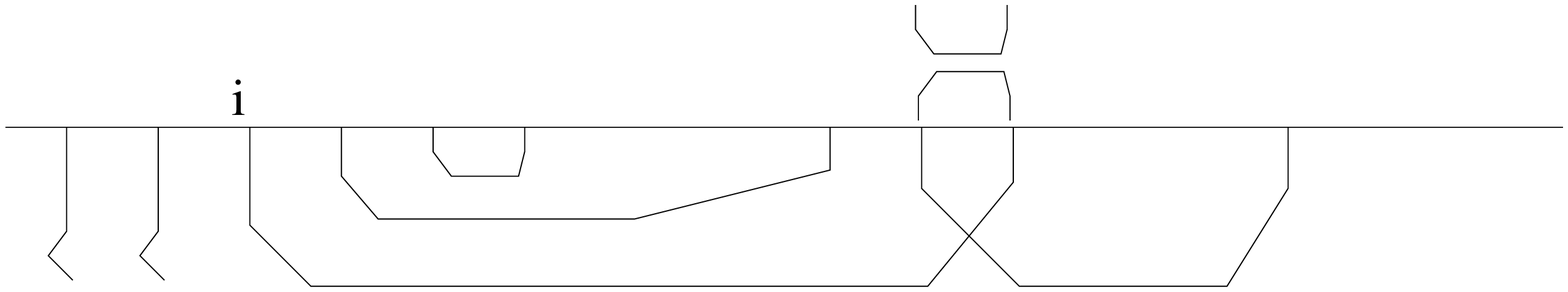}
= \;\;\;
\includegraphics[width=2in]{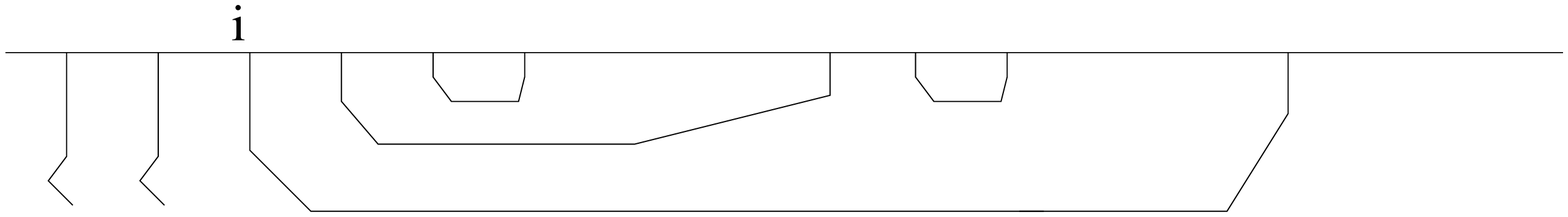}
\]
Here $c_k$ has one fewer link and one fewer link region so exclusivity
is preserved. 
\\
(iv)  $U_j$ touches two chains:
\[
\includegraphics[width=2in]{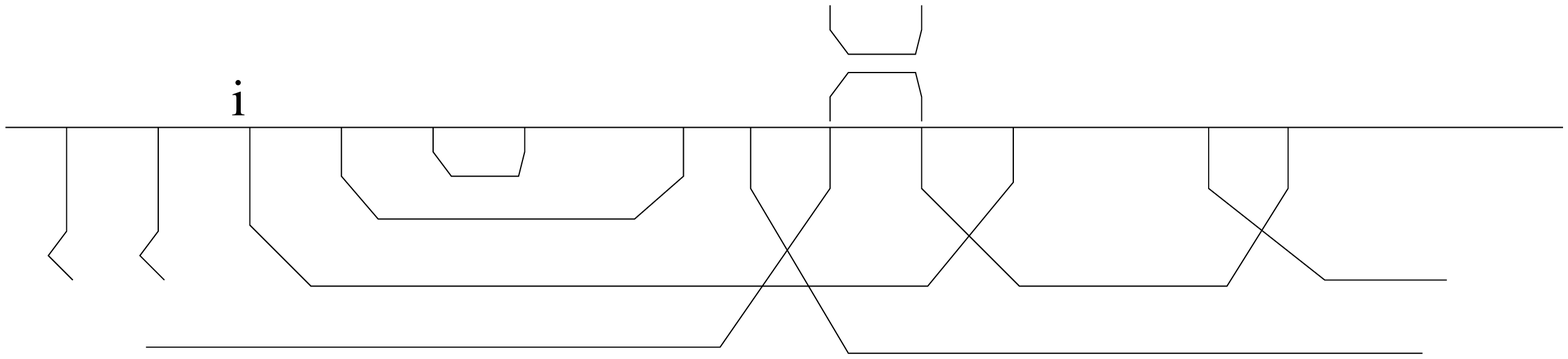}
= \;\;\;
\includegraphics[width=2in]{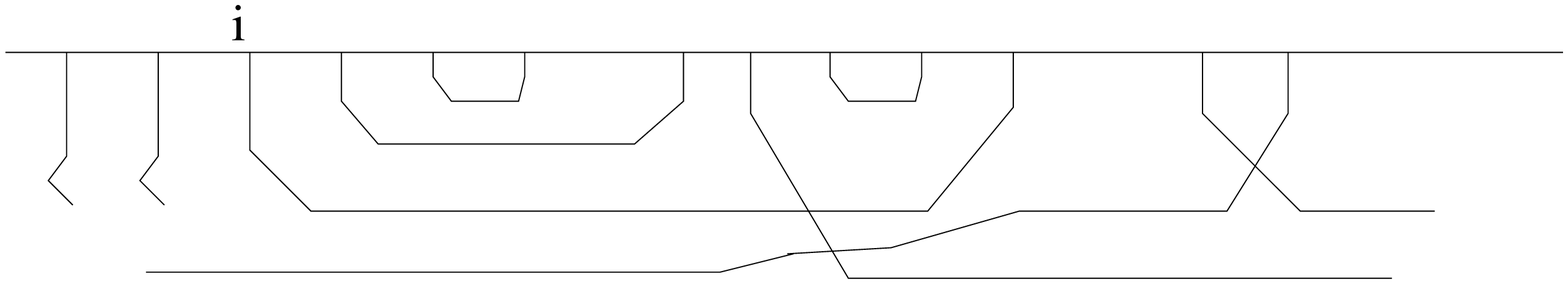}
\]
Here two adjacent link regions (for the two touched chains)
become combined as a single link
region.
Since they are adjacent there is no other link region between them in
$p$, and the combination does not affect exclusivity.
\qed
}


In Appendix~\ref{ss:colour} we give explicit examples of $\sigma_j$ actions
in figures with colour-coded chains. If you view in colour they may
help to reinforce the Lemma. 



\subsection{The Chain-basis Theorem} \label{ss:ChainbT}

%
%

For given $k$, and $i <m$, let
$\Jimi \defeq  k { J}_{\leq i-1}^i(m,m) $.
Recall the Brauer algebra $B_m = (k J(m,m), *)$ as defined in (\ref{eq:Bn}).

{\theo\label{mt3}
Fix $k$ and $\delta\in k$.  
Fix $i$ and consider $m \geq i$.  
Then  
$\Biim  \; = \; \Jimi$ 
as a $k$-module and hence as a subalgebra of $B_m $.}



\proof
It is clear that $(J_{i-1,i}^i,*)=(kS_i,*)=B_{i+1,i+1,i}$. 
So 
consider $k\Jimi$ with $ i<m$. 
Note that the Coxeter generating set of $\Biim$ 
(from (\ref{de:coxeter})) is in $\Jimi$. 
Thus it is enough to show $\Jimi  \; \subseteq \; \Biim$.

We  work by induction on the number of crossings 
$\chi_p$. 
\ignore{{
For this we establish a set $[p]''$ of pictures
of low height and low crossing number 
in  \S\ref{chi}. 
}}
Let 
\[
\Jic \; \defeq \; \{ p \in \Jim \; | \; \chi_p \leq c \}
\]
\newcommand{\Jico}{J^{i,0}_{\leq i-1}(m,m)}%
The base  of the induction concerns $p\in \Jico$. 
By the $Li$ condition this means that $p=1_i \otimes p'$ for some
non-crossing $p'$.
Then $p$ is generated by 
$\{U_j, j\geq i\}$ so it is clear that $p\in \Biim$. 
For the inductive step
we assume that $ \Jic \subseteq \Biim$ 
and require to show that this implies that 
 $\Jica \subseteq \Biim$.

\ignore{{ 
We then suppose \red{what?????}
that if $\pi(\sigma)\in B_{i+1,i+1,m}$ and  
$\pi(d)\in \Jimi$ with 
$d$ having $c$ crossings then also $\pi(\sigma|d)\in B_{i+1,i+1,m}$.
}}


Consider $p \in \Jica$.
Note that if the pair from vertex $i$ is non-crossing in $p$
(in the sense of (\ref{de:crossing}))
then it must
be $\{ i, i' \}$ and it is clear that $p$ lies in $\Biim$.
So suppose the pair from vertex $i$ is crossing.
If this pair crosses a pair from a lower vertex then
the pairs from some adjacent pair $j,j+1$ with $j<i$ cross.
(Note that the \lic\ condition implies that the pairs starting in
$\{ 1,...,i \}$ all pass out of this set --- to some set
$\{ p(1),...,p(i) \}$ say.
This permutation may be considered as generated from the
non-crossing one by the natural action of $S_{i-1}$ on $\{ 1,...,i \}$.
Since this action can be expressed in a reduced form in the Coxeter
generators it includes, if any crossing, a crossing of an adjacent pair.)
In this case $\sigma_{j} p$ has
crossing number $c$ and hence lies in $\Biim$ by assumption.
But then $\sigma_{j} \sigma_{j} p = p$ also lies in $\Biim$ and we
are done.

\begin{figure}
  \[
  \includegraphics[width=1.72in]{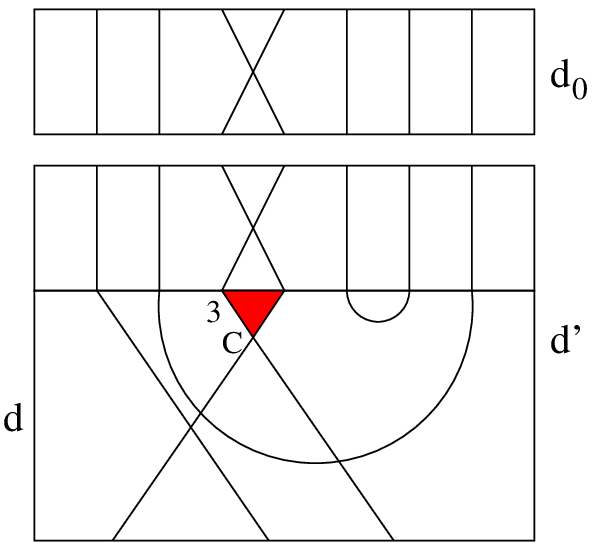}
  \]
\caption{\label{deco}
  The `decomposition' of a $d\in [p]''$
  for a $p\in J_{\leq 2}^3(7,3)$ into
  $d=d_0 |d'$,
with a $C^{ik}$-crossing 
  with $i=3$, $k=4$.
\redx{
We can strip off some crossings with low $i$ and expose ...}} 
\end{figure}

\begin{figure}
  \[
  \includegraphics[width=2.38in]{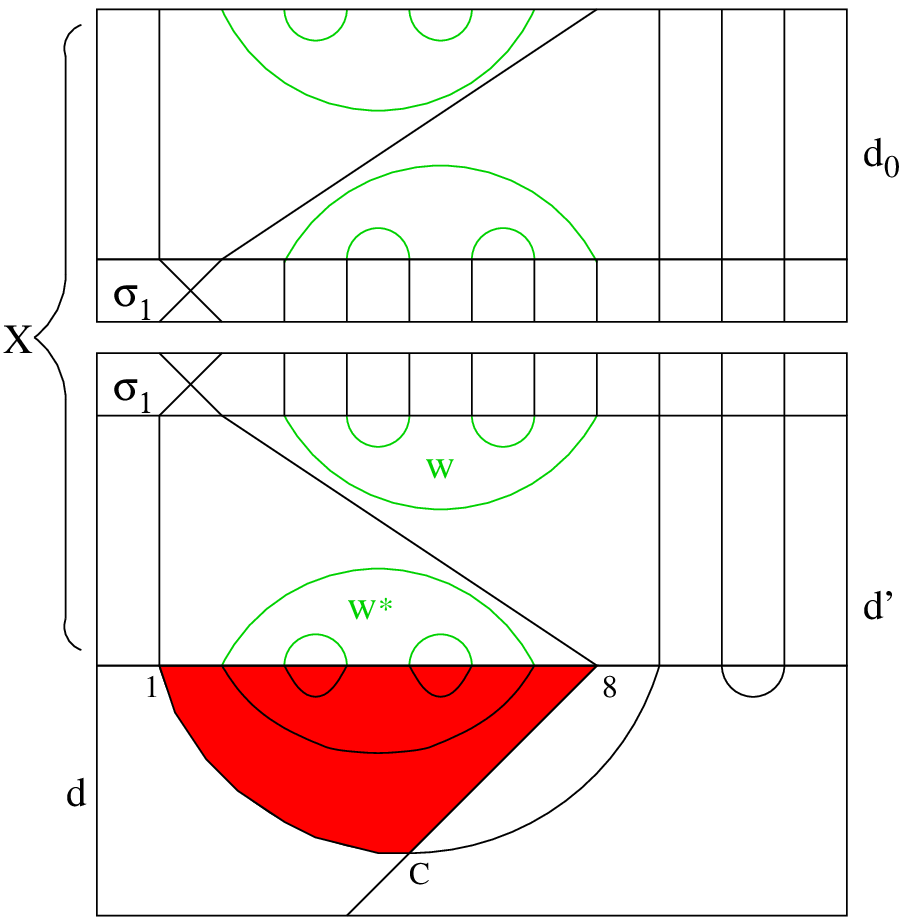}
  \]
\caption{\label{deco2}
  The `decomposition' $d\sim X|d$   
  or $d\sim d_0|d'$   
  with
a $C^{ik}$-crossing with 
  $i=1$, $k=8$.  
Here $d \in [p]$, $X \in [\tau^* \tau]$, $d_0 \in [\tau^*]$,
$d' \in [\tau p]$.
}
\end{figure}

It remains to consider the case in which the pair from $i$ crosses a
pair from some $k > i$.
Consider in particular the lowest such $k$. Call this crossing
$C^{ik}$.
For simple examples see Fig.\ref{deco} and \ref{deco2}. 
By Lem.\ref{lem:ppnc}
the vertices between $i$ and $k$ form a set $w$ of non-crossing
pairs.

\ignore{{
 \red{Now jump down!...}.


In Lem.\ref{lem:A} we assert that 
(on the `north edge' or else on the south edge) 
either there is a $j \leq i$ so
that the pairs involving $j$ and $k=j+1$ cross in $p$
(in the sense of (\ref{de:pcross})); 
or $j=i$ and $k$ cross for
some $k>i$ with intermediate vertices paired and non-crossing. 
Consider the `north' case, and 
call this crossing $C^{jk}$ (the south case is exactly analogous).
We will assert in Lemma~\ref{mj} that $\Jimi$ is closed under
left (and right) action of $\sigma_j$ with $j\leq i$.
Thus in the case $k=j+1$ 
we have $\sigma_j p \in \Jic$
(note that $\sigma_i p$ has crossing number at most $c+1+1$,
and that we have a cancellation with $C^{jk}$). 
But then 
$\sigma_j \sigma_j p \in \Biim$, so $p \in \Biim$.

The case $j=i$ and $k>i+1$ is more complicated, but we will 
argue analogously, passing
from $p$ to a product $\tau p$ such that $\tau p \in \Jic$;
and hence to $\tau' \tau p \in \Biim$ with $\tau' \tau p =p$. 
First note by Lemma~\ref{lem:ppnc} that between $i$ and $k$ there is only a
collection $w$ of non-crossing pairs in $p$.


\red{...Jump to here!...}
}}

\newcommand{\taub}{\overline{\tau}} 
\newcommand{\taubb}{\overline{\tau^*}} 

Consider the partition
$
\mu = 1_i \otimes w^* \otimes 1 \otimes w \otimes 1_{m-k} 
$
in $B_m$.
See the middle layer in Fig.\ref{deco2} for an example.
(Since pairs from $j<i$ do not cross the pair from $i$ here,
pairs from $j<i$ do not need to be tracked closely.
Our example has $i=1$ to preclude such clutter.)
%
Note that $\mu$
lies in $\Biim$.
But then $\tau = \sigma_i \mu$ lies in $\Biim$ 
{\em and} has a single crossing, which is a crossing of the pairs from
$i'$ and $k'$.
Hence $\tau$ cancels $C^{ik}$ in the partition given by the  product $\tau p$. 
(In general 
there are some loops in $\tau p$, but each loop includes pairs from 
$\tau$, so there is a `garden path' modification (\ref{de:gardenpath})
$\taub$ of $\tau$, also in $\Biim$, 
so that $\taub p$ gives the same partition as $\tau p$ but without loops.)

\begin{figure}
  \[
  \includegraphics[width=2.38in]{xfig/ind2x}
  \qquad \includegraphics[width=2.21in]{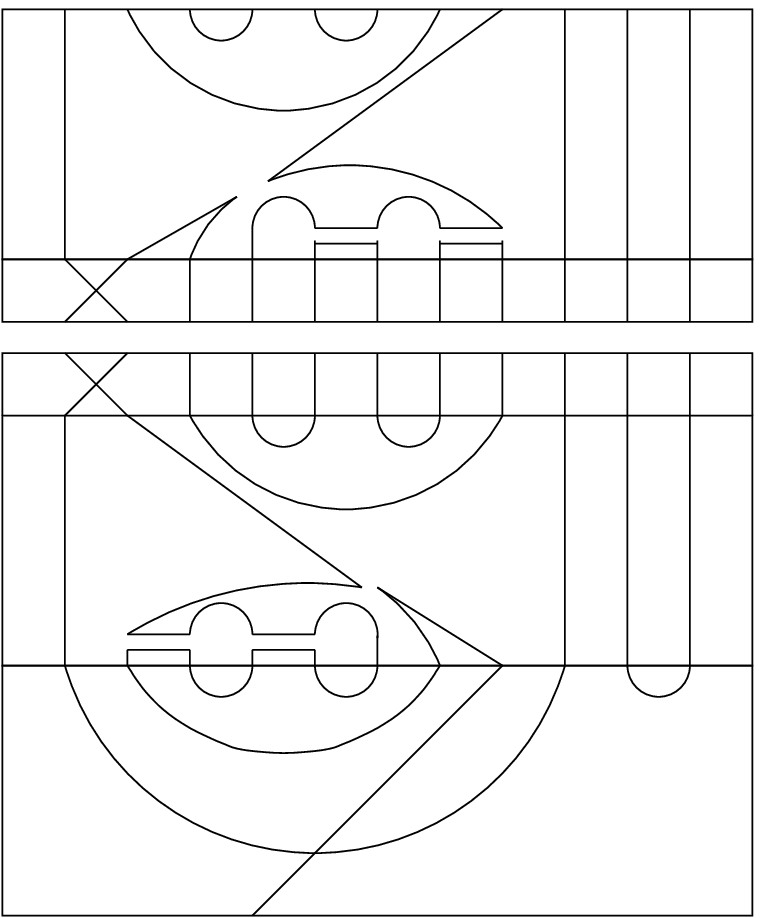}
  \]
\caption{\label{deco2gp}
(Right)  A picture  of the same partition modified  
  such that there
  are no loops in it.   
}
\end{figure}

By Lem.\ref{mj} we have that $\tau p $ is in $\Jimi$.
By the cancellation
we have $\chi_{\tau p} = \chi_p + \chi_{\tau} -2 = c+1+1-2$.
Thus $\tau p$ lies in $\Jic$ and hence,
by the inductive assumption, in $\Biim$.
Since $\tau^*$ also lies in $\Biim$ we have that
$\tau^* \tau p$ lies in $\Biim$. But
$\tau^* \tau = 1_i \otimes w \otimes w^* \otimes 1_{m-k+1}$
so
$\tau^* \tau p = \delta^{2 |w|} p$
--- see Fig.\ref{deco2} for an illustration.
Unless $\delta=0$ this directly implies $p \in \Biim$ completing the
inductive step. 
(Alternatively there is an analogue of $\taub$ so that  
$\taubb \taub p =p$, which completes the inductive step in general
--- see Fig.\ref{deco2gp} for an illustration.)
\qed

\newpage

\subsection{The category $\BBl{l-1}^l$}

Let  $u$, $u^*$ be the unique elements in $J(0,2), J(2,0)$,
respectively;
and $1_1 \in J(1,1)$. 

{\lem \label{lem:modulex}
  Let $q \in J(m,n)$. 
  Suppose $p= q \otimes u$ or $p=q\otimes u^*$ or $p=q\otimes 1_1$.
  Then $p$ is $Li$-chain if and only if $q$ is.
}
\proof Note that such extensions affect neither the $ht$ nor the chain
structure. \qed

{\theo There is a subcategory of $\BBl{l-1}$ given by 
\[
\BBl{l-1}^l = (\N_0 , k \JJ^{\,l}_{\leq l-1}(m,n), *)
\]

}
{\proof
  (There are different ways to approach this.)
Let us consider the Hom set bases $J_{\leq i-1}^i(m,n)$
and $J_{\leq i-1}^i(n,o)$. 
It is enough to check that  the restricted composition on 
$ J_{\leq i-1}^i(m,n)\times J_{\leq  i-1}^i(n,o)$
lies in
$J_{\leq  i-1}^i(m,o)$.
We will do this by
embedding the Hom sets into $J_{i-1,r}^i=B_{i+1,i+1,r}$ with
$r=\max(m,n,o)$
as follows.
The embedding uses monoidal composition with 
powers of $u$, $u^*$,
as illustrated in Fig.\ref{fig:cat}.
%
%
\redx{delete?:Since the composition in the algebra is closed by Th.\ref{mt3}, 
this will ensure that the category composition is closed.}
%
%
For $p\in J(s,t)$ and $l\in {\mathbb N}$
congruent to $s$ mod.2
define 
\[
\overline{p}^l
=p\otimes u^{\otimes\,\theta(l-s)}\otimes
                 (u^*)^{\otimes\,\theta(l-t)}\ ,
\]
where $\theta(x)=\max(x/2,0)$.

Now let $p\in J_{\leq i-1}^i(m,n)$, $ p'\in J_{\leq  i-1}^i(n,o)$.
Thus by Lem.\ref{lem:modulex} and Theorem~\ref{mt3}
$\overline{p}^r, \overline{p'}^{r}\in B_{i+1,i+1,r}$.
\redx{
But the $Li$-chain property is also 
preserved under monoidal composition with an element from
$J_{-1}(-,-)$ on the right: it corresponds to adding a noncrossing partition 
in the disk intervals $[m,n']$ of $p$ and $[n,o']$ of $p'$.   } 

By the same token
there is (provided $\delta\in k^*$) a unique $p''\in J_{\leq i-1}^i(m,o)$ such that
\beq \label{eq:catify}
\overline{p}^r*\overline{p'}^{r}   \propto   \overline{p''}^r\ .
\eq
Thus by
Theorem~\ref{mt3} and 
Lem.\ref{lem:modulex} again we have $p * p' \in k J^i_{\leq i-1} (m,o)$
as required
(at least provided that $\delta\in k^*$).

Note that the product on the left in (\ref{eq:catify}) can be zero
if $\delta=0$.
Note however that 
we can eliminate  
loops from the product by
a mild generalisation of $\overline{p}^r  $ using
suitable garden paths.
\redx{
in the concatenation in computing 
the product ($*\to *^g$) in the algebra and define $p*p'$
by the unique element $p''\in kJ_{\leq i-1}^i(m,o)$ satisfying 
$\overline{p}^r*^g\overline{p'}^{r}=\overline{p''}^r$. 
}
\redx{
Note that any number $r+2w, w\in {\mathbb N}$
can be used above without changing the definition of the product. Using this 
observation one can explicitly verify associativity by a direct 
calculation.}
\qed

\begin{figure}
\begin{center}
  \includegraphics[width=8cm]{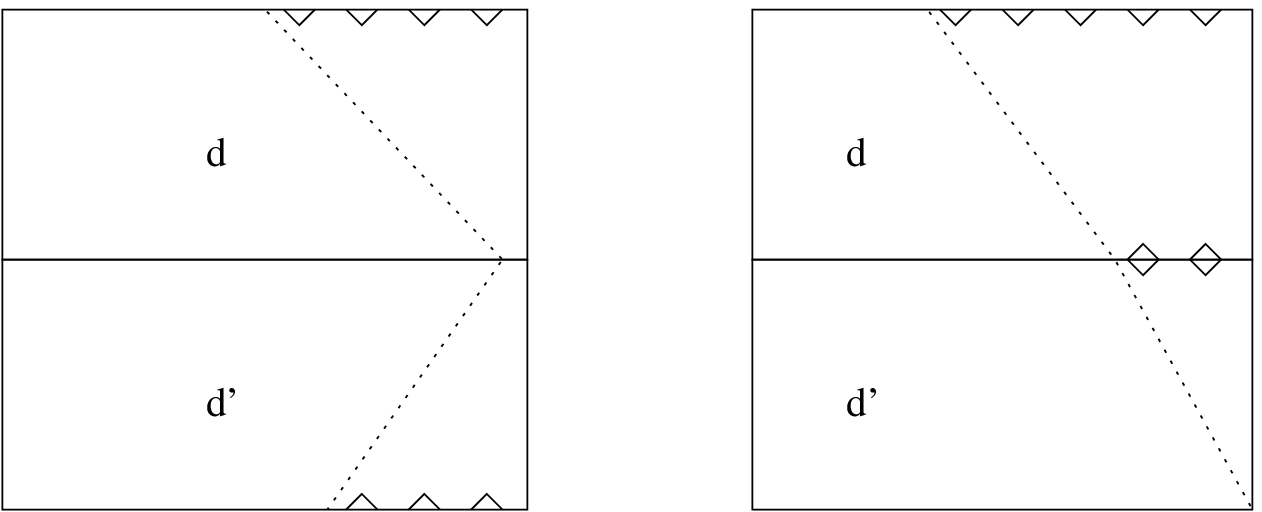}
  \caption{
\label{fig:cat}
    Pictures of $\overline{p}^r$ and $\overline{p'}^r$
    from partitions 
$p\in J(m,n)$ and $ p'\in J(n,q)$ with $n>m>o$ (left) and $m>n>o$ (right).}
\end{center}
\end{figure}
}






\newpage

\section{The blob isomorphism Theorem}  \label{ss:main01} 




We start with some notation. 

Given a pair partition $p\in J(m,n)$
denote by $p^+$ the relabeling obtained by adding +1 to each label. 
That is, 
$p^+$
is a pair partition of the set
$
\{2,3,\dots m+1,2',3',\dots (n+1)'\}
$.
Write $p^-$ for the partition in $J(m,n)$ obtained by the 
inverse relabeling, changing each element of the
shifted underlying set $\{ 2,3,...,(n+1)' \}$  by -1.

\newcommand{\ls}{l} 

For $x $ a chain in the form
\[
x = \{\{ 1,r_1+1 \}, \{ l_1+1,r_2 +1 \},
      \dots  \{ l_{\ls-1}+1,r_{\ls}+1 \}, \{ l_{\ls},1' \} \}
\]
(recall for example that every $p\in J^1_0(m,n)$ is $p=c_1 \cup p'$ where $c_1$
is such a chain)
let
\[
\bar{x} = \{\{l_1,r_1\},\{l_2,r_2 \}, ... ,\{l_{\ls},r_{\ls} \}\}
\]
Schematically this is
\beq \label{eq:schembx}
\includegraphics[width=3cm]{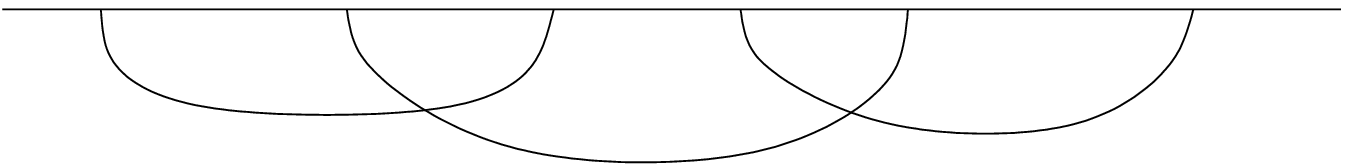}
\mapsto
\includegraphics[width=3cm]{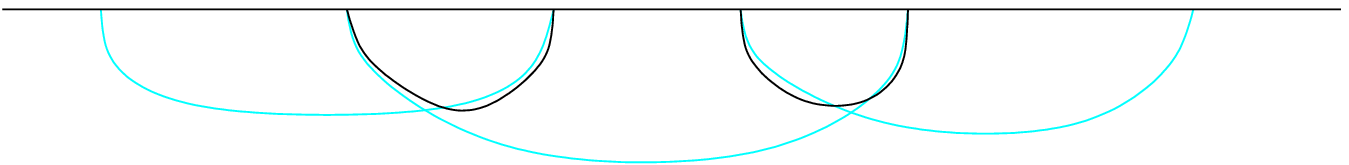}
\mapsto
\includegraphics[width=3cm]{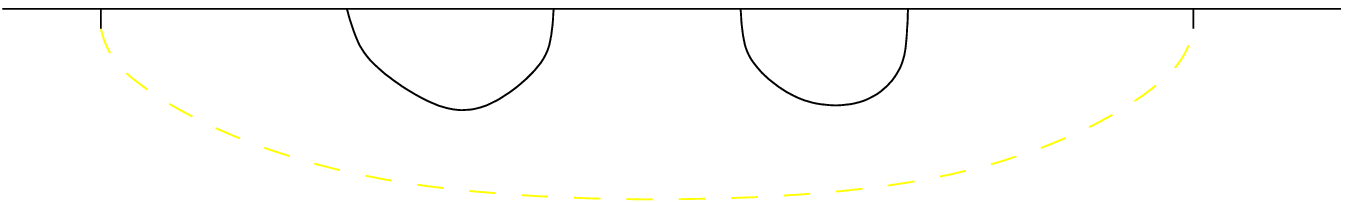}
\eq
(at least up to ambient isotopy).
Note that the sequence  
$l_1,r_1,l_2,r_2,\dots,l_{\ls},r_{\ls}$ is a subsequence of
$1,2,\dots,m,n',(n-1)',\dots,1'$.
That is $l_1 < r_1 < l_2 < r_2 <...$ in the disk order.

\newcommand{\xx}{{\mathsf x}} 

Let $x\mapsto \xx(x) $ be the map that inverts $x\mapsto \bar{x}$.
Note that this simply reverses the arrows in the schematic.

\ignore{{

We will first define a bijection 
$
\Psi_{}: \bBd(m,n)\to  \Jd^{1}_0(m+1,n+1) .
$ 
Then, using $\Psi_{}$,  
we will construct another family of maps $\Phi_{}$, which 
satisfies  
$$
\Phi_{}((p_1,s_1) \bx (p_2,s_2)) = \Phi_{}(p_1,s_1)*\Phi_{}(p_2,s_2)
$$
with $(p_1,s_1)\in \bB(m,n)$ and $(p_2,s_2)\in \bB(n,q)$,
provided that   
$$
\delta' = (1+\delta)/2 .
$$ 

In the latter case $\Phi_{m,m}$ will be shown to be an isomorphism of the algebras 
$J^{\bullet(\delta,(\delta+1)/2)}(m)$ and $J_0^{1,\delta}(m+1,m+1)$
(Theorem~\ref{th:main1}). 
\qq{what algebras?? is there more to show??}
\qq{what about the categories??}

}}
\subsection{The initial set map
  $\Psi: \bB(m,n)\to J^1_0(m+1,n+1)$ \label{setmap}} \label{smap}
\begin{figure}
\begin{center}
\includegraphics[width=4.5cm]{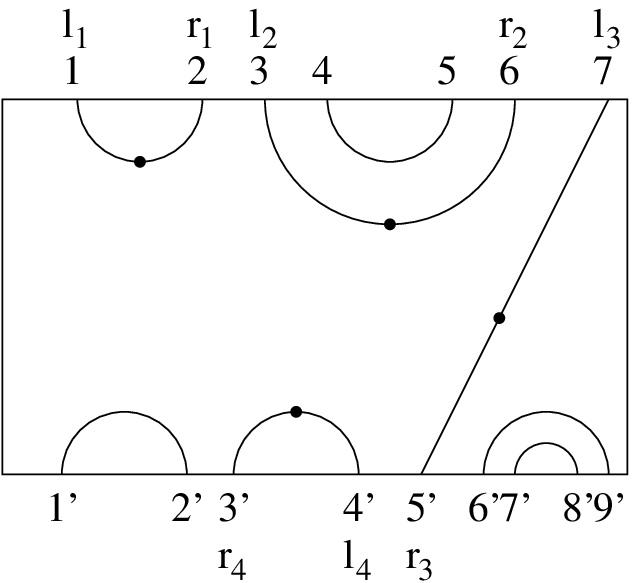}$\quad$
\includegraphics[width=4.55cm]{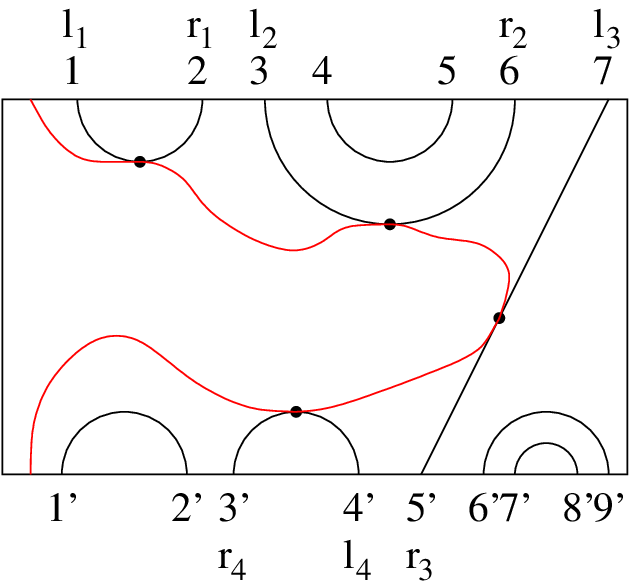} \quad
\includegraphics[width=4.55cm]{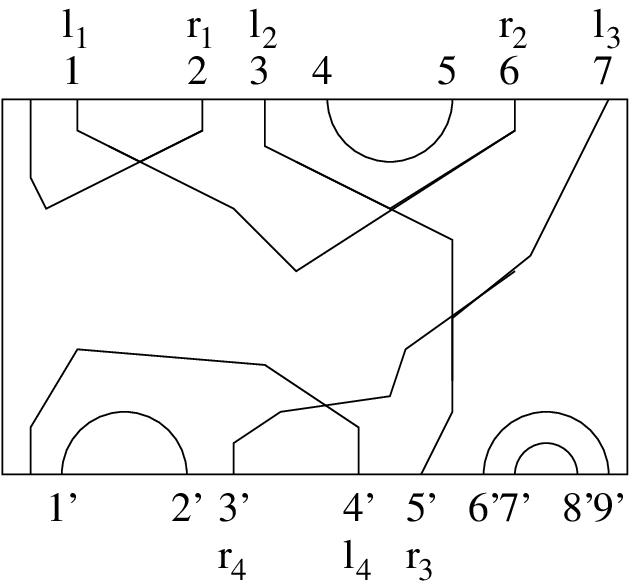}
\end{center}
\caption{
  (a) A picture $d\in [p]'$ for $ p\in J^\bullet(7,9)$.
  (b) Insert a path $c$ into $d$.  
(c) Cross at blobs.
  \label{fig:fug01}} 
\end{figure}
\begin{figure}
\begin{center}\includegraphics[width=7cm]{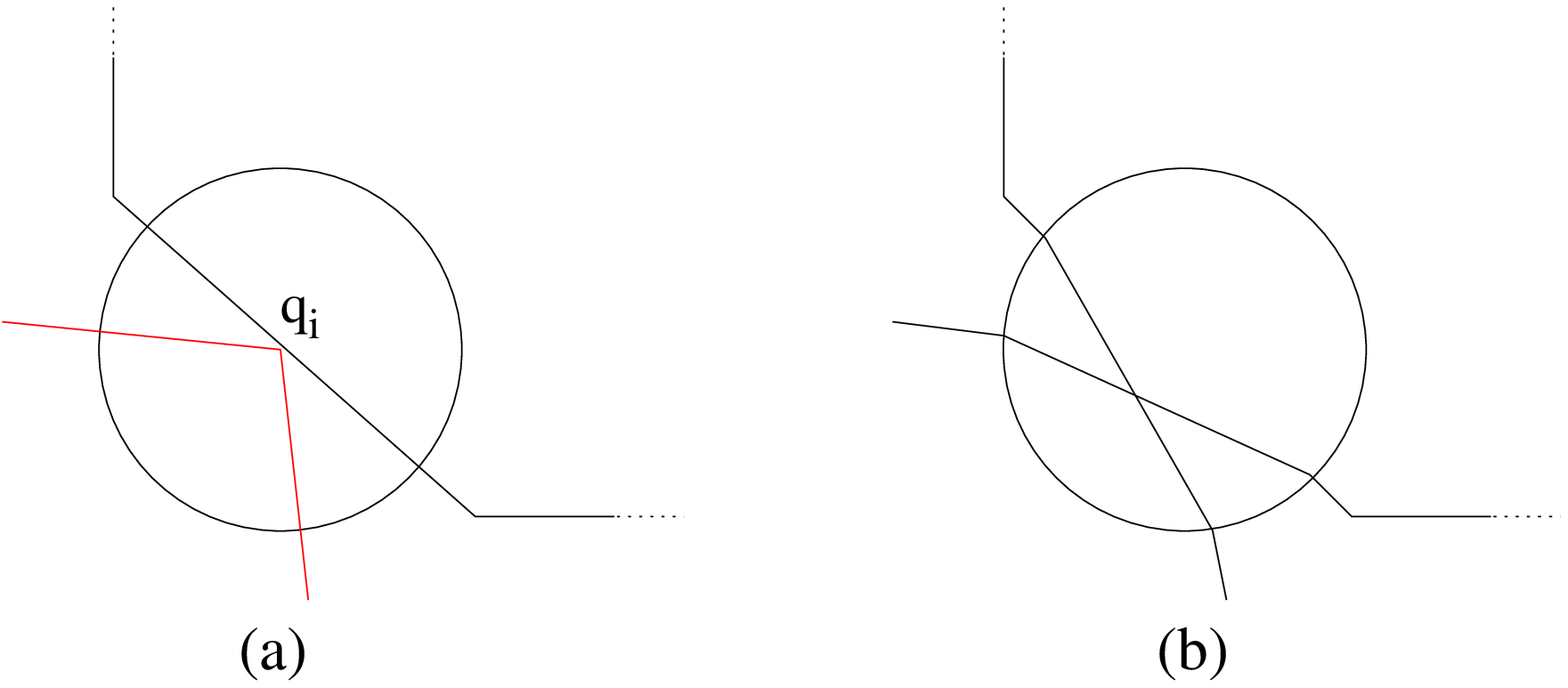}\end{center}
\caption{\label{fig:nec}(a): a small neighbourhood of $q_i$, with the
  path $c$ in red.  
(b) reconnection of the line segments in the small neighbourhood
  according to the map 
$\Psi$.}
\end{figure}

Consider $(p,s) \in \bB(m,n)$, as in \S\ref{ss:blob1}.  
See Fig.\ref{fig:fug01}(a) for a picture $d$ of an example,
with $d \in [p]'$.
(Note that here such pictures are unique up to {\em ambient} isotopy.) 
We can use such a picture $d$ of $(p,s)$ to describe a modification
$\Psi(p,s) \in J^1_0(m+1,n+1)$
with the following steps: \\
(1)
taking $d \in [p]'$, add new boundary points $0,0'$ and draw a path between them
 that touches the lines of $d$ only at each of the blobs,
 as in Fig.\ref{fig:fug01}(b).
\\ (2)
noting that each blob now has four line segments incident,
replace this with two crossing lines
as in Fig.\ref{fig:nec},
to obtain $d'$ (Fig.\ref{fig:fug01}(c)).
\\ (3)
Set $\Psi(p,s) = \pi(d')$.

\medskip

Note that (1) is possible by the left-exposed property of blobs.

{\lem{ \label{lem:Psi}
(I)    The procedure $\Psi$ gives a map
$$
\Psi : \bB(m,n)\to J^1_0(m+1,n+1)
$$
(II) $\Psi$ is a bijection.
\\ (III) For $(p,s) \in \bB(m,n)$, writing
$s = \{\{ l_1 , r_1 \}, \{l_2,r_2 \},...\}$
where $l_k < r_k$ and this is the $k$-th pair,
both orderings with respect to the disk-order,
then
\[
\Psi(p,s) = \xx(s) \cup (p\setminus s)^+
\]
where
\[
\xx(s) = \{\{ 1,r_1+1 \}, \{ l_1+1,r_2 +1 \},
      \dots  \{ l_{|s|-1}+1,r_s+1 \}, \{ l_{|s|},1' \} \}
\]
(IV) The map $\bar\Psi : J^1_0(m+1,n+1)\to \bB(m,n)$
given by $x \cup p \mapsto \bar{x} \cup p^-$ is inverse to $\Psi$.
\\ (V) If $\Psi(X) = Y\otimes u$ then $X$ takes the form
$X = X' \otimes u$ (in the obvious sense).
\qed
}}
\proof
(I)
It will be clear that the procedure gives an element of
$ J^1_0(m+1,n+1) $. Thus it remains to show that this is independent of
the choice of $d \in [p]'$ representing $(p,s) \in \bB(m,n)$.
This follows since the construction of picture $d'$ in $\Psi$,
both in the choice of $d$ and 
at steps (1-2), is unique up to {\em ambient} isotopy and produces a {\em
  canonical}  picture.
(It also follows that we can recast the  procedure at the original
set-theoretic level.
See (III).)

\newcommand{\appfug}{
\begin{figure}
\begin{center}
\includegraphics[width=4.5cm]{\figd/xfig/blob1}$\quad$
\includegraphics[width=4.55cm]{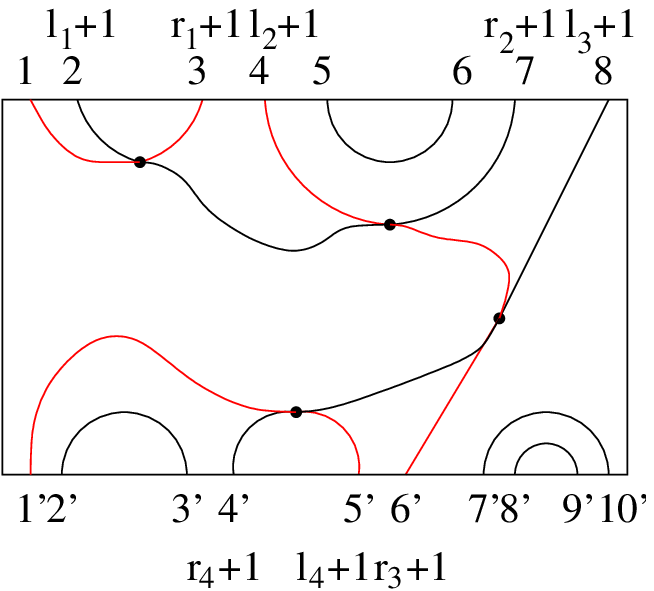}$\quad$
\includegraphics[width=4.5cm]{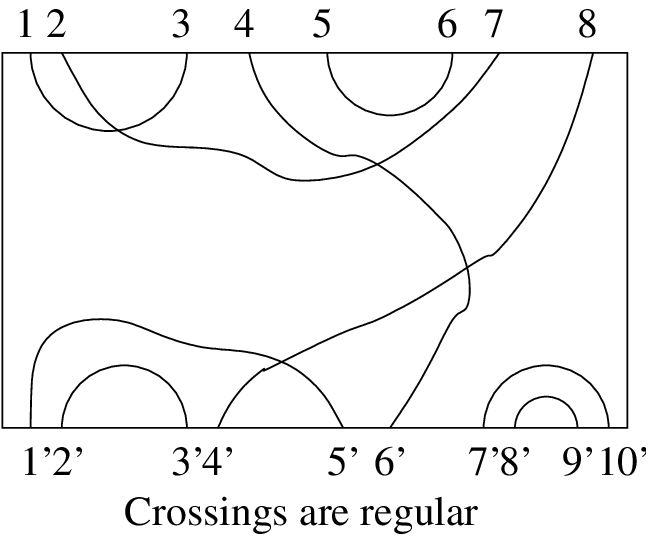}
\end{center}
\caption{(a)
  insertion of the path $c$ into $d\in [p]', p\in J^\bullet(7,9)$,  (b)
  renumbering the boundary vertices and modifications of the small
  neighbourhoods    of the blobs $q_i$,
  (b)
  the resulting object is a picture of $\Psi(p,s)\in J_0^1(8,10)$. 
\label{fig:fug1}} 
\end{figure}

Given $(p,s) \in \bB(m,n)$ 
we write the pairs in $s$ as 
$\{l_k,r_k\}, k=1,...,|s|$,
where $l_k < r_k$ and this is the $k$-th pair,
both orderings with respect to the disk-order. 

See Fig.\ref{fig:fug1}(a) for an example: the set of lines with blobs 
project to the set of pairs, which correspond to
$\{ l_1 =1, r_1 =2 \}$,
$\{ l_2 =3, r_2 =6 \}$,  
$l_3=7$, $r_3=5'$, $l_4=4'$, $r_4=3'$.

Note that the sequence  
$l_1,r_1,l_2,r_2,\dots,l_{|s|},r_{|s|}$ is a subsequence of
$1,2,\dots,m,n',(n-1)',\dots,1'$.


\begin{lem} \label{lem:fugue1}
There is a well-defined  map 
$$
\Psi : \bB(m,n)\to J^1_0(m+1,n+1)
$$ 
given by
\beq\label{setmap}
 \Psi(p,s)  =  
    x \cup (p\backslash s)^+\ 
\eeq
where, with $s$ as above, 
\[
x = \{\{ 1,r_1+1 \}, \{ l_1+1,r_2 +1 \},
      \dots  \{ l_{|s|-1}+1,r_s+1 \}, \{ l_{|s|},1' \} \}
\]
%
\end{lem}
\proof
We will verify $\Psi(p,s)\in J_0^1(m+1,n+1)$ by constructing the
picture of it. 
Let $d\in [p]'$.
Choose  
$|s|$ points, which are not from the $r$-set of the lines corresponding to the
left-exposed pairs in $s$. 
Denote the chosen $|s|$ points by $q_k,k=1,2,\dots,|s|$ according to the
disk-ordering of the left points $l_i$ of the pairs. 

Denote the middle point between the top left (bottom left) corner of
$R$ and the most left northern (southern) endpoint of a line 
by $0$ ($0'$, respectively). 

We continue the proof of Lemma~\ref{lem:fugue1} 
after the proof of the following lemma. 


\lem{\label{path}
Consider $d$ as  
above. 
There exists a connected path $c$ from $0$ to $0'$
  without intersecting any line of the picture,  
which connects $0,q_1,q_2,\dots,q_{|s|},0'$ in this order with the following 
property. At a neighbourhood of the point $q_i$, $c$ consists of two consecutive
straight segments, whose common point is $q_i$, see Fig~\ref{fig:nec}(a) and 
also Fig.~\ref{fig:fug1}(a) for an example of the full path $c$ in
red.
}
 
\proof{
The points $0,q_1,q_2,\dots,q_{|s|},0'$ are in the boundary of the
${\cal A}_0$, which is a connected set, and they are  
oriented naturally as $0,q_1,q_2,\dots,q_{|s|},0'$. Therefore, the
line $c$ can be drawn along the boundary of ${\cal A}_0$, but
slightly inwards, except at the points.
\qed}

\medskip


Now we  
conclude the proof of Lemma \ref{lem:fugue1}. 
Consider the union of $d$ and a path $c$ as in Lemma \ref{path}. 
We relabel the vertices of $d$ in $\partial R$ 
as $2,3,\dots m+1$ ($2',3',\dots, (n+1)'$) in 
the northern (southern) edge and  
the beginning
and endpoint of $c$ by $1$ and $1'$, respectively. 
Take a small neighbourhood of $q_k$ 
($k\in\{1,2,\dots,|s|\}$ and change it slightly such that the connection 
of segments change according to the formula for the map $\Psi$, 
see Fig.~\ref{fig:nec}b. The resulting object has height 
$0$ since the path $c$ is in ${\cal A}_0$ and all 
intersections are in $c$. It is also 
left-simple by construction. Hence, the resulting object is a picture for 
$\Psi(\{p,s\})\in J_0^1(m+1,n+1)$ (see figure~\ref{fig:fug1}).
\qed

\medskip

}


\redx{THIS should just follow quickly from the chain construction!}

(II)
Note also that steps (1) and (2) are reversible.
By
(\ref{pa:once})
we can pass from $p \in J^1_0(m+1,n+1) $
to a canonical picture.
Thus $\Psi$ is invertible.

(III) This is the promised formal version of the picture
manipulation.
Confer Lemma~\ref{le:ppnc}, (\ref{eq:schembx}) and Fig.\ref{fig:nec}. 

(IV) This is a disjoint combination of manifest inverses.

(V) Follows from (IV).
  \qed

\subsection{An algebra isomorphism}
\redx{WHere is the theorem?? Is it one of the Lemmas here? 17? 19?}

\mdef
There is a map $\Theta$ on the $\bb_n$ generator set $U^e$ 
(cf. Theorem~\ref{th:blobgen1}) to $kJ_0^1(n+1,n+1)$ given by:
\beq \label{eq:Theta} 
U_i \mapsto U_{i+1}  \; = \Psi(U_i),
\qquad\qquad
e \mapsto \frac{1}{2} (1+\sigma_1 )  \; = \frac{1}{2}(\Psi(1)+\Psi(e))
\eq
(note the uses of notation $U_i$ distinguished by context).

{\theo{ \label{th:Theta}
Fix a commutative ring $k$ and $\delta \in k$, and set $\delta'=\frac{\delta+1}{2}$.
Map $\Theta$ extends to an algebra isomorphism
$\Theta: \bb^{}_n
     \rightarrow {J}^{1}_{0,n+1}$
{
(with  ${J}^{1}_{0,n+1}$ as defined in \S\ref{ss:ChainbT}).  
}
}}

\proof{
To verify that $\Theta$ extends to an algebra homomorphism
it is enough, by Theorem~\ref{th:blobgen1},
to check that the images obey 
the relations $\tau'$ from that Theorem.  
The relations on the $\{ U_i \}$ are all clear.
The image of $ee=e$ is clear. Next $U_1 e U_1 = \delta' U_1$:
\[
U_2 \frac{1}{2}(1+\sigma_1 ) U_2 = \frac{\delta+1}{2} U_2
\]
The remaining commutation relations will be clear.

This homomorphism is surjective on
$B_{2,2,n+1}$   
since it hits the generators. On the other hand
$B_{2,2,n+1} = \Jonip$    
by
Theorem~\ref{mt3}. But $\bb_n$ and $\Jonip$ are isomorphic as vector
spaces by
Lem.\ref{lem:Psi}.
\qed
}


\subsection{Category version}

\ignore{{
Let
  $
\Psi:\bBd(m,n)\to \Jd_0^{1}(m+1,n+1)
$ 
\qq{I think it is $\Psi$ not $\Phi$ here???}
\qq{domain/range defined??}
\qq{do we need the $k$??}
the set map defined in section \ref{setmap}. 
}}

\mdef
For $m,n \in \N_0$ define 
\begin{eqnarray*}
\Phi_{}:
\bBd(m,n) & \to & k\Jd_0^1{}(m+1,n+1)\\
(p,s)&\mapsto&2^{-|s|} \sum_{z\subseteq s} \Psi(p,z)
\end{eqnarray*}
and extend $k$-linearly.
We write $\Phi_{m,n}$ to indicate cases of $\Phi$ where convenient.


\mdef \label{pa:wds}
Recall
(Thm.\ref{th:blobgen1}(II))
that every (blob-)partition $(p,s)$ can be expressed as a product of
generators, and write $\#(p,s)$ for the length (the minimum number of
factors in such a word).
\\
Remark: it is not in general easy to determine length from $(p,s)$
itself. It is easy to determine the number of factors $e$, since this
is just $|s|$.
Note also that disk order induces a total order on $s$.
For example if $\{a,b\}$ is the last pair in $s$ with $a<b$ in disk
order then the image under $\Psi$ includes $\{ a+1, 1' \}$.


\begin{lem} \label{le:PhiTheta}
The map $\Phi$ restricts to an algebra isomorphism
$\bb_n(\delta,(\delta+1)/2) \rightarrow J_{0,n+1}^{1}$.
That is, $\Phi = \Theta$ when $m=n$.
\end{lem}



\noindent Proof:
We work by induction on word length, as in (\ref{pa:wds}).  
The base case is trivial.


Consider $(p,s)$ with $|s| =j>0$ and consider
the same partition with one fewer blobs, $(p,s\setminus a)$ say
(here $a$ is one of the pairs in $s$). Then
\beq \label{eq:Phex}
\Phi(p,s) = \frac{1}{2^{|s|}} \left( \sum_{z\subseteq s\setminus a} \Psi(p,z)
                 +  \sum_{z\subseteq s\setminus a} \Psi(p,z\cup\{a\}) \right)
\eq


We assume $\Phi(w) = \Theta(w)$ for words $w$ of length $\leq l$, and aim to
show for words of length $l+1$.
Every such word may be written in the form
$wx$ where $w$ has length $l$ and $x$ is a generator.
First consider $x=e$.
If $(p,s) = we$ in $\bb_n$
then $s$ has an element $a=\{1',j\}$ for
some $j$, and $w = (p,s\setminus a)$.
Note that in this case $\Psi(p,z\cup\{a\})=\Psi(p,z)\Psi(e)$.
So in this case 
(\ref{eq:Phex}) becomes  
\[
\Phi(p,s) = \frac{1}{2} \left( \Phi(p,s\setminus a)
                 + \sum_{z\subseteq s\setminus a} \Psi(p,z)\Psi(e) \right)
=  \Phi(p,s\setminus a) \frac{1}{2} \left( \Psi(1)
                 + \Psi(e) \right)
= \Phi(p,s\setminus a) \Theta(e)
\]
by (\ref{eq:Theta}).
By inductive assumption, $\Phi(w)=\Theta(w)$,
so we have made the inductive step in this case.

Similarly if $(p,s) = wU_i$ for some $i$, with $w U_i$ longer than
$w = (p_w, s_w)$, then either $w \mapsto wU_i$ does not change $s$,
$p = p_w U_i$ and 
\[
\Phi(p,s) =  \frac{1}{2^{|s|}}  \sum_{z\subseteq s} \Psi(p,z)
=  \frac{1}{2^{|s|}}  \sum_{z\subseteq s} \Psi(p_w U_i,z)
=  \frac{1}{2^{|s|}}  \sum_{z\subseteq s} \Psi(p_w ,z) \Psi (U_i)
=  \Phi(p_w,s) \Theta(U_i)
\]
or  $w \mapsto wU_i$ changes a single element of $s$ by
$\{i',j \} \mapsto \{k ,j \}$
(where $k$ is the element in a pair with $i+1'$ in $p_w$)
and
similarly
$\Psi(p,s) = \Psi(p_w, s_w) \Theta(U_i)$.
This completes the inductive step.
\qed


\ignore{{

\noindent
Proof:
We work by induction on length. That is, we suppose
$\Phi(w) = \Theta(w)$ for partitions $w$ of length less than $n$, and
aim to deduce equality for $wx$ for $x$ a generator such that
$wx$ has length $n$.

The base case can be $n=1$, which is easy to check.
We must consider cases $x=e$ and $x=U_i$. 

First consider $x=e$.
We require to compare $\Phi(we)$ with $\Phi(w)\Phi(e)$.
We may write
$$
w = (p,s) = (\{\{ 1', j \}, \{ 2', k \},...\},...)
$$
where we have omitted pairs that do not intersect $\{ 1', 2' \}$.
In particular we may suppose neither vertex `lies in' $s$, since
if $1'$ lies in $s$ then $we$ is not longer than $w$; and then
$2'$ cannot lie in $s$ by definition of blob-partition.
Meanwhile
$we = (p',s') = (\{\{ 1', j \}, \{ 2', k \},...\},\{\{1', j \},...\})$.
We have
\[
\Phi(p,s) = \sum_{z\subseteq s} \Psi(p,z)
  = \{ \{1',a \}, \{ 2', j+1 \}, ...  \}  +...
\]

Before completing the inductive step note that if $w$ has no blobs
then $\Phi(w) = \Theta(w)$ by the construction and Lemma~\ref{}.

Now consider $x=U_i$.

\[ \]

}}

\medskip

\begin{lem}\label{otempora}
Let $(p,s)\in J^\bullet(m,n)$ and $p_0\in J_{-1}(m_0,n_0)$.
Then we have the following (confer \paref{de:2.18} and \paref{de:moncat}).
\\
(i)
the operation 
$(p,s) \mapsto (p\otimes p_0,s)$ defines an injection
$f_{-\otimes p_0} : J^\bullet(m,n) \hookrightarrow
J^\bullet(m+m_0,n+n_0)$,
and similarly on the corresponding vector spaces.
\\ 
(ii)
\beq
\Phi_{m+m_0,n+n_0}(p\otimes p_0,s)
  =\Phi_{m,n}(p,s)\otimes p_0 .  
\eeq
\redx{+ meaning? can I just delete all these +'s?}
\end{lem}
{\proof  
(i) $(p\otimes p_0,s)\in J^\bullet(m+m_0,n+n_0)$
since the height of
a picture is not increased by concatenating a non-crossing piece to the right and 
$s\subseteq S_p^L$ implies
$s\subseteq S_{p\otimes p_0}^L$.
Injectivity is clear. 
(ii) Now  
\beq
\Phi(p\otimes p_0,s)=
2^{-|s|}\sum_{z\subseteq s} \Psi(p\otimes p_0,z)=
2^{-|s|}\sum_{z\subseteq s}\Psi(p,z)\otimes p_0
=\Phi(p,s)\otimes p_0 \label{sq}\ ,
\eeq
where the second equality holds by construction
since
the $p_0$ part has no blobs.
\qed
}

\medskip


{\theo
The collection of maps $\Phi = \Phi_{m,n}$ yields an equivalence between category $\bb$
(with $\delta' = (1+\delta)/2$)
and category ${\cal J}^1_0$
(the map on objects is $\Phi(n)=n+1$).}
\medskip



\noindent
{\em Proof.} It is enough to show:
\\
 (1) Map $\Phi$ is  
 a $k$-vector space isomorphism 
$\Phi: k\bBd(m,n) \rightarrow k\Jd_0^1(m+1,n+1)$.
\\
(2)
For any  $(p,s)\in \bB(m,n)$ and $(p',s')\in \bB(n,q)$
%
then
$
\Phi_{}(p,s)*\Phi_{}(p',s')=\Phi_{}((p,s)\bx (p',s')) .
$



\medskip 
For (1) consider first  a case with  $m<n$.
By construction there exists $x$ such that
$n-m=2x$. Then $p_0 = u^{\otimes x}$ in  \paref{otempora} gives an
embedding into the rank $n$ `algebra' case.
Now suppose  (for a contradiction) that $\Phi$ is not injective.
%
By \paref{otempora} this would induce non-injectivity in the algebra
case, contradicting \paref{le:PhiTheta}.

For surjectivity we may proceed as follows. Suppose $X$ lies in
$k J^1_0(m-1,m+1)$. Then $X\otimes u \in J^1_0(m+1,m+1) = J^1_{0,m}$
and so by Lemma~\ref{le:PhiTheta} there exists $Y \in b_m$ such that
$\Phi(Y)= X\otimes u$.
From Lemma~\ref{lem:Psi}(V)
and the construction
we see that $Y=Y'\otimes u$ and $\Phi(Y') = X$.
Thus $\Phi$ is surjective in this case.
Other cases are similar.

For (2)
let $t=\max(m,n,q)$.
For any $p_0\in J_{-1}(t-m,t-n)$ and 
$p_0'\in J_{-1}(t-n,t-q)$
note that
$(p\otimes p_0,s),\; (p'\otimes p_0',s')  \; \in \; J^\bullet(t,t)$.
Since $\Phi=\Theta$ is an algebra
isomorphism we have
\beq
\Phi((p\otimes p_0,s)\circ (p'\otimes p_0',s'))
=\Phi(p\otimes p_0,s)*\Phi(p'\otimes p'_0,s')
\label{isomt}\ .
\eeq
Let us expand the left-hand side first.
By construction 
\[
(p\otimes p_0,s)\circ (p'\otimes p_0',s')
=((p,s)\circ (p',s'))\otimes (p_0*p_0')
\]
hence, using Lemma \ref{otempora}, we have
\beq\label{equate1}
\Phi((p\otimes p_0,s)\circ (p'\otimes p_0',s'))=\Phi((p,s)\circ (p',s'))\otimes (p_0*p_0')\ .\eeq
Now, the r.h.s. of (\ref{isomt}) using Lemma \ref{otempora} for both terms reads
\beq\label{equate2}
\Phi(p\otimes p_0,s)*\Phi(p'\otimes p'_0,s')=
(\Phi(p,s)\otimes p_0)*(\Phi(p',s')\otimes p_0')=(\Phi(p,s)*\Phi(p',s'))\otimes (p_0*p'_0)\ ,
\eeq
where the right equality holds by construction.
The 
statement (2) of the lemma follows from the equality of the r.h.s. of (\ref{equate1}) 
and (\ref{equate2}) by Lemma~\ref{lem:moncancel}(II).
\qed


\[ \]

\section{On representation theory consequences for short Brauer algebras} \label{ss:other1}
\newcommand{\Specht}[2]{{\mathcal S}^{#1}_{#2}} 
\newcommand{\Deltam}[2]{\Delta^{#1}_{#2}} 
\newcommand{\Deltamd}[2]{D^{#1}_{#2}}   
\newcommand{\Cheb}{Chebyshev}


\subsection{Summary of relevant results for $\bb_n$}

Let us restrict to the case $k=\C$. 
From a representation theory perspective the natural parameterisation
of $\bb_n$ is $\delta=[2]$
(recall $[n]=(q^n - q^{-n})/(q-q^{-1})$)
and $\delta' = \frac{[m+1]}{[m]}$. 
Then if $m \not\in \Z$ we know that $\bb_n$ is semisimple, with a
well-known structure \cite{\ms}. 
If $m \in \Z $ but $q$ not a root of unity then the algebras are no
longer semisimple (for sufficiently large $n$), but the 
structure is still relatively simple to describe.
The most interesting case is $m \in \Z$ and $q$ a root of unity.
The structure in this case is quite complicated. 
See e.g. \cite{CoxGrahamMartin03} for a full description. 
With this summary in mind, note that due to (\ref{th:Theta})
we are interested in the cases when 
\[
\frac{[m+1]}{[m]} \; = \;  \frac{[2] + 1}{2}
\]
This is solved for example by $m=1$ when $[2]=1$.


For our present purposes the key point here  
has a precursor already
even from the arithmetically simpler Temperley--Lieb case, as follows.

\mdef \label{de:cheby}
Recall (see e.g. \cite{Martin91}) 
that the \Cheb\ polynomials are the polynomials $d_n$ determined by
the recurrence 
$ d_{n+2}=x d_{n+1}-d_n$,
with
initial conditions $d_0=1, d_1=x$.
(We write $x$ for $\delta$ here, simply for reasons of familiarity.)
The first few are
$d_n = 1,x , x^2 -1, x^3 -2x, x^4 -3x^2 +1, ...$
($n=0,1,2,3,...$).

These arise for example as the determinants of gram matrices such as:
\[
\Deltam{n}{n-2} = 
\mat{cccccccc}
\delta & 1 & 0 & 0 & 0\\
1 &  \delta & 1 & 0 &0\\
0 & 1 & \delta &  1 &0\\
0 & 0 & 1 & \delta & 1 \\
0 & 0 & 0 & 1    & \delta
\tam
\]
The obvious translational symmetry of this structure
(arising from the
local geometrical translational symmetry - the monoidal structure - of the TL
diagram `particles')
gives rise to the natural fourier parameterisation $d_{n-1} = [n]$. 
Loosely speaking, the geometrical boundary conditions here pick out a pure fourier sine
series
(fixing one end); and then the $n$ value (fixing the other end
--- hence the special behaviour at roots of unity of $q$).
The blob algebra generalises this essentially by changing the boundary
conditions.
Next we look for evidence of similar phenomena in the short Brauer
case. 


\ignore{{

\subsection{Summary of results for $\BBl{l,n}$}
...
\subsection{Standard Bratteli diagrams}
The difference between $\BBl{l,n}$ and $\BBBl{l,n}$ is that
...

}}

\subsection{Gram matrices, towers of recollement}

\begin{figure}
  \includegraphics[width=2.87643in]{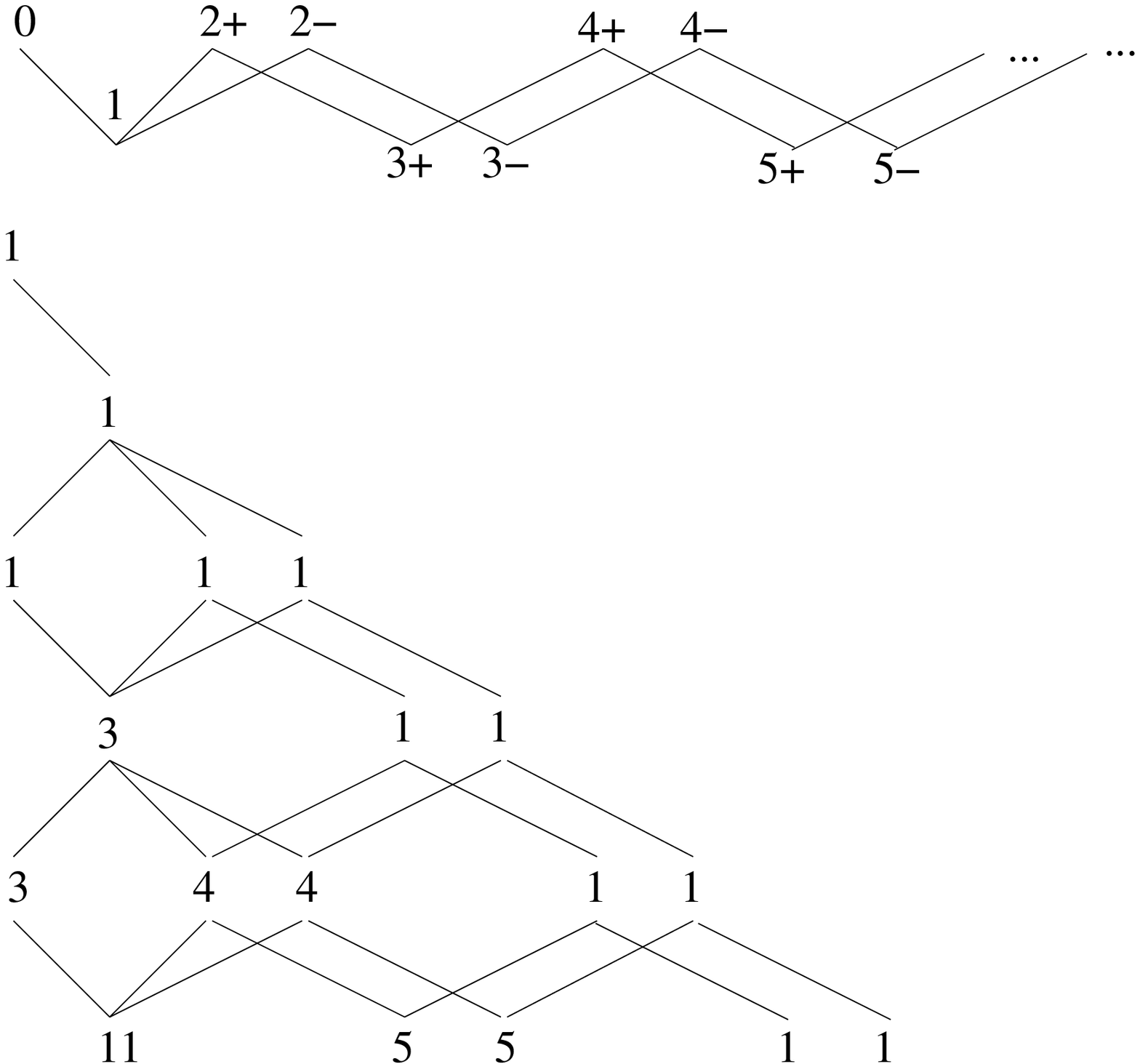}
  \caption{(a) Indicative labelling scheme for standard modules for
    height $l=0$ Brauer algebras.
    \quad
    (b) Bratteli diagram with dimensions of standard modules up to
    $n=5$.
    \label{fig:bratt1}}
\end{figure}

We assume familiarity with the representation theory as treated in
\cite{\kamy}, including the construction of standard modules. 

Here we restrict consideration to height 0.
Our labeling scheme for  Gram matrices
$\Deltam{n}{\lambda}$
of the
standard modules $\Specht{n}{\lambda}$
is 
$\Deltam{n}{\lambda} = \Deltam{n}{m,\pm}$
(superscript:  algebra rank $n$;  
subscript: number $m$ of propagating lines
and (for $m>1$) $\pm$ is the
symmetric / antisymmetric label from $S_2$).
See Fig.\ref{fig:bratt1}.
For example, 
the diagram basis for the $n=6$  
standard module
corresponding to $\lambda = (4,+)$
can be drawn as:
\[
\includegraphics[width=5.4in]{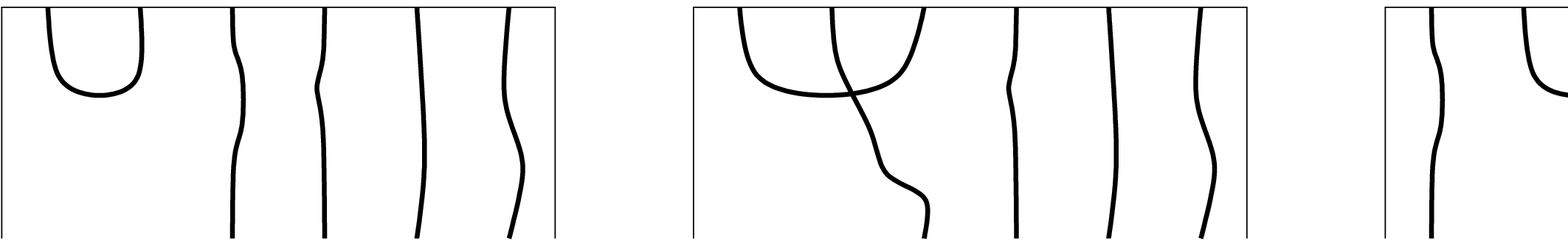}
\]
where we omit to draw the $(2)$-symmetrizer sitting on the first two
propagating lines
(thus we can draw the $\lambda=(4,-)$ case similarly, provided we keep
in mind the omission, which affects calculations). 
Note that the basis (so drawn) contains one extra diagram compared to the
$l=-1$/Temperley--Lieb case.

The extra diagram has an interesting effect on the gram matrix
of the natural contravariant form (see \cite{\kamy}).
As for
the TL case this can be computed in terms of \Cheb\ polynomials
(or equivalently fourier transforms). But here the initial conditions
are different. We have
\[
\Deltam{3}{1} = 
\mat{cccccccc}
\delta & 1 & 1 \\
1 & \delta & 1 \\
1 & 1 & \delta 
\tam,
\hspace{.351in}
\Deltam{4}{2,\pm} = 
\mat{cccccccc}
\delta & 1 & 1 & 0 \\
1 & \delta & 1 & \pm 1 \\
1 & 1 & \delta & 1 \\
0 & \pm 1 & 1 &  \delta 
\tam,
\hspace{.351in}
\Deltam{n}{n-2,+} = 
\mat{cccccccc}
\delta & 1 & 1 & 0 & 0 & 0 \\
1 & \delta & 1 & 1 & 0 & 0 \\
1 & 1 & \delta & 1 & 0 & 0 \\
0 & 1 & 1 &  \delta & 1 & 0 \\
0 & 0 & 0 & 1 & \delta &  1 \\
0 & 0 & 0 & 0 & 1 & \delta
\tam
\]
(we give the $n=6$ example, but the general pattern will be clear).
Laplace explanding $\Deltamd{n}{\lambda} = |\Deltam{n}{\lambda}|$
with respect to the bottom row we get a \Cheb\ recurrence
\[
\Deltamd{n}{n-2,\pm} = \delta \Deltamd{n-1}{n-3,\pm} - \Deltamd{n-2}{n-4,\pm}
\]
where the initial conditions are
\newcommand{\x}{\delta}
$\Deltamd{3}{1} = (\x-1)^2 (\x+2) $ and
$\Deltamd{4}{2,+}=  \x(\x-1)(\x^2+\x-4)$
and
$\Deltamd{4}{2,-}=  (\x-1)(\x+1)(\x-2)(\x+2)$.



Note from Theorem~1.1(ii) of \cite{CoxMartinParkerXi06}
(the tower-of-recollement method)
and Proposition~5.3 of \cite{\kamy}
(standard restriction rules)
that the other gram determinants and
indeed the `reductive' representation theory can be determined from
this subset of gram determinants.
We will address this task in a separate paper.
Here we restrict to some of the
key
preliminary observations. 

\medskip


The \Cheb\ polynomials $d_n$ from (\ref{de:cheby})
are a basis for the space of polynomials; and the
recurrence is linear, so we can express our recurrence in terms of
them,
and hence make use of their more `fourier-like' formulations:
$d_{n-1} = [n] = \frac{q^n -q^{-n}}{q-q^{-1}}$,
where $\delta = x = q+q^{-1}$.
The determinants $D_n^\pm$ of the key subset of Gram matrices 
of form $\Deltam{n}{n-2,\pm}$
%
can be expressed
as
%
\begin{eqnarray}
  D_n^+&=&(x-1) \left[ (x+2) (x-1) d_{n-3}-2x d_{n-4}\right]\\
D_n^-&=&(x-1) (x+2) \left[ (x-1) d_{n-3} -2 d_{n-4} \right]
\end{eqnarray}

Explicitly, the low rank cases of  all the Gram matrices are as follows:
\begin{eqnarray*}
D_1^3&=&(x-1)^2(x+2)\\
D_0^4&=&(x-1)^2x^3(x+2)\\
D_2^{4+}&=&(x-1)x(x^2+x-4)\\
D_2^{4-}&=&(x-1)(x+1)(x-2)(x+2)\\
D^5_1&=&(x-1)^{12}(x+1)(x-2)(x+2) {{}^6  (x^2 +x-4)}\\
D^{5+}_3&=&(x-1)(x^4+x^3-5x^2-x+2)\\
D^{5-}_3&=&(x-1)(x+2)(x^3-x^2-3x+1)\\
D^6_0&=&(x-1)^{12}x^{11}(x+1)(x-2)(x+2)^6(x^2+x-4)\\
D^{6+}_2&=&(x-1)^8x^5(x+1)(x-2)(x+2)(x^2+x-4)^6(x^4+x^3-5x^2-x+2)\\
D^{6-}_2&=&(x-1)^8(x+1)^6(x-2)^6(x+2)^7(x^2+x-4)(x^3-x^2-3x+1)\\
D^{6+}_4&=&(x-1)^2 x(x^3+2x^2-4x-6)\\
D^{6-}_4&=&(x-1)^2(x+2)(x^3-4x-2)
\end{eqnarray*}
%
(the cases not computed by recursion may be computed by brute force,
see below).


A key point to take from this is that the short Brauer algebras
manifest some similarities with the root-of-unity paradigm
for non-semisimplicity,
but move beyond it.
As noted, taken in combination with tower of recollement methods these results
`seed' the reductive representation theory
(the determination of decomposition matrices).
We address this analysis fully in a separate paper,
but the programme may be illustrated as follows.


\newcommand{\ing}[1]{\includegraphics[width=.31in]{xfig/#1.eps}}

\medskip


This form corresponds to the map from the standard module
$\Specht{n}{\lambda}$ to its
contravariant dual which, on general grounds,
maps the simple head to the socle \cite{\kamy}.
Thus when the form is non-singular we deduce that the standard module
is simple.
And on the other hand when it is singular the standard
module will have a corresponding submodule.
It is not generally easy to determine the rank of the form and hence
the dimension of the simple head from the gram determinant.
For example the rank of $\Deltam{5}{1}$ is easily seen to be 1,
while the
dimension of $\Specht{5}{1}$ is 11 (see Fig.\ref{fig:bratt1} or below)
and the determinant factor is $(x-1)^{12}$. 

To illustrate 
first consider  
$D^3_1$.
The basis here is
$\{ \ing{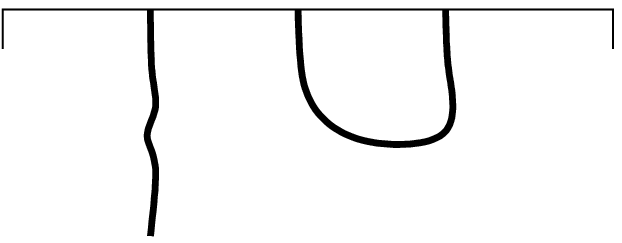} , \ing{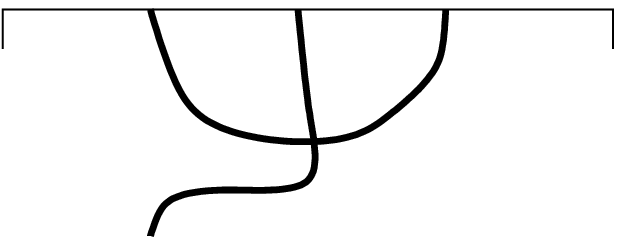} , \ing{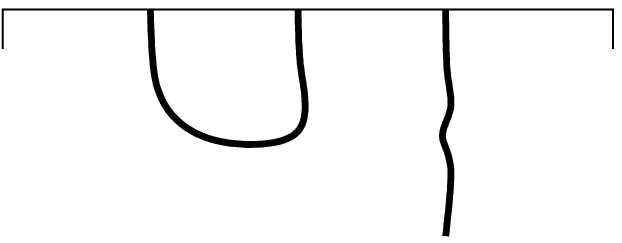} \}$.
For example 
the action of generators on the element 
$
\ing{S31bb} - \ing{S31a}
$
at the singular point $\delta = x =1$
is: 
\[
\stackrel{\ing{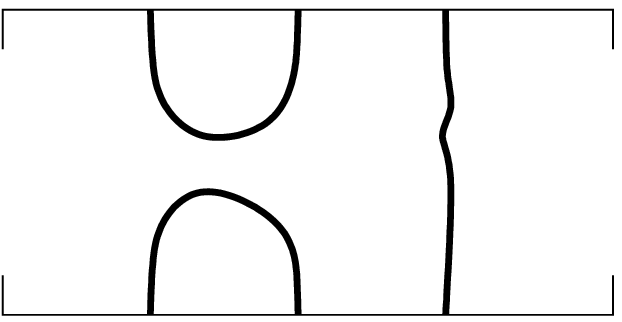} }{ \ing{S31bb}}  - \stackrel{\ing{S31d} }{ \ing{S31a}} 
\quad = \;  0,
\mbox{  }\quad\quad
\stackrel{\ing{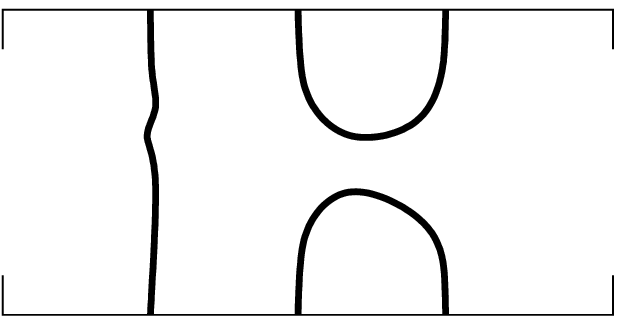} }{ \ing{S31bb}}  - \stackrel{\ing{S31e} }{ \ing{S31a}} 
\quad = \;  0,
\mbox{ and }\quad
\stackrel{\ing{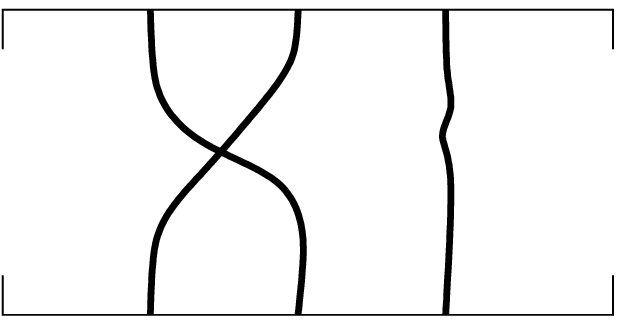} }{ \ing{S31bb}}  - \stackrel{\ing{S31ff} }{ \ing{S31a}} 
\quad = \;  - (\ing{S31bb} - \ing{S31a} )
\]
That is, when $\delta=x=1$ this element spans a submodule isomorphic to
$\Specht{3}{3,-}$.
Meanwhile for the element
$\ing{S31a} + \ing{S31bb} -2 \ing{S31c}$:
\[
  \stackrel{\ing{S31d} }{ \ing{S31a}}
+ \stackrel{\ing{S31d} }{ \ing{S31bb}} 
-2\stackrel{\ing{S31d} }{ \ing{S31c}} 
\quad =\quad
  \stackrel{\ing{S31e} }{ \ing{S31a}}
+ \stackrel{\ing{S31e} }{ \ing{S31bb}} 
-2\stackrel{\ing{S31e} }{ \ing{S31c}} 
\quad = \;  0,
\] \[
  \stackrel{\ing{S31ff} }{ \ing{S31a}}
+ \stackrel{\ing{S31ff} }{ \ing{S31bb}} 
-2\stackrel{\ing{S31ff} }{ \ing{S31c}} 
\quad = \;  \ing{S31a} + \ing{S31bb} -2 \ing{S31c}
\]
So this element spans a submodule isomorphic to $\Specht{3}{3,+}$. 
We deduce that the simple head is one-dimensional.

On the other hand consider
$\ing{S31a} + \ing{S31bb} + \ing{S31c}$
in case $\delta=x=-2$. This spans a submodule isomorphic to
$\Specht{3}{3,+}$.
Here the simple head is two-dimensional.

By the module-category embedding property \cite[(4.26)]{\kamy}
these standard module morphisms
have images in higher ranks,
thus when $x=1$ our map
$\Specht{3}{3,-} \rightarrow \Specht{3}{1}$
gives a map
$\Specht{5}{3,-} \rightarrow \Specht{5}{1}$
and so on.
The embedding functor is not exact so 
we cannot tell {\em directly} from the gram matrix if
an image map has a kernel. 
So (comparing also with the dimensions
from Fig.\ref{fig:bratt1}),
a naive lower bound on the exponent
in the factor $(x-1)^{12}$ in $D^5_1$
is 4+4,
corresponding to the dimensions of the simple heads of
$\Specht{5}{3,+}$ and $\Specht{5}{3,-}$
when $x=1$.
It is intriguing to compare with the blob case \cite{CoxGrahamMartin03}.
There the embedded
standard module morphisms are injective, but if that is the case here
the naive bound is still only lifted to 5+5,
so we see that there will be some nice subtleties here.

\begin{figure}
\[
\includegraphics[width=5.94in]{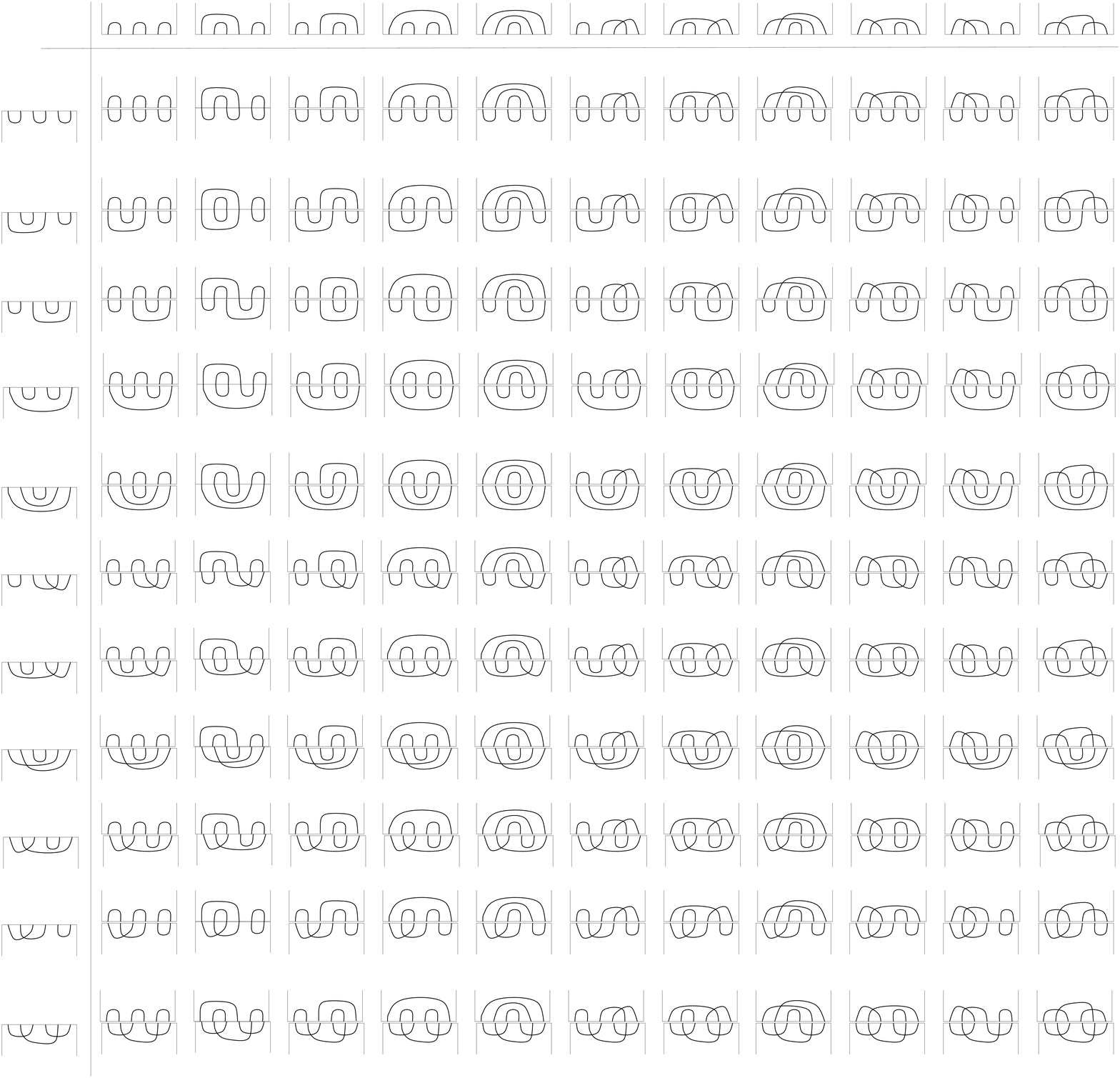}
\]
\caption{Gram matrix calculation for $n=6$ and $\lambda=0$. \label{fig:arrayz}}
\end{figure}

As a further illustration, the basis for $n=6$ and $\lambda=0$ is:
\[
\includegraphics[width=5.94in]{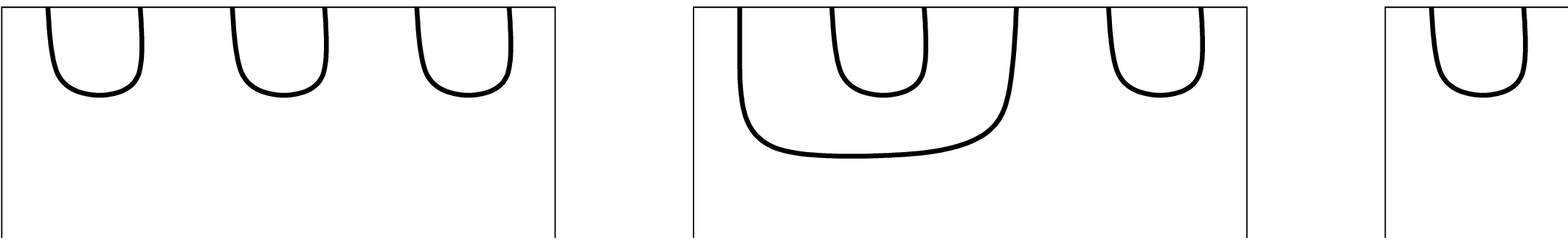}
\]
(N.B. the basis for $n=5$, $\lambda=1$ is combinatorially identical).
(As noted, we do not strictly need such cases for the `Cox criterion'.
It is enough to use $\lambda=n-2$.
We include it for 
curiosity's sake.)
The corresponding gram matrix then comes from the
array in Fig.\ref{fig:arrayz}.
Thus, writing $j$ for $\delta^j$
(with $j$ the number of connected components in a diagram),
the gram matrix is given by
\[
\mat{ccccc|ccccc|cc}
3 & 2 & 2 & 1 & 2 & 2 & 1 & 2 & 1 & 2 & 1
\\
2 & 3 & 1 & 2 & 1 & 1 & 2 & 1 & 1 & 2 & 1
\\
2 & 1 & 3 & 2 & 1 & 2 & 1 & 1 & 2 & 1 & 1
\\
1 & 2 & 2 & 3 & 2 & 1 & 2 & 1 & 2 & 1 & 1
\\
2 & 1 & 1 & 2 & 3 & 1 & 1 & 2 & 1 & 1 & 2
\\ \hline
2 & 1 & 2 & 1 & 1 & 3 & 2 & 1 & 1 & 1 & 2
\\
1 & 2 & 1 & 2 & 1 & 2 & 3 & 2 & 1 & 1 & 2
\\
2 & 1 & 1 & 1 & 2 & 1 & 2 & 3 & 2 & 1 & 1
\\
1 & 1 & 2 & 2 & 1 & 1 & 1 & 2 & 3 & 2 & 2
\\
2 & 2 & 1 & 1 & 1 & 1 & 1 & 1 & 2 & 3 & 2
\\ \hline
1 & 1 & 1 & 1 & 2 & 2 & 2 & 1 & 2 & 2 & 3
\tam
\]
The determinant here 
can still be computed by brute force.


\section{Discussion} \label{ss:discusstar}
Some notable open questions follow.

\noindent
Q1. How to generalise the `short Brauer' construction to the BMW algebra
\cite{BirmanWenzl89,Murakami87}?


\noindent
Q2.
How to 
relate the
usual two-parameter version of the blob algebra to the
short 
Brauer algebras --- which by the original construction have only a
single parameter.

Recall that there is, essentially trivially, a two-parameter version of $T_n$.
First recall that $T_n$ has a basis of non-crossing Brauer diagrams 
\cite{Weyl46o,brown} up to ambient isotopy
(see \S\ref{ss:pre1} for a summary of Brauer diagram concepts
--- ambient isotopy does not include, for example, 
the Reidemeister moves included
in general Brauer diagram equivalence, but it is sufficient in the
non-crossing case, and this is key here). 
The elements of the basis can be
seen as partitioning the interval into alcoves. 
These alcoves can be
shaded black or white with the property that 
\\
(A1) the colour changes across each boundary; and 
\\
(A2) the leftmost alcove is white, say.
\\
(NB Another way of saying this is that arcs have a well-defined `height'
in the sense of this paper, which is either odd or even.)
\\
Thus in composition both black and white loops may form.
The number of each separately is an invariant of ambient isotopy.
It follows that we may associate a different parameter to each.

Thus we have an algebra $T_n(\delta_b, \delta_w)$, say. It is easy to
see that 
$T_n(\delta_b, \delta_w) \cong T_n(\alpha\delta_b, \delta_w /\alpha)$
for any unit $\alpha$, so the difference can usually be scaled away.
For example recall the following.
%

{\theo{\label{th:TLgen1} {\rm \cite{Martin91}}
Consider the algebra defined by generators 
$U = \{ U_1,U_2,...,U_{n-1} \}$ and relations 
$\tau = \{ \mbox{
$U_i^2 = \delta U_i$,
$U_i U_{i\pm 1} U_i = U_i$, 
$U_i U_j = U_j U_i$, $j\neq i\pm 1$ } \}$.
The map 
$$
U_i \mapsto u_i = 
\{\{1,1' \} \{2,2' \}, ... , 
\{i , i+1 \},\{ i',i+1' \}, ...,
  \{n,n' \}\}    \qquad (i=1,2,...,n-1)
$$
extends to an algeba isomorphism 
$ k\langle U \rangle/\tau \cong T_n$.
\qed
}}

\medskip
To see the isomorphic two-parameter version consider the 
effect on the relations of the map 
$U_i \mapsto \alpha U_i$ ($i$ odd), 
$U_i \mapsto \alpha^{-1} U_i$ ($i$ even). 

\medskip

The blob algebra $\bb_n$ can be seen as the subalgebra of 
$T_{2n}(\delta_b, \delta_w)$ generated by diagrams with a lateral-flip
symmetry. In this subalgebra, however, it is {\em not } possible to
scale away the second parameter.

The short Brauer algebras are, from one perspective, generalisations
of $T_n$. It is interesting to consider if there are analogous
generalisations of the two-parameter version that (like the blob)
have the property that the second parameter becomes material. 
This is not obvious. The generalisation 
destroys the
two-tone alcove construction.

\medskip

How does the two-tone construction look in the categorical setting?
Here we write $T(n,m)$ for the subset $J_{-1}(n,m)$ of $J(n,m)$ of non-crossing
pair-partitions. We fix $\delta \in k$ and 
note that $\TTT = (\N_0, k T(n,m),*)$ is a subcategory of $\BB$. 
Indeed $\TTT = \BB^{-1}$. 
The inclusion is of $k$-linear categories, and also of monoidal
$k$-linear categories. 

As in the algebra case we note that in the non-crossing setting 
we can count the number of black and white loops separately
(i.e. these numbers are separately well-defined). 
Note however that the monoidal structure on $\TTT$ does not preserve
this property. It is the axiom (A2) that is the problem.




\vspace{1cm}

\noindent {\bf Acknowledgements.} We thank EPSRC for funding under
grant 
EP/I038683/1.
We thank Shona Yu, Azat Gaynutdinov and Peter Finch for useful
conversations.

\medskip

\appendix  \section*{Appendix} 
\section{Colour pictures for Lemma~\ref{mj}} \label{ss:colour}

\begin{figure} \newlength{\ggg} \setlength{\ggg}{10cm}
\begin{center}
\includegraphics[width=\ggg]{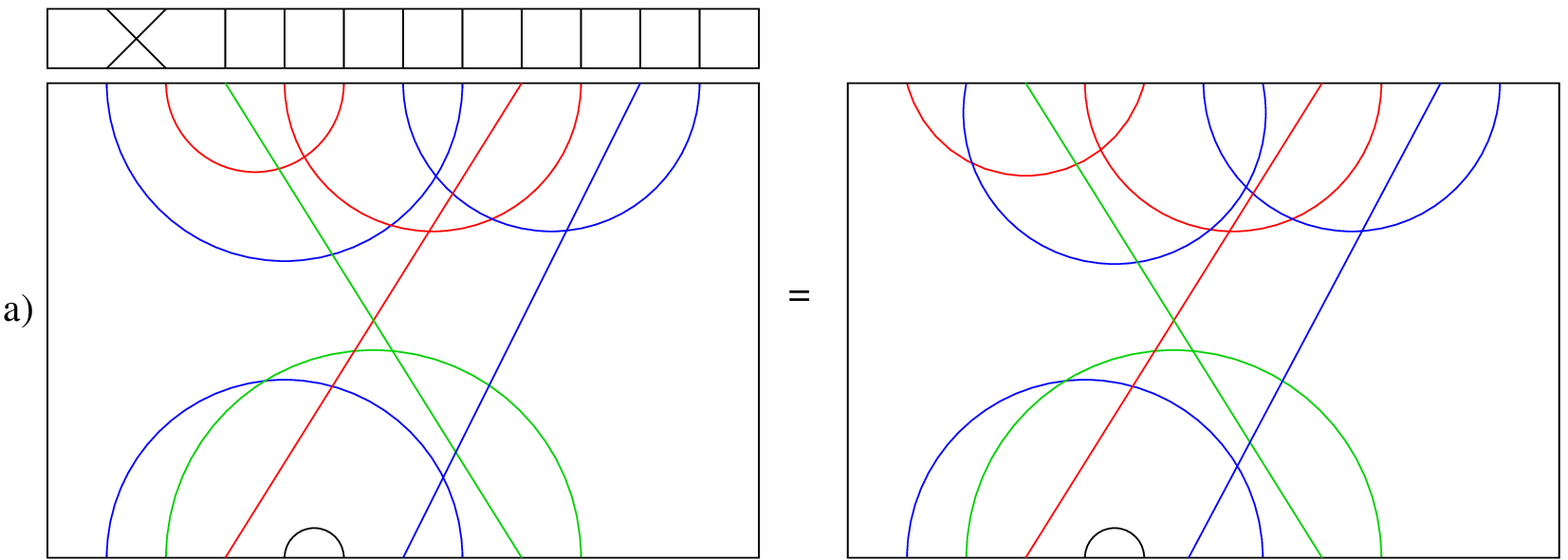}
\vspace{0.5cm}

\includegraphics[width=\ggg]{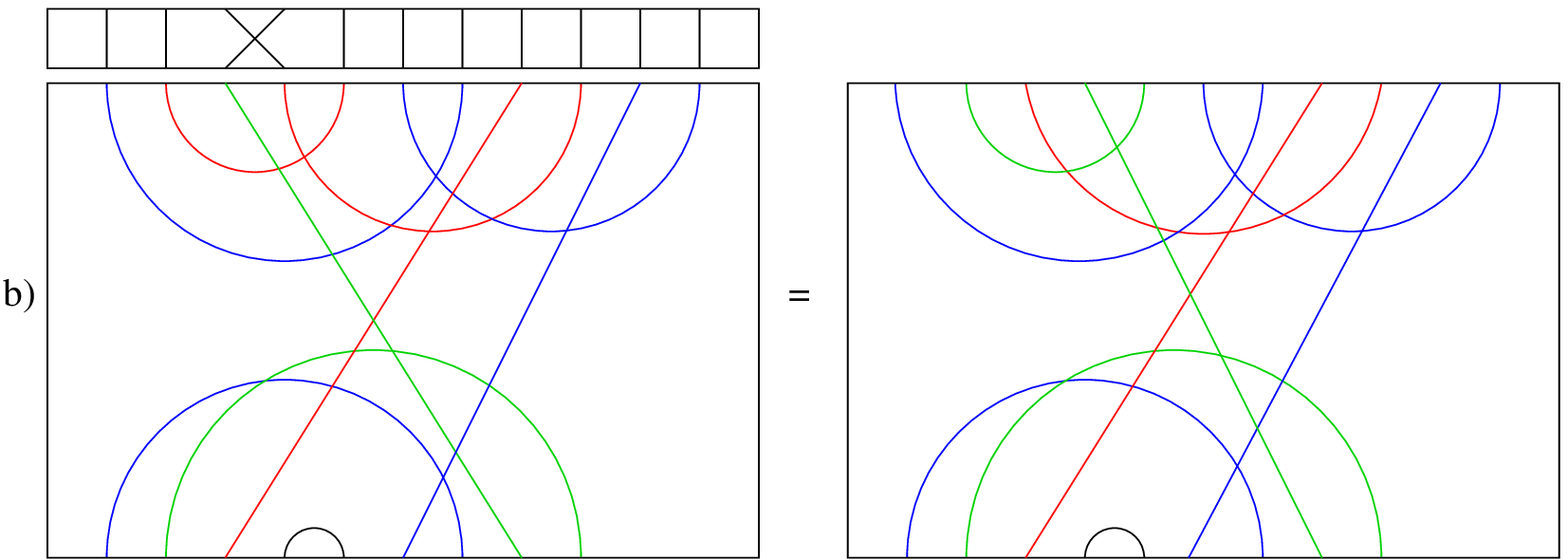}
\vspace{0.5cm}

\includegraphics[width=\ggg]{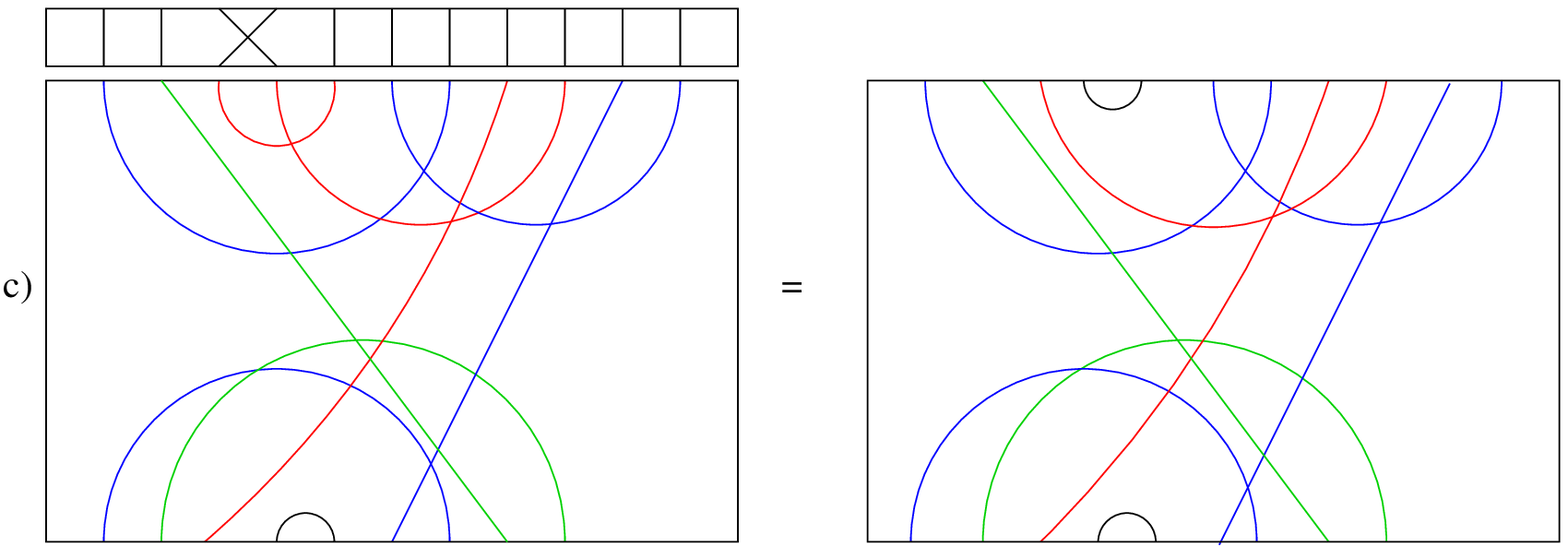}
\vspace{0.5cm}

\includegraphics[width=\ggg]{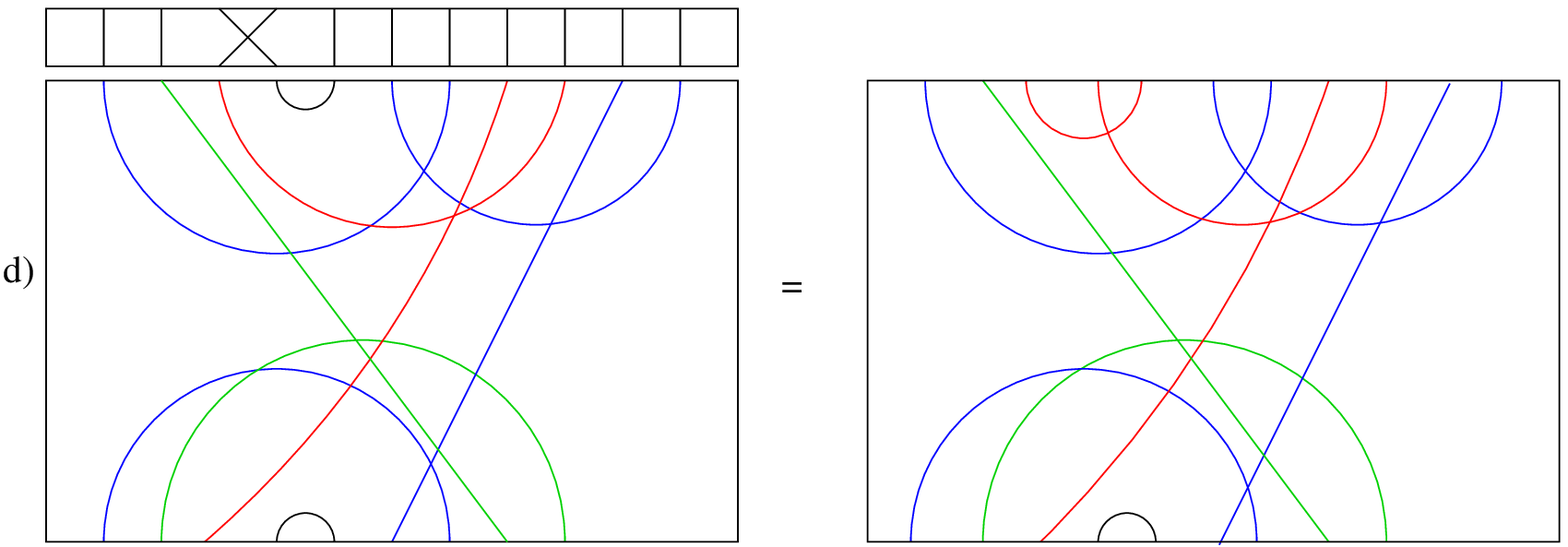}
\end{center}\caption{\label{smod}Examples of left actions of $S_4$ elements ((a): $\sigma_1$, b),c),d): $\sigma_3$) 
on elements of $J_2^3(11,9)$. They correspond to the four prototype cases discussed in the proof of Lemma (\ref{mj}).}
\end{figure}

{
Consider Lemma~\ref{mj}. 
When $j<i$, the initial points of chain from $i$ and $i+1$ are
interchanged (see Fig.~\ref{smod}a)). When $j=i$ we have  
three different cases: (i) The line from $i+1$ is part of a chain from
$[1,i]$, distinct from that from $i$, (ii) both  
lines from $i$ and $i+1$ belong to the same chain, (iii) the line from
$(i+1)$ is non intersecting.
Observe from  
Fig.~\ref{smod}b),c),d) that the resulting partitions are $Li-$simple
with $i$ exclusive chains from $[1,i]$ to $[1',i']$ 
and standalone pairs with no intersecting region with any other pair
of chains.
}





\bibliographystyle{amsplain}
\bibliography{\bibd/new31,bib/local}

\end{document}